\documentclass[12pt]{article}
\usepackage[utf8]{inputenc}
\title{\textcolor{black}{\textbf{{Normal forms, Lyapunov exponents, and pluripotential theory on {$\mathbb{P}^k(\cmplex)$}}}}}
\author{Virgile Tapiero}
\date{Friday March 13, 2026}

\usepackage[autolanguage]{numprint}
\usepackage[utf8]{inputenc}
\usepackage[T1]{fontenc}
\usepackage[english, french]{babel}
\usepackage{fontenc}
\usepackage{graphicx}
\usepackage{hyperref}
\usepackage{mathrsfs}
\usepackage{amsmath, amssymb, amscd, amsthm}
\usepackage{mathrsfs}
\usepackage{amsthm}
\usepackage[all]{xy} 
\usepackage{esint}
\usepackage{color}
\usepackage{pst-all}
\usepackage{upgreek}
\usepackage{comment}

\newtheorem{thm}{Theorem}[section]
\newtheorem{lemme}[thm]{Lemma}
\newtheorem{prop}[thm]{Proposition}

\newtheorem{defn}[thm]{Definition}
\newtheorem{rmq}[thm]{Remark}

\newtheorem*{theoremA*}{Theorem A}
\newtheorem*{theoremA*k}{Theorem A on $\mathbb{P}^k$}
\newtheorem*{theoremB*}{Theorem B}
\newtheorem*{theoremB*k}{Theorem B on $\mathbb{P}^k$}

\newcommand{\reels}{\mathbb{R}}
\newcommand{\cmplex}{\mathbb{C}}

\newcommand{\LLog}{\ \mathrm{Log}}

\newcommand{\Dciteequa}[1]{{Equation (#1)}}

\newcommand{\Dcitethm}[1]{{Theorem #1}}
\newcommand{\Dcitelem}[1]{{Lemma #1}}
\newcommand{\Dciteprop}[1]{{Proposition #1}}

\providecommand{\keywords}[1]
{
  {\small \textbf{\textit{Keywords---}}} #1
}

\providecommand{\remerciements}[1]
{
  {\small \textbf{\textit{Acknowledgements---}}} #1
}

\usepackage[left=3cm,right=3cm,top=2cm,bottom=3cm]{geometry}

\newenvironment{changemargin}[2]{\begin{list}{}{%
\setlength{\topsep}{0pt}%
\setlength{\leftmargin}{0pt}%
\setlength{\rightmargin}{0pt}%
\setlength{\listparindent}{\parindent}%
\setlength{\itemindent}{\parindent}%
\setlength{\parsep}{0pt plus 1pt}%
\addtolength{\leftmargin}{#1}%
\addtolength{\rightmargin}{#2}%
}\item }{\end{list}}

\begin{document}
\hypersetup{pdfborder=0 0 0}

\maketitle

\begin{center}
    \begin{changemargin}{0.000005cm}{0.000005cm}
        \begin{otherlanguage}{english}
            \begin{abstract}
            {
            We study the dynamical properties of endomorphisms $f$ of $\mathbb{P}^k$ of algebraic degree $d \geq 2$. We investigate the relationships between the Green current $T$ of $f$, the equilibrium measure $\mu = T^k$, and the Lyapunov exponents $\lambda_1 \geq \cdots \geq \lambda_k$ of $\mu$. The latter are bounded below by $\frac{1}{2} \mathrm{Log} \ d$. Dujardin proved in \cite{Duj12} that if $\mu \ll T^r \wedge \omega_{\mathbb{P}^k}^{k-r}$ for some $1 \leq r \leq k-1$, then $\lambda_{r+1} = \cdots = \lambda_k = \frac{1}{2} \mathrm{Log} \ d$. In this article we prove that, conversely, if $\lambda_r>\lambda_{r+1} = \cdots = \lambda_k = \frac{1}{2} \mathrm{Log} \ d$, then $\mu\ll T^r\wedge\omega_{\mathbb{P}^k}^{k-r}$, answering a question asked by Dujardin. Our arguments rely on pluripotential theory, ergodic theory, and normal forms for the inverse branches of the endomorphism. We also use normal forms to provide another proof of Dujardin’s result.
            } 
            \end{abstract}
        \end{otherlanguage}
    \end{changemargin}
\end{center}

{\section{Introduction}\label{sec:introsliceproperty}}

\subsection{Overview} 

\subsubsection{Classical results on Lyapunov exponents of rational maps}

Let $f$ be a rational map of $\mathbb{P}^1$ of degree $d \geq 2$, and let $\omega_{\mathbb{P}^1}$ be the spherical $(1,1)$-form on $\mathbb{P}^1$. According to Brolin, Lyubich, and Freire-Lopes-Ma\~né {\cite{Bro65,Lyu83,FLM83}}, the probability measure $\mu = \lim_n \frac{1}{d^n}{f^n}^{*}\omega_{\mathbb{P}^1}$ is the unique invariant measure of maximal entropy for $f$, equal to $h_{\mathrm{top}}(f)=\mathrm{Log}\ d$. 
This measure has a constant Jacobian $f^*\mu = d\mu$, is mixing, and integrates functions with logarithmic singularities. Therefore, the measure $\mu$ has a Lyapunov exponent $\lambda_f = \int_{\mathbb{P}^1} \mathrm{Log}|f'(z)|\ \mathrm{d}\mu(z)$, where $|f'(z)|$ is computed with respect to the metric $\omega_{\mathbb{P}^1}$. By Margulis-Ruelle inequality, $\lambda_f$ satisfies the inequality
$$\lambda_f\geq\frac{1}{2}\mathrm{Log}\ d.$$
The map $f \mapsto \lambda_f$ is continuous (see Ma\~n\'e \cite{Mane88}) and plurisubharmonic on $\mathrm{Rat}_d$, the Zariski open subset of $\mathbb{P}^{2d+1}$ that parametrizes all rational maps of degree $d$ (see DeMarco, Przytycki \cite{Przyt85, DeMarco03}). The extremal value $\lambda_f = \frac{1}{2}\mathrm{Log}\ d$ is attained on $\mathrm{Rat}_d$, and the corresponding rational maps can be characterized in several different ways. The articles by Ledrappier {\cite{led84}} and Zdunik {\cite{zdu90}} (see also Mayer {{\cite{Mayer02}}}) indeed show that the following properties are equivalent\footnote[1]{They are also equivalent to $\mathrm{HD}(\mu)=2$, where $\mathrm{HD}(\mu)$ is the Hausdorff dimension of $\mu$.}
\footnote[2]{(a)$\Rightarrow$(b): Ledrappier \cite{led84}; (b)$\Rightarrow$(c): Zdunik \cite{zdu90} and Mayer {{\cite{Mayer02}}}.}:
\begin{equation}\label{eq:equivI}
    \mathrm{(a)}\ \lambda_f=\frac{1}{2}\mathrm{Log}\ d;\ \mathrm{(b)}\ \mu\ll\mathrm{Leb}_{\mathbb{P}^1}=\omega_{\mathbb{P}^1};\ \mathrm{(c)}\ f\ \mathrm{is\ Latt\text{è}s\ on}\ \mathbb{P}^1.
\end{equation}
Lattès maps were introduced by Samuel Lattès in 1918, we refer to Milnor's article \cite{Mil06}, which provides a modern and detailed exposition of Lattès rational maps. Roughly speaking, a Lattès map is a map on $\mathbb{P}^1$ that can be lifted to an affine map on a complex torus $\mathbb{C}/\Lambda$ via a finite branched cover $\sigma: \mathbb{C}/\Lambda \to \mathbb{P}^1$. These maps have played an increasing role in complex dynamics, particularly in studying the bifurcation locus of $\mathrm{Rat}_d$ and of the moduli space $\mathcal{M}_d=\mathrm{Rat}_d/\mathrm{PGL}(2,\cmplex)$, see for example McMullen \cite{McMullen87}, DeMarco \cite{DeMarco01,DeMarco03}, Berteloot \cite{BertelootCetraro13}, Buff-Gauthier \cite{BuffGauthier13}.

\subsubsection{Generalization for endomorphisms of \texorpdfstring{$\mathbb{P}^k$}{TEXT}}
 
This article aims to study and generalize the equivalences \eqref{eq:equivI} in higher dimensions, particularly we propose a generalized version of (a) $\Leftrightarrow$ (b) for endomorphisms of $\mathbb{P}^k$, with $k \geq 2$. For a holomorphic map $f$ of $\mathbb{P}^k$ with algebraic degree $d \geq 2$, the probability measure $\mu$ can be defined by $\mu := \lim_n \frac{1}{d^{kn}}{f^n}^{*}\omega_{\mathbb{P}^k}^{k}$, where $\omega_{\mathbb{P}^k}$ is the normalized Fubini-Study $(1,1)-$form of $\mathbb{P}^k$, we refer to the books {\cite{sib99,dinsib}} by Dinh and Sibony. $\mu$ is called the equilibrium measure of $f$, which is the unique invariant measure of maximal entropy for $f$, equal to $h_{\mathrm{top}}(f)=\mathrm{Log}\ d^k$ according to Gromov {\cite{Gro03}} and Misiurewicz-Przytycki \cite{MisiuPrzy77}. This measure has a constant Jacobian $f^*\mu = d^k \mu$, is mixing, and integrates functions with logarithmic singularities. In particular, $\mu$ has Lyapunov exponents $\lambda_1 \geq \cdots \geq \lambda_k$.

The measure $\mu$ is also equal to $T^k = T \wedge \cdots \wedge T$, where $T$ is the Green current of $f$ defined by $T := \lim_n\frac{1}{d^n}{f^n}^{*} \omega_{\mathbb{P}^k}$. It is a positive closed $(1,1)$-current of mass 1 satisfying $f^*T = dT$. Again, we refer to {\cite{sib99,dinsib}}. $T$ can be thought as an invariant singular metric on $\mathbb{P}^k$.

Briend-Duval {\cite{BriDuv99}} showed that the Lyapunov exponents of $\mu$ on $\mathbb{P}^k$ are bounded below by the same constant as in dimension $1$:
\begin{equation}\label{eq:Briend-Duval}
    \lambda_1\geq\cdots\geq\lambda_k\geq\frac{1}{2}\mathrm{Log}\ d.
\end{equation}
But this is no longer an application of Margulis-Ruelle inequality, the proof relies on more difficult arguments involving pluripotential theory. Another proof, using only local unstable manifolds and entropy, was later given by De Thélin {\cite{det08}}.

A remarkable fact that motivates the present article is that the equivalences \eqref{eq:equivI} hold in higher dimensions by involving the current $T$. More precisely, the works of Berteloot, Dupont, and Loeb {\cite{berloe98,BL01,berdup05,dup06}} show that the following assertions are equivalent (where $U$ is an nonempty open set)\footnote[2]{They are also equivalent to $\mathrm{HD}(\mu)=2k$, as shown by Dinh-Dupont \cite{DinhDup04}.}\footnote[3]{(i)$\Rightarrow$(ii): Dupont \cite{dup06}; (ii)$\Rightarrow$(iii): Berteloot--Dupont \cite{berdup05}; (iii)$\Rightarrow$(iv): Berteloot--Loeb \cite{BL01}. Dupont \cite{dup10} also gave another proof of (i)$\Rightarrow$(ii) using a Central Limit Theorem for the observable $\mathrm{Log\ det}\ df$.}\footnote[4]{In (iii), if $T|_U$ is smooth and $>0$ for some $U\ne\emptyset$, then $f$ is Lattès and thus $T$ is smooth and $>0$ on $\mathbb{P}^k\backslash E$, where $E$ is a codimension $1$ analytic subset of $\mathbb{P}^k$.}:
\begin{equation}\label{eq:equivII}
    \mathrm{(i)}\ \forall j,\ \lambda_j=\frac{1}{2}\mathrm{Log}\ d;\ \mathrm{(ii)}\ \mu\ll\omega_{\mathbb{P}^k}^k;\ \mathrm{(iii)}\ \exists U : T|_U\ \mathrm{is\ smooth} >0;\ \mathrm{(iv)}\ f\ \mathrm{is\ Latt\text{è}s}.
\end{equation}
In this article, we develop techniques that allow us to study the measure $\mu$ in the presence of inhomogeneity in the Lyapunov spectrum, i.e., $k\geq2$ and there exists $1 \leq r \leq k-1$ such that
\begin{equation*}
	\lambda_1\geq\cdots\geq\lambda_r>\lambda_{r+1}=\cdots=\lambda_k=:\lambda_{\mathrm{min}}.
\end{equation*}
If moreover $\lambda_{\mathrm{min}} = \frac{1}{2} \mathrm{Log}\ d$, we will say that $f$ is an $r$-extremal map of $\mathbb{P}^k$, in contrast to the case where all exponents are minimal, which corresponds to extremal mappings in the literature. The cases $k=2$ and $r=1$ correspond to semi-extremal mappings of $\mathbb{P}^2$. The terminology of semi-extremal maps was introduced in \cite{dupont_rogue_2020,DupRog20}. These mappings have been studied by Bianchi, Bonifant, Dabija, Dupont, Milnor, Rogue, Taflin and the author \cite{BoDaMil07,BiaTaf17,dupont_rogue_2020,DupRog20,DT,duptap23,tap24}.  

The proof of (i)$\Rightarrow$(ii) in the equivalences \eqref{eq:equivII} relies essentially on standard tools from differentiable ergodic theory, and the fact that $f$ is holomorphic is not central. In contrast, in this article, our techniques are based on the joint use of ergodic theory with two main tools: holomorphic normal forms for the inverse branches of the iterates of $f$, and pluripotential theory through the use of the current $T$.

\subsection{New results}

From now on we will always assume $k\geq2$. In his work on Fatou directions {\cite{Duj12}}, R. Dujardin showed the following property for each $r \in\{1, \cdots, k-1\}$:
\begin{equation}\label{eq:DujardinTheorem}
	\mu\ll T^r\wedge\omega_{\mathbb{P}^k}^{k-r}\Longrightarrow\lambda_{r+1}=\cdots=\lambda_k=\frac{1}{2}\mathrm{Log}\ d.
\end{equation} 
Actually, R. Dujardin proved this result for $k=2$, but his arguments also apply to $k \geq 2$. A different proof is provided in \cite{duptap23} for $\mathbb{P}^2$ using normal forms for the inverse branches of $f$. We reproduce this proof, adapted to $\mathbb{P}^k$, in Theorem \ref{thm:ProofofDujardinTheorem} (Section \ref{sec:normalforms}). 

Dujardin \cite[Question 3.7]{Duj12} asked about the converse of \eqref{eq:DujardinTheorem} for $k=2$. The first partial answer was given by Dupont-Taflin \cite{DT}, assuming an additional geometric hypothesis, namely that $f$ preserves a fibration. In particular, their result provides examples of $r$-extremal maps for each $1\leq r\leq k-1$. In \cite{tap24}, the author provided a converse result for $k=2$, showing that certain constraints on $\mu$ and $T$, satisfied by some Dupont-Taflin examples, lead to the existence of a local invariant foliation $\mathcal{F}$ for (an iterate of) $f$ on a neighborhood of $\mathrm{Supp}(T)\backslash\mathcal{E}_f$, where $\mathcal{E}_f$ is the exceptional set of $f$. When $\mathcal{F}$ extends to a neighborhood of $\mathrm{Supp}(T)$, $\mathcal{F}$ extends to $\mathbb{P}^2$ and it is an invariant pencil of lines. 

We now state our new results, which provide an answer to Dujardin's question. The first result concerns the normal forms mentioned earlier. Let $\widehat{\mathbb{P}}^k:=\{\widehat{x}=(x_n)_{n\in\mathbb{Z}}:f(x_n)=x_{n+1}\}$ be the space of orbits, $\pi_0:\widehat{x}\mapsto x_0$ the natural projection, $\widehat{f}$ the left shift, and $\widehat{\mu}$ the unique probability measure invariant under $\widehat{f}$ such that $(\pi_0)_*\widehat{\mu}=\mu$. Then the following diagram commutes for $\displaystyle{n\!\geq\! N(\widehat{x})}$ ($\mathrm{\displaystyle{cf.\ Theorem\ \ref{thm:normalforms}}}$):
\begin{equation}
\label{eq:diagrameformesnormales_intro}
    \xymatrix{
  B(x_{-n}, 2\eta(\widehat{x}_{-n})) \ar[d]_{\sigma_{\widehat{x}_{-n}}} & &  B(x_0, 2\eta(\widehat{x})) \ar[ll]_{f^{-n}_{\widehat{x}}} \ar[d]^{\sigma_{\widehat{x}}=(Z_{\widehat{x}},W_{\widehat{x}})=(Z^1_{\widehat{x}}, \cdots, Z^r_{\widehat{x}}, W^{r+1}_{\widehat{x}}, \cdots, W^k_{\widehat{x}})}  \\
  \mathbb{D}^k(\rho(\widehat{x}_{-n}))  & & \mathbb{D}^k(\rho(\widehat{x})) \ar[ll]^{R_{n,\widehat{x}}}
}\end{equation}
The holomorphic map $f^{-n}_{\widehat{x}}$ is the inverse branch of $f^n$ that maps $x_0$ to $x_{-n}$. The vertical maps $\sigma_{\widehat{x}_{-n}}$ and $\sigma_{\widehat{x}}$ are injective, holomorphic, and closed to the identity in $C^1$-norm. The map $R_{n,\widehat{x}}$ is an invertible polynomial map of $\cmplex^k$ such that the derivative $d_0R_{n,\widehat{x}}$ has eigenvalues of the form $e^{-n(\lambda_1\pm\varepsilon)},\cdots,e^{-n(\lambda_k\pm\varepsilon)}$.

A set $W^u(\widehat{x},R):=\{\widehat{y}=(y_n)_{n\in\mathbb{Z}}\ : \ y_{-n}=f^{-n}_{\widehat{x}}(y_0),\ y_0\in B(x_0,R)\}$, with $R\leq 2\eta(\widehat{x})$, is called an unstable manifold. A Pesin box $P:=\bigsqcup_{\widehat{x}\in\mathcal{T}}W^u(\widehat{x},R)$ is a disjoint collection of these manifolds (see Definition \ref{defn:pesinboxes}). Let us consider the column vector $dW_{\widehat{x}}:=\left[dW_{\widehat{x}}^{r+1}\cdots dW_{\widehat{x}}^k\right]^{\mathrm{T}}$ of $W$-normal forms given by \eqref{eq:diagrameformesnormales_intro}. Our first result shows that these normal forms glue together on Pesin boxes (See Theorem \ref{thm:localfoliations} for more details):

\begin{thm}\label{thm:Patching} 
    Assume that there exists $1\leq r\leq k-1$ such that ${\lambda_{r}>\lambda_{r+1}}$.
    Let $P=\bigsqcup_{\widehat{z}\in\mathcal{T}}W^u(\widehat{z},R)$ be a Pesin box of positive measure for $\widehat{\mu}$, $\displaystyle{\mathcal{T}\!\!\subset}\{\widehat{z}:\eta(\widehat{z})\geq2R\}$. Then for $\widehat{\mu}-$almost every $\widehat{x}\in P$, denoting $\widehat{z}\in\mathcal{T}$ such that $\widehat{x}\in W^u(\widehat{z},R)$, we have:
	\begin{enumerate}
        \item For each $r+1\leq l\leq k$, $W_{\widehat{x}}^l$ extends as an holomorphic submersion on $B(z_0,R)$.
		\item There exists $\displaystyle{C\!\in\!\mathrm{GL}_{k-r}(\mathcal{O}(B(z_0,R)))}$ such that $\displaystyle{dW_{\widehat{x}}\!=\!CdW_{\widehat{z}}}$ on $B(z_0,R)$. 
        \item If moreover $\lambda_{r+1}=\cdots=\lambda_k$, $C\in\mathrm{GL}_{k-r}(\cmplex)$ is a constant matrix.
	\end{enumerate}
\end{thm}

Theorem \ref{thm:Patching}, coupled with the use of the Green current $T$, allows us to prove our second result, which generalizes the classical Pesin formula in the case $\mu \ll \mathrm{Leb}_{\mathbb{P}^k}$:

\begin{thm}\label{thm:TheoremB} Assume that $\lambda_r > \lambda_{r+1}$ and that $\lambda_{r+1}\geq\cdots\geq\lambda_k$ are not resonant, with $1 \leq r \leq k-1$. Then, there exists a Pesin box $\mathcal{A}$ of $\widehat{\mu}$-positive measure, a measurable partition $\xi^u$ of $\widehat{\mathbb{P}}^k$, and a family $(q_{\widehat{x}})_{\widehat{x}\in\widehat{\mathbb{P}}^k}$ of probabilities $q_{\widehat{x}}\in\mathrm{Prob}(\xi^u_{\widehat{x}})$, such that the following holds : by denoting $\widehat{q}(A):=\int q_{\widehat{x}}(A\cap\xi^u_{\widehat{x}})\ \mathrm{d}\widehat{\mu}(\widehat{x})$, $\widehat{q}$ is a probability on $\widehat{\mathbb{P}}^k$, and for all $n \geq 0$ and for $\widehat{\mu}$-almost every $\widehat{x}$ we have:
\begin{enumerate}
    \item $(\pi_0)_{*}\widehat{q}\ll T^r\wedge \omega_{\mathbb{P}^k}^{k-r}$ and $(\pi_0)_*q_{\widehat{x}}\ll T^r\wedge\omega_{\mathbb{P}^k}^{k-r}$.
    \item If $\widehat{x}\in\mathcal{A}$, $q_{\widehat{x}} = \frac{1}{L(\widehat{x})}{(T^r\wedge dd^c|W_{\widehat{x}}|^2)\circ\pi_0\lfloor_{\xi_{\widehat{x}}^u}}$, $L(\widehat{x})=(T^r\wedge dd^c|W_{\widehat{x}}|^2)(\pi_0(\xi^u_{\widehat{x}}))>0$.
    \item If $\lambda_{r+1}=\cdots=\lambda_k$ then $q_{\widehat{x}}$ is the conditional measure of $\widehat{q}$ on $\xi^u_{\widehat{x}}$.
    \item $\int-\LLog\ q_{\widehat{x}}\left(\widehat{f}^{-n}\xi^u\right)_{\widehat{x}}\ \mathrm{d}\widehat{\mu}(\widehat{x}) = \LLog\ d^{nr} + 2n(\lambda_{r+1}+\cdots+\lambda_k)$.
    \item $\int-\LLog\ \mu_{\widehat{x}}\left(\widehat{f}^{-n}\xi^u\right)_{\widehat{x}}\ \mathrm{d}\widehat{\mu}(\widehat{x}) =  \LLog\ d^{nr} + 2n(k-r)\frac{1}{2}\mathrm{Log}\ d$.
    \item If $\lambda_{r+1}=\cdots=\lambda_k=\frac{1}{2}\mathrm{Log}\ d$, then $\widehat{q}=\widehat{\mu}$.
\end{enumerate}
\end{thm}

In this statement, we have denoted $dd^c |W_{\widehat{x}}|^2 := dd^c |W_{\widehat{x}}^{r+1}|^2 \wedge \cdots \wedge dd^c |W_{\widehat{x}}^k|^2$, $\xi^u_{\widehat{x}}$ is the atom of the partition $\xi^u$ containing $\widehat{x}$, and $\mu_{\widehat{x}}$ is the conditional measure of $\widehat{\mu}$ on $\xi^u_{\widehat{x}}$.

\newpage

From Theorem \ref{thm:TheoremB}, we finally deduce our third result, which answers R. Dujardin's question \cite[Question 3.7]{Duj12} :

\begin{thm}\label{thm:TheoremA} Let $f$ be an endomorphism of $\mathbb{P}^k$ of degree $d\ge2$ and let $1\le r\le k-1$. Then
\begin{equation}\label{eq:DujardinTheoremII}
\lambda_r>\lambda_{r+1}=\cdots=\lambda_k=\frac{1}{2}\mathrm{Log}\ d\Longrightarrow\mu\ll T^r\wedge\omega_{\mathbb{P}^k}^{k-r}.
\end{equation} 
For the case $r=1$ we have moreover a characterization
\begin{equation}\label{eq:equivalencek-1exposants}
\lambda_2=\cdots=\lambda_k=\frac{1}{2}\mathrm{Log}\ d\Longleftrightarrow\mu\ll T\wedge\omega_{\mathbb{P}^k}^{k-1}.
\end{equation} 
\end{thm}

Let us make two remarks on this result. First, note that the converse implication of \eqref{eq:equivalencek-1exposants} comes from Dujardin's theorem \eqref{eq:DujardinTheorem}, and that the direct implication of \eqref{eq:equivalencek-1exposants} is a particular case of \eqref{eq:DujardinTheoremII} if $\lambda_1 > \lambda_2$. In the case where $\lambda_1 = \lambda_2$, the direct implication of \eqref{eq:equivalencek-1exposants} follows from (i)$\Rightarrow$(ii) and (ii)$\Rightarrow$(iii) in \eqref{eq:equivII} (because in this case $\omega_{\mathbb{P}^k}^k \ll T^r \wedge \omega_{\mathbb{P}^k}^{k-r}$, see the footnote 4 on page 3).

Second, the converse of \eqref{eq:DujardinTheoremII} is false in general. For example, if $f$ is a Lattès map, we have $\mu \ll \omega_{\mathbb{P}^k}^k \ll T^r \wedge \omega_{\mathbb{P}^k}^{k-r}$ by \eqref{eq:equivII}, but $\lambda_r = \lambda_{r+1}$ since all the Lyapunov exponents are minimal, proving that the converse of \eqref{eq:DujardinTheoremII} does not hold in general. 

\subsection{Outlines of the proofs}

\subsubsection{Proof of Items \textit{1.} to \textit{5.} of Theorem \ref{thm:TheoremB}}

The construction of a measurable partition $\xi^u$ satisfying Item \textit{5.} is due to Ledrappier \cite{led84} (see also Dupont \cite{dup06}), the novelty here is the formula in Item \textit{4.} We now explain our construction of the measures $q_{\widehat{x}}$ in order to obtain $\widehat{q}$ satisfying the formula in Item \textit{4.} Our general strategy follows the approach of  \cite{led84,dup06}: we construct $\widehat{q}$ with the same Jacobian as $\mu$ (namely $d^k$) and absolutely continuous with respect to $T^r\wedge\omega_{\mathbb{P}^k}^{k-r}$. For this purpose, we introduce the measures 
\begin{equation}\label{eq:Twedgedd^c|W_x|^2}
    T^r\wedge dd^c|W_{\widehat{x}}|^2,\ \mathrm{with}\ dd^c|W_{\widehat{x}}|^2:=dd^c|W_{\widehat{x}}^{r+1}|^2\wedge\cdots\wedge dd^c|W_{\widehat{x}}^k|^2.
\end{equation}
We then define $q_{\widehat{x}}$ from the measure \eqref{eq:Twedgedd^c|W_x|^2} by observing the following fact: a crucial property of the commutative diagram \eqref{eq:diagrameformesnormales_intro} (see Theorem \ref{thm:normalforms}) is that when $\lambda_{r+1} \geq \cdots \geq \lambda_k$ are not resonant, the components $R_{n,\widehat{x}}^{l}$, $r+1 \leq l \leq k$, of $R_{n,\widehat{x}}$ are linear and satisfy
$$R^{l}_{n,\widehat{x}}(z,w)=\beta_{n,\widehat{x}}^lw\ ,\ e^{-n(\lambda_l+\varepsilon)}\le|\beta_{n,\widehat{x}}^l|\le e^{-n(\lambda_l-\varepsilon)}.$$ 
Hence, heuristically, by defining $q_{\widehat{x}}$ by
\begin{equation}\label{eq:q_x=Twedgedd^c|W_x|^2}
    q_{\widehat{x}}:=\frac{1}{L(\widehat{x})}(T^r\wedge dd^c|W_{\widehat{x}}|^2)\circ\pi_0\lfloor_{\xi^u_{\widehat{x}}},
\end{equation}
with
\begin{equation}\label{eq:L(x)}
    L(\widehat{x}):=(T^r\wedge dd^c|W_{\widehat{x}}|^2)(\pi_0(\xi_{\widehat{x}}^u)),
\end{equation}
\newpage
\noindent we obtain a measure whose Jacobian is $d^{rn} e^{2n(\lambda_{r+1}+\cdots+\lambda_k \pm \varepsilon)}$. If $\lambda_{r+1}=\cdots=\lambda_k=\tfrac12\mathrm{Log}\, d$, we recover the Jacobian of $\mu$ up to a factor $e^{\pm 2n\varepsilon}$. This error term is removed by Ledrappier-Strelcyn Lemma \ref{lemma:Ledrappier-Strelcynlemma}, which yields the integral formula in Item \textit{4.} However, the precise definition of $q_{\widehat{x}}$ and the implementation of this strategy present several difficulties and require new arguments.\\ 

\textbullet The first difficulty is that we cannot always define $q_{\widehat{x}}$ by Formula \eqref{eq:q_x=Twedgedd^c|W_x|^2}, because this formula requires that $\pi_0(\xi^u_{\widehat{x}})\subset\mathrm{Dom}(W_{\widehat{x}})$, which is the domain of definition of the measure \eqref{eq:Twedgedd^c|W_x|^2}. To define $q_{\widehat{x}}$ in general, we use a special Pesin box $\mathcal{A}$ introduced in Section \ref{sec:constructionofA}, Equation \eqref{eq:realborelianA}. This Pesin box satisfies $\pi_0(\xi^u_{\widehat{x}})\subset\mathrm{Dom}(W_{\widehat{x}})$ when $\widehat{x}\in\mathcal{A}$, which allows us to define $q_{\widehat{x}}$ by Formula \eqref{eq:q_x=Twedgedd^c|W_x|^2} in this case. Then, for a general $\widehat{x}$, we use the Poincaré recurrence theorem to obtain $n\geq0$ such that $\widehat{x}_{-n}\in\mathcal{A}$. This allows us to define $q_{\widehat{x}}$ proportionally to $q_{\widehat{x}_{-n}}\circ \widehat{f}^{-n}$ (see Definition \ref{defn:64}), and we check that this definition does not depend on the choice of $n$ such that $\widehat{x}_{-n}\in\mathcal{A}$ (see Proposition \ref{lemma:66}).\\  

\textbullet The second difficulty is that in this definition of $q_{\widehat{x}}$, we first need to check that $L(\widehat{x})$ (defined by \eqref{eq:L(x)} when $\widehat{x}\in\mathcal{A}$) is not equal to $0$. For this purpose, in Section \ref{sec:pluripotentialaspects} we study the measures defined on a germ $(\mathbb{C}^k,0)$ of the form $T^r\wedge dd^c|w_{r+1}|^2\wedge\cdots\wedge dd^c|w_k|^2$, where $(w_{r+1},\cdots,w_k)$ are the last $(k-r)$ coordinates of $\mathbb{C}^k$. We prove Theorem \ref{thm:SupportofT^rwedgedd^c|w|^2}, which shows that if $T$ is a current with Hölder continuous local potential, and if $0$ lies in the support of the Monge-Ampère mass $T^k$, then $0$ also lies in the support of $T^r\wedge dd^c|w_{r+1}|^2\wedge\cdots\wedge dd^c|w_k|^2$.   

This theorem was already known in the case $r=1$, it is contained in the work of Briend-Duval \cite{BriDuv99}, see also the survey of Sibony \cite[Cor. A.10.3, p.~181]{sib99}. Theorem \ref{thm:SupportofT^rwedgedd^c|w|^2} extends this result to $r\geq 1$. Since this is a local problem, in Theorem \ref{thm:SupportofT^rwedgedd^c|w|^2} we prove the following implication 
$$(T^r\wedge dd^c|w_{r+1}|^2\wedge\cdots\wedge dd^c|w_k|^2)(\mathbb{D}^k)=0\ \Longrightarrow\ T^k(\mathbb{D}^k)=0,$$
where $\mathbb{D}^k$ is a polydisc centered at $0$. The proof of this implication is as follows. Let $G$ be the potential of $T$. We then show that $(T^r \wedge dd^c |w_{r+1}|^2 \wedge \cdots \wedge dd^c |w_k|^2)(\mathbb{D}^k) = 0$ implies, for all $w \in \mathbb{D}^{k-r}$, that $(dd^c G_w)^r = 0$ on $\mathbb{D}^r$, where $G_w(z) := G(z, w)$. Then, for any ball $B(a, s) \Subset \mathbb{D}^k$ and any point $p_0 = (z, w) \in B(a, s)$, denoting by $B(a_0, s_0)$ the orthogonal projection of the ball $B(a,s)$ onto $(\mathbb{D}^r \times \{w\})\cap B(a,s)$, we construct a holomorphic map $h : B_{\mathbb{C}^r} \to B(a_0, s_0)$ from the unit ball of $\mathbb{C}^r$ such that $(dd^c G \circ h)^r = 0$ on $B_{\mathbb{C}^r}$ and $h(0) = p_0$. Using the comparison principle of Bedford–Taylor (see Theorem \ref{thm:strongcomparisonprinciple}) applied to $G \circ h$, we deduce that $G$ coincides with the solution of the Dirichlet problem on $B(a, s) \subset \mathbb{D}^k$ (see Section \ref{sec:Dirichlet}). It then follows that $T^k(\mathbb{D}^k) = 0$, which proves Theorem \ref{thm:SupportofT^rwedgedd^c|w|^2}.\\  

\textbullet The third difficulty is to prove, from our definition of $q_{\widehat{x}}$ (see Definition \ref{defn:64}), that the measure $q_{\widehat{x}}$ has the property that when $\lambda_{r+1} = \cdots = \lambda_k$, $q_{\widehat{x}}$ depends only on the atom $\xi^u_{\widehat{x}}$ containing $\widehat{x}$, and not on the particular choice of $\widehat{x}$ in this atom. This property is necessary for the measures $q_{\widehat{x}}$ to be the conditional measures of a probability measure $\widehat{q}$ on the atoms $\xi^u_{\widehat{x}}$ of $\xi^u$, which is necessary to obtain $\widehat{q} = \widehat{\mu}$ in Item \textit{6.} of Theorem \ref{thm:TheoremB}. We prove Theorem 1.1 for this purpose.

\subsubsection{Proof of Theorem \ref{thm:Patching} ; \texorpdfstring{$q_{\widehat{x}}$}{TEXT} depends only on \texorpdfstring{$\xi^u_{\widehat{x}}$}{TEXT}}

We show that the normal forms of Berteloot-Dupont-Molino \cite{bdm07} satisfy new properties that were not present in the original statement. In Theorem \ref{thm:Patching}, we prove that on a Pesin box one can extend the domain of definition of the normal forms (see Theorem \ref{thm:localfoliations}). In particular, this allows us to give a common domain of definition for all coordinates $W_{\widehat{x}}$ with $\widehat{x}\in\mathcal{A}$ (see Proposition \ref{eq:pi(eta)includedinDom(f-n)}), allowing to define the measures $q_{\widehat{x}}$ as explained in the previous paragraph. It then remains to prove that when $\lambda_{r+1}=\cdots=\lambda_k$, $q_{\widehat{x}}$ depends only on $\xi_{\widehat{x}}^u$ and not on the particular choice of $\widehat{x}$ in this atom.  

Theorem \ref{thm:Patching} allows us to check this last property by proving that the $1$-forms $dW_{\widehat{x}}$ patch along unstable manifolds of Pesin boxes. More precisely, if $\widehat{x}$ and $\widehat{y}$ lie in the same unstable manifold $W^u(\widehat{z}, R) \subset P$ of a Pesin box $P$, then we prove that there exists a matrix $C \in \mathrm{GL}_{k-r}(\mathcal{O}(B(z_0, R)))$ (depending on $\widehat{x}$ and $\widehat{y}$) such that
$$ dW_{\widehat{y}} = C dW_{\widehat{x}} \ \mathrm{on} \ B(z_0, R), $$
and in the case $\lambda_{r+1} = \cdots = \lambda_k$, we prove that $C \in \mathrm{GL}_{k-r}(\mathbb{C})$ is constant. So, in this case, this implies the following equality of measures
$$ (T^r \wedge dd^c |W_{\widehat{y}}|^2) = |\mathrm{det}(C)|^2 \times (T^r \wedge dd^c |W_{\widehat{x}}|^2), \quad |\mathrm{det}(C)|^2 \in \mathbb{C}^*.$$ 
Then taking $P = \mathcal{A}$ we have $W^u(\widehat{z}, R) \supset \xi_{\widehat{x}}^u$, so if $\widehat{y} \in \xi_{\widehat{x}}^u$, we deduce that $q_{\widehat{x}} = q_{\widehat{y}}$ (see Proposition \ref{prop:q_x=q_ysieta_x=eta_y}). 

The proof of Theorem \ref{thm:Patching} is dynamical: we write the elements of $dW_{\widehat{y}_{-n}}$ as linear combinations of the elements of $dZ_{\widehat{x}_{-n}}$ and $dW_{\widehat{x}_{-n}}$, and pull them back by $f^{-n}_{\widehat{z}}$. Then, since ${\lambda_r>\lambda_{r+1}}$, the strongest contractions in the $Z$-directions force the corresponding terms to vanish which gives $dW_{\widehat{y}}=CdW_{\widehat{x}}$ for some $C\in\mathrm{GL}_{k-r}({\mathcal{O}(B(z_0,R))})$. If moreover $\lambda_{r+1} = \cdots = \lambda_k$, we can show that each coefficient of $C$ is a constant function, and thus $C \in \mathrm{GL}_{k-r}(\mathbb{C})$, see Lemma \ref{lemma:laconvergenceuniformede_un|Fn|2}.

\subsubsection{The equality \texorpdfstring{$\widehat{q}=\widehat{\mu}$}{TEXT} and the proof of Theorem \ref{thm:TheoremA}}

The proof of $\widehat{q}=\widehat{\mu}$ when $\mu$ has $k-r$ minimal exponents relies on the formulas in Items \textit{4.} and \textit{5.} of Theorem \ref{thm:TheoremB}, on Jensen's inequality, and on the properties of $\xi^u$, which allow one to show that the conditional measures $\mu_{\widehat{x}} = q_{\widehat{x}}$ coincide, implying the equality $\widehat{q} = \widehat{\mu}$, see \cite{led84,dup06}. We reproduce these arguments in Theorem \ref{thm:q_x((g-neta)_x)=mu_x((g-neta)_x)}.

The proof of Theorem \ref{thm:TheoremA} then proceeds as follows. First, as explained in the remark below Theorem \ref{thm:TheoremA}, the equivalence \eqref{eq:equivalencek-1exposants} follows from the implication \eqref{eq:DujardinTheoremII}, the equivalences \eqref{eq:equivII}, and Dujardin's theorem \eqref{eq:DujardinTheorem}. It therefore remains to prove the implication \eqref{eq:DujardinTheoremII}. To do this, we use the equality $\widehat{q}=\widehat{\mu}$ and the fact that $(\pi_0)_*\widehat{q}\ll T^r\wedge\omega_{\mathbb{P}^k}^{k-r}$, see Lemma \ref{lemma:pi_0_*q_x<<Twedgeomega}.\\

\remerciements{The author thanks his Ph.D advisor C. Dupont for introducing him to the problem of absolute continuity, for teaching him the classical methods on partitions and entropy, and also for his kindness and constant support. This work was conducted within the the France 2030 framework programme, Centre Henri Lebesgue ANR-11-LABX-0020-01.}

\section{Pluripotential aspects}\label{sec:pluripotentialaspects}

\subsection{Comparison of \texorpdfstring{$T^r \wedge dd^c |w|^2$}{TEXT} with the Monge-Ampère mass \texorpdfstring{$T^k$}{TEXT}}\label{sec:Bedford-Taylor-Theory}

For general accounts on currents on complex manifolds, we refer to Bedford-Taylor \cite{bedtay76,bedtay82}, Demailly \cite{agbook}, Dinh-Sibony \cite{dinsib} and Sibony \cite{sib99}. Let us simply recall the results from Bedford-Taylor Theory and from the pluripotential aspects of endomorphisms of $\mathbb{P}^k$ that will be use full for our purposes. For an open set $\Omega$, we denote by $\mathrm{Psh}_{C,\alpha}(\Omega)$ the set of plurisubharmonic functions on $\Omega$ that are $(C,\alpha)$-Hölder continuous. 

\begin{prop}[Dinh-Sibony {\cite[Proposition 1.18]{dinsib}}]\label{prop:HolderpotentialsfortheGreenCurrent} 
    For each point of $\mathbb{P}^k$, there exists a chart $\Omega$ containing this point on which the Green current $T$ of $f$ has a potential $G$ belonging to $\mathrm{Psh}_{C,\alpha}(\Omega)$ for some constants $C,\alpha>0$.
\end{prop}

In the sequel of this section, we will mainly be interested in studying locally, on a chart $\Omega$, the support of $T^r$. Thanks to Proposition \ref{prop:HolderpotentialsfortheGreenCurrent} it is sufficient to work with the local model of the polydisc $\Omega:=\mathbb{D}^k$ on $\cmplex^k$, with $T$ being given by $T=dd^cG$, with $G\in\mathrm{Psh}_{C,\alpha}(\Omega)$ (note that $G$ is also Hölder continuous on $\overline{\Omega}$ and hence it is bounded on $\overline{\Omega}$). In this setting, we will denote $(z,w)=(z_1,\cdots,z_r,w_{r+1},\cdots,w_k)$ the standard coordinates of $\cmplex^k$. Our goal in the following subsections is to prove the following theorem. Let us fix now $1\leq r\leq k-1$, and let $dd^c|w|^2$ be defined by:
$$dd^c|w|^2:=dd^c|w_{r+1}|^2\wedge\cdots\wedge dd^c|w_k|^2.$$

\begin{thm}\label{thm:SupportofT^rwedgedd^c|w|^2}
    Let $G \in \mathrm{Psh}_{C, \alpha}\left(\mathbb{D}^k\right)$ for some $C, \alpha > 0$, and let us fix $1 \leq r \leq k - 1$. Let us denote $T := dd^c(G)$. If $(T^r \wedge dd^c |w|^2)(\mathbb{D}^k) = 0$, then $T^k(\mathbb{D}^k) = 0$.
\end{thm}

This theorem was already known for $r=1$ thanks to the work of Briend-Duval \cite{BriDuv99}, see also Sibony \cite[Cor. A.10.3, p.~181]{sib99}. Here, we extend this result using the solution of the Dirichlet problem on balls (Theorem \ref{thm:Dirichlet}) and the Comparison Principle \ref{thm:strongcomparisonprinciple} of Bedford-Taylor as main tools. The proof is divided into three steps, which are described in the following three subsections.

\subsection{Step 1 : Rewriting of the assumption \texorpdfstring{$(T^r \wedge dd^c |w|^2)(\mathbb{D}^k)=0$}{TEXT}.}

We write $\Omega=\Omega'\times\Omega''$ with $\Omega':=\mathbb{D}^r$ and $\Omega'':=\mathbb{D}^{k-r}$, and fix $G\in\mathrm{Psh}_{C,\alpha}(\Omega)$. For each $w\in\Omega''$, let $G_w:\Omega'\to\mathbb{R}$ be defined by $G_w(z)=G(z,w)$. Then $G_w\in\mathrm{Psh}_{C,\alpha}(\Omega')$, so we can define $T_w:=dd^cG_w$ and $T_w^r:=(dd^cG_w)^r$ on $\Omega'$.

\begin{prop}\label{prop:Fubini} The following two points hold:
    \begin{enumerate}
        \item Let $\psi=\chi(w)\varphi(z)$ be a test function, with $\varphi\in C^{\infty}_{c}(\Omega')$ and $\chi\in C^{\infty}_c(\Omega'')$. Then one has  
        \begin{equation}\label{eq:Fubini}
            \langle T^r\wedge dd^c|w|^2,\psi\rangle = \int_{\Omega''}\chi(w)\langle T^r_{w},\varphi\rangle\ \mathrm{dLeb}_{\mathbb{C}^{k-r}}(w).
        \end{equation}
        \item If $(T^r\wedge dd^c|w|^2)(\mathbb{D}^k)=0$, then $T^r_{w}=0$ on $\Omega'$ for all $w\in\Omega''$.
    \end{enumerate}
\end{prop}
\noindent\textbf{\underline{Proof:}}\\
\noindent\textit{1.} We have by definition of $dd^c$ in the sense of currents (integration by parts):
\begin{align*}
    \int_{\Omega}\psi\ T^r\wedge dd^c|w|^2 & = \int_{\Omega} G\ T^{r-1}\wedge dd^c\psi\wedge dd^c|w|^2.
\end{align*}
In the system of coordinates $(z, w) \in \Omega' \times \Omega''$, by wedging a current with $dd^c |w|^2$, only the part of the current computed in the $z$ coordinates remains. Thus, we deduce that the following equality of currents holds on $\Omega$: 
$$ G\left(T^{r-1} \wedge dd^c \psi \wedge dd^c |w|^2\right) = \chi(w)\left(G_w\ T_w^{r-1} \wedge dd^c \varphi \wedge dd^c |w|^2\right).$$
So, integrating and using Fubini's Theorem, and then integrating by parts, we obtain
\begin{align*}
    \int_{\Omega}\psi\ T^r\wedge dd^c|w|^2 & = \int_{\Omega''}\chi(w)\left(\int_{\Omega'}G_{w}\ T_{w}^{r-1}\wedge dd^c\varphi\right)\ \mathrm{dLeb}_{\mathbb{C}^{k-r}}(w)\\
                                                        & = \int_{\Omega''}\chi(w)\left(\int_{\Omega'}\varphi \ T_{w}^{r-1}\wedge dd^cG_{w}\right)\ \mathrm{dLeb}_{\mathbb{C}^{k-r}}(w),
\end{align*}
and so we get $\int_{\Omega}\psi\ T^r\wedge dd^c|w|^2=\int_{\Omega''}\chi(w)\langle T_{w}^r,\varphi\rangle\ \mathrm{dLeb}_{\mathbb{C}^{k-r}}(w)$.\\

\noindent\textit{2.} Let us now fix $U \Subset \Omega''$ a non-empty open subset, and let $\chi \geq0$ be a test function on $\Omega''$ such that $\chi = 1$ on $U$. Let $\varphi \geq 0$ be an arbitrary test function on $\Omega'$. Applying \eqref{eq:Fubini} to $\psi(z, w) := \chi(w) \varphi(z)$ gives $\int_{\Omega''}\chi(w)\langle T_{w}^r,\varphi\rangle\ \mathrm{dLeb}_{\mathbb{C}^{k-r}}(w)=0$.
Since $\chi(w) \langle T_w^r, \varphi \rangle \geq 0$ and $U \subset \Omega''$, it follows that there exists a Borel set $A \subset U$ of full Lebesgue measure in $U$ such that for all $w \in A$, we have $\chi(w) \langle T_w^r, \varphi \rangle = 0$. Since $\chi = 1$ on $U$, we deduce that $\langle T_w^r, \varphi \rangle = 0$ for all $w \in A$.

Next, we prove that $\langle T_w^r, \varphi \rangle = 0$ for all $w \in U$. Fix $w \in U$. Since $A$ is dense in $U$ (because $A$ has full Lebesgue measure in $U$ and $U$ is open), there exists a sequence $w_n \in A$ such that $\lim_n w_n = w$. By integration by parts we have
\begin{align*}
    \langle T^r_w, \varphi \rangle & = \int_{\Omega'}G_w(dd^cG_w)^{r-1}\wedge dd^c\varphi\\
    & = \int_{\Omega'}G_{w_n}(dd^cG_w)^{r-1}\wedge dd^c\varphi + \int_{\Omega'}(G_w-G_{w_n})(dd^cG_w)^{r-1}\wedge dd^c\varphi.
\end{align*}
By the Chern-Levine-Nirenberg inequality, see \cite[p.146]{agbook}, there exists a constant $C_{\varphi} > 0$ (depending only on $\mathrm{Supp}(\varphi)$ and $\Omega'$) such that the second integral is less than or equal to $||G_w - G_{w_n}||_{L^{\infty}} C_{\varphi} ||\varphi||_{C^2} ||G_w||_{L^{\infty}}^{r-1}$. Therefore, by using the fact that $G$ is $(C, \alpha)$-Hölder, we deduce that
\begin{align*}
    \langle T^r_w, \varphi \rangle & \leq \int_{\Omega'}G_{w_n}(dd^cG_w)^{r-1}\wedge dd^c\varphi + C||w-w_n||^{\alpha}C_{\varphi}||\varphi||_{C^2}||G_w||_{L^{\infty}}^{r-1}.
\end{align*}
By repeating these arguments inductively, we obtain that
\begin{align*}
    \langle T^r_w, \varphi \rangle & \leq \int_{\Omega'}G_{w_n}(dd^cG_{w_n})^{r-1}\wedge dd^c\varphi + C||w-w_n||^{\alpha}C_{\varphi}||\varphi||_{C^2}\sum_{l=0}^{r-1}||G_{w_n}||_{L^{\infty}}^{l}||G_w||_{L^{\infty}}^{r-1-l}.
\end{align*}
We have, $\int_{\Omega'}G_{w_n}(dd^cG_{w_n})^{r-1}\wedge dd^c\varphi=\langle T^r_{w_n}, \varphi \rangle=0$ because $w_n\in A$. We then observe that the right hand side of the last inequality tends to 0 as $n$ goes to infinity, thus we get $\langle T^r_w, \varphi \rangle=0$. Since this is true for every test function $\varphi\geq0$, we have proved that $T_w^r=0$ on $\Omega'$. So, we have proved that $T^r_w=0$ for each $w\in U$, and for every open set $U\Subset\Omega''$. Therefore, we have proved $T^r_w=0$ for every $w\in \Omega''$.\qed 

\subsection{Step 2 : Construction of Monge-Ampère flat balls}

Let us denote $H_w:=\left\{(Z,W)\in \Omega'\times\Omega''\ :\ W=w\right\}=\Omega'\times\{w\}$ for each $w\in\Omega''$, and let us denote $B'\subset\Omega'$ the open unit ball of $\cmplex^r$. 

\begin{defn}\label{defn:boulesplates} 
    Let $a\in\Omega$, and let $s>0$ such that $B(a,s)\Subset\Omega$. Let us fix a point $p_0\in B(a,s)\cap H_w$ with $w\in\Omega''$. A Monge-Ampère flat ball centered at $p_0$ inside $B(a,s)$, is a continuous map $h:\overline{B'}\to \overline{B(a,s)}\cap H_w$ satisfying the following items:
    \begin{enumerate}
        \item $h(0)=p_0$ and $h(\partial B')\subset\partial B(a,s)$ and $h(B')\subset B(a,s)$.
        \item $h$ is holomorphic on $B'$.
        \item The Monge-Ampère mass of $G\circ h$ is null on $B'$, i.e., $(dd^cG\circ h)^r=0$ on $B'$.
    \end{enumerate}
\end{defn}

The interest of this definition is that we have the following existence lemma of Monge-Amp\`ere flat balls. The proof is elementary, it relies essentially on the transitivity of the action of the group $\mathrm{Aut}(B')$ of holomorphic automorphisms of the unit ball $B'$.

\begin{lemme}\label{Lemme:constructionbouleplatespluriharmoniques}
    Let $B(a,s)\Subset\Omega$ be a ball, and let $p_0\in B(a,s)$ be an arbitrary point. If $w\in\Omega''$ such that $p_0\in H_w$ satisfies $T_w^r=0$ on $\Omega'$, then there exist a Monge-Ampère flat ball centered at $p_0$ inside $B(a,s)$.
\end{lemme}
\noindent\textbf{\underline{Proof:}} Let $a':=(a_0,w)\in B(a,s)$ be the orthogonal projection of $a$ on $H_w$, then by putting $s_0:=\sqrt{s^2-||a-a'||^2}>0$, $B'(a_0,s_0):=B_{\mathbb{C}^{r}}(a_0,s_0) \times \{w\}\subset B(a,s)$ satisfies $\partial B'(a_0,s_0)=\partial B_{\mathbb{C}^{r}}(a_0,s_0) \times \{w\}\subset \partial B(a,s)$, and $p_0\in B'(a_0,s_0)$. Let us now construct a map $h:\overline{B'}\to\overline{B(a,s)}\cap H_w$ satisfying the three items of Definition \ref{defn:boulesplates}.

The map $\widetilde{h}(z):=(s_0z+a_0,w)$ is a biholomorphism from $B'$ to $B'(a_0,s_0)$ extending to a homeomorphism from $\overline{B'}$ to $\overline{B'(a_0,s_0)}$, and $h$ maps $\partial B'$ onto $\partial B'(a_0,s_0)$. Let us denote $\widetilde{p}_0:=\widetilde{h}^{-1}(p_0)$.
It is known that the group $\mathrm{Aut}(B')$ of biholomorphisms of $B'$ is transitive: for all $p,q\in B'$, there exists $\gamma\in\mathrm{Aut}(B')$ such that $\gamma(p)=q$ and $\gamma$ extends on $\overline{B'}$ such that $\gamma:\overline{B'}\to\overline{B'}$ is a homeomorphism (see for example Rudin \cite{rudin2008}). Thus there exists $\gamma\in\mathrm{Aut}(B')$ such that $\gamma(0)=\widetilde{p}_0$, and such that $\gamma$ extends to a homeomorphism on $\overline{B'}$. Then by defining $h:=\widetilde{h}\circ\gamma$, we obtain a homeomorphism $h:\overline{B'}\to\overline{B'(a_0,s_0)}$ such that $h(0)=p_0$ and such that $h(\partial B')=\partial B'(a_0,s_0)$. The map $h$ thus satisfies the two first item of Definition \ref{defn:boulesplates}. 

It remains to show that $(dd^cG\circ h)^r=0$ on $B'$. By denoting $\pi(z,w)=z$, we can observe that 
$\int_{B'}(dd^cG \circ h)^r = \int_{B'}\left[dd^cG_w \circ (\pi\circ h)\right]^r = \int_{B_{\mathbb{C}^{r}}(a_0,s_0)}|\mathrm{Jac}_{\mathbb{C}}(\pi\circ h)|^{2} (dd^cG_w)^r$. This last integral is equal to $0$
since $(dd^cG_w)^r=T^r_w = 0$ on $\Omega'\supset B_{\cmplex^r}(a_0,s_0)$ by assumption, thus $h$ finally satisfies the three required items of Definition \ref{defn:boulesplates}.\qed 

\subsection{Step 3 : Solving the Dirichlet problem on balls}\label{sec:Dirichlet}
  
We will now complete the proof of Theorem \ref{thm:SupportofT^rwedgedd^c|w|^2}. To do this, we will show that under the hypothesis of the theorem, every point of $\Omega$ is the center of a Monge-Ampère flat ball inside an open ball of $\Omega$, and from this, we will deduce that the potential $G$ is a solution to the Dirichlet problem on each of these open balls. Thus, the Monge-Ampère mass of $G$ will be zero on the open set $\Omega$, proving the theorem.

Let us begin by recalling the Dirichlet problem in our context. The Dirichlet problem on a ball $B \Subset \Omega$ associated with $G$ is the following PDE, where the unknown is a plurisubharmonic function $u$ on $B$:
$$
(\mathrm{Dir}_B): \left\{
    \begin{array}{ll}
        (dd^cu)^k=0 & \mbox{on } B\\
        u=G & \mbox{on } \partial B\\
        u &\mathrm{psh\ on}\ B\ \mathrm{and\ continuous\ on}\ \overline{B}
    \end{array}
\right.
$$

\begin{thm}[Bedford-Taylor {\cite{bedtay76}}, Klimek {\cite[Theorem 4.4.1]{Klimek91}}]\label{thm:Dirichlet}
    For each ball $B\Subset\Omega$, the solution to $(\mathrm{Dir}_B)$ is $u=\widehat{G}$ defined by:
    $$\widehat{G}=\sup\{v\ |\ v\in \mathrm{Psh}(B)\cap C^{0}(\overline{B}),\ v\leq G\ \mathrm{on}\ \partial B\}.$$
    Moreover, since $G$ is $\alpha-$Hölder, $\widehat{G}$ is $(\alpha/2)-$Hölder. 
\end{thm}

We will combine this theorem with the following Comparison Principle, originally due to Bedford--Taylor \cite{bedtay82} (see also Klimek \cite[Theorem 3.7.4]{Klimek91}).

\begin{thm}[Comparison Principle]\label{thm:strongcomparisonprinciple}
    Let $B\subset\cmplex^r$ be an open ball of $\cmplex^r$, $r\geq1$, and assume that $u$ and $v$ are two plurisubharmonic functions on $B$, continuous on $\overline{B}$. If $u\geq v$ on $\partial B$ and if $(dd^cu)^r=0$, then $u\geq v$ on $B$.
\end{thm}

\noindent\textbf{\underline{Proof of Theorem \ref{thm:SupportofT^rwedgedd^c|w|^2}:}} The proof follows the ideas of Sibony in \cite[Corollaire A.10.3, p. 181]{sib99}. Let us assume that $(T^r\wedge dd^c|w|^2)(\Omega)=0$, and let us prove that $T^k(\Omega)=0$. According to Proposition \ref{prop:Fubini}, we have $T^r_w=0$ on $\Omega'$, for all $w\in\Omega''$. Let us fix an open ball $B(a,s)\Subset\Omega$ and a point $p_0 \in B(a,s)$. Let $w\in\Omega''$ such that $p_0\in H_w$. We have $T^r_w=0$ on $\Omega'$, thus according to Lemma \ref{Lemme:constructionbouleplatespluriharmoniques}, there exists a continuous map $h: \overline{B'} \to \overline{B(a,s)}\cap H_w$ satisfying the three items of Definition \ref{defn:boulesplates}.
 
Let now $\widehat{G}_{a,s}$ be the solution of the Dirichlet problem $(\mathrm{Dir}_{B(a,s)})$ introduced above. According to Theorem \ref{thm:Dirichlet}, $\widehat{G}_{a,s} = \sup_v v$, where $v$ ranges over the set of functions $v \in \mathrm{Psh}(B(a,s)) \cap C^0(\overline{B(a,s)})$ satisfying $v \leq G$ on $\partial B(a,s)$. For such a function $v$, since $h(\partial B')\subset\partial B(a,s)$, we have $v \circ h \leq G \circ h$ on $\partial B'$. Applying the Comparison Principle \ref{thm:strongcomparisonprinciple}, it implies $v \circ h \leq G \circ h$ on $B'$. In particular, since $h(0) = p_0$, we have $v(p_0) \leq G(p_0)$. By taking the supremum over the functions $v$, we deduce that $\widehat{G}_{a,s}(p_0) \leq G(p_0)$. By the definition of $\widehat{G}_{a,s}$, we also have $G(p_0) \leq \widehat{G}_{a,s}(p_0)$, so finally $\widehat{G}_{a,s}(p_0) = G(p_0)$. Since the preceding arguments are valid for all the points $p_0\in B(a,s)$, we deduce that $G=\widehat{G}_{a,s}$. In particular, $G$ is the solution of the Dirichlet problem on $B(a,s)$, for all the open balls $B(a,s)\Subset\Omega$, so we get $(dd^cG)^k=0$ on $\Omega$, completing the proof.\qed

\section{Ergodic aspects}\label{sec:Oseledecnormalforms}

We assume in what follows that $k \geq 2$. Let $f$ be an endomorphism of $\mathbb{P}^k$ of degree $d \geq 2$. Let $\mu$ denote its equilibrium measure, and let $\lambda_1 \geq \cdots \geq \lambda_k$ be the Lyapunov exponents of $\mu$. Denote by $\mathrm{Crit}(f)$ the critical set of $f$, and set $\mathcal{C} := \bigcup_{n \in \mathbb{Z}} f^n(\mathrm{Crit}(f))$. It is known that $\mu(\mathcal{C}) = 0$, see \cite{dinsib,sib99}.

{\subsection{Oseledec Theorem}}

\begin{thm}\label{thm:oseledec}
Assume that there exists $r \in\{1, \dots, k-1\}$ such that $\lambda_r > \lambda_{r+1}$. Then there exists a totally invariant Borel set $A_{os}$ of full $\mu$-measure, disjoint from $\mathcal{C}$, and a measurable bundle $V_s \subset T\mathbb{P}^k$ (depending on $r$) defined over $A_{os}$ such that $\mathrm{dim}_{\mathbb{C}} V_s(x) = k - r$ and satisfying
$$\forall x \in A_{os} \ , \ \forall \vec v\in V_s(x)\backslash\{0\},\ \lim_{n\to\infty}\frac{1}{n}\LLog||d_xf^n\cdot\vec v|| \in\{\lambda_{r+1},\cdots,\lambda_k\}.$$
Moreover, $d_xf^n:V_s(x)\to V_s(f^nx)$ is a $\cmplex-$linear isomorphism for each $n$, and if $\vec v\in T_x\mathbb{P}^k\backslash V_s(x)$ then the limit above belongs to $\{\lambda_1,\cdots,\lambda_r\}$. We denote 
$$\chi_f(x):=\LLog\ \vert\mathrm{det}_{\cmplex}\left(d_xf:V_s(x)\to V_s(fx)\right)\vert.$$
The following points hold :
\begin{enumerate}
    \item $\chi_f\in L^1(\mu)$,
    \item $\forall x\in A_{os}$, $\frac{1}{n}\sum_{i=0}^{n-1}\chi_f(f^i(x))=\frac{1}{n}\chi_{f^n}(x)$. In particular, $\int\chi_f\ \mathrm{d}\mu = \lambda_{r+1}+\cdots+\lambda_k$. 
    \item $\chi_{f^n}\circ\pi_0\in L^1(\widehat{\mu})$ and $\int\chi_{f^n}(x_0)\ \mathrm{d}\widehat{\mu}(\widehat{x}) = n(\lambda_{r+1}+\cdots+\lambda_k)$.
    \end{enumerate}
\end{thm}
The construction of $V_s$ is classical, we refer to Katok-Hasselblatt's book {\cite{kathas95}}. The proof of the first item can be found in  {\cite[Section 3.7]{sib99}}. The second one comes from the chain rule formula and Birkhoff ergodic theorem. The third item relies on $(\pi_0)_*\widehat{\mu}=\mu$.
   
{\subsection{Normal forms for inverse branches}\label{sec:normalforms}}

Let $X := \mathbb{P}^k \setminus \mathcal{C}$ and let $\widehat{X} := \{\widehat{x} \in \widehat{\mathbb{P}^k} : x_n \in X, \ \forall n \in \mathbb{Z}\}$. This is a totally $\widehat{f}-$invariant subset of $\widehat{\mathbb{P}^k}$ of full $\widehat{\mu}-$measure. For each $\widehat{x} \in \widehat{X}$, there exists $\eta>0$ and a sequence of holomorphic maps $f^{-n}_{\widehat{x}} : B(x_0, \eta) \to \mathbb{P}^k$ such that $f^{-n}_{\widehat{x}}(x_0) = x_{-n}$ and $f^n \circ f^{-n}_{\widehat{x}} = \mathrm{Id}$ on $B(x_0, \eta)$ for every $n \geq 0$. The existence of these inverse branches has been proved by Briend and Duval in the article \cite{BriDuv99}.

The following theorem was proved by Berteloot-Dupont-Molino {\cite{bdm07}}, it gives normal forms for the iterated inverse branches of Briend-Duval. We refer to the articles {\cite{jonvar02}} by Jonsson-Varolin  and {\cite{bertelootdupont2019}} by Berteloot-Dupont  for related results. A real-valued measurable function $\varphi$ on $\widehat{\mathbb{P}^k}$ is $\varepsilon-$tempered if $e^{-|n|\varepsilon}\varphi \leq \varphi \circ \widehat{f}^n \leq  e^{|n|\varepsilon} \varphi$ for every $n \in \mathbb Z$. For each $1\leq j \leq k$, we say that the Lyapunov exponents $\lambda_{j}\geq\cdots\geq\lambda_k$ are resonant, if there exists a $(k-j)-$tuple $(c_i)_{i=j+1,\cdots,k}\in\mathbb{N}^{k-j}$ of length $\geq 2$ such that $\sum_{i=j+1}^kc_i\lambda_i=\lambda_j$. Note that if $\lambda_j=\cdots=\lambda_k$ then these exponents are not resonant. 

\begin{thm}\label{thm:normalforms} Assume that there exists $r\in\{1,\cdots,k-1\}$ such that $\lambda_{r}>\lambda_{r+1}$ and such that $\lambda_{r+1}\geq\cdots\geq\lambda_k$ are not resonant. Then there exists $\varepsilon_0>0$ small enough such that for all $0<\varepsilon\leq \varepsilon_0$ the following holds. There exist a totally $\widehat{f}-$invariant Borel set $\mathcal{FN}\subset\widehat{X}$ of full $\widehat{\mu}-$measure, $\varepsilon-$tempered functions ${\rho},\eta : {\mathcal{FN}} \to \left]0,1\right]$, $\varepsilon-$tempered functions $\beta,\Gamma : {\mathcal{FN}} \to [1,+\infty[$, and a function $N : {\mathcal{FN}} \to \mathbb{N}$ satisfying the following properties for every $\widehat{x}\in \mathcal{FN}$ :
\begin{enumerate}
    \item There is an injective holomorphic map $\sigma_{\widehat{x}} : B(x_0,2\eta(\widehat{x})) \to \mathbb{D}^k(\rho(\widehat{x}))$ such that :
    \begin{enumerate}\label{enum:propertiesofd(xi)}
        \item $\sigma_{\widehat{x}}(x_0)=0$, $x_0\in A_{os}$ and $V_s(x_0)=d_{0}\sigma_{\widehat{x}}^{-1}(\{0\}^{r}\times\cmplex^{k-r})$,
        \item $\forall p,q \in B({x_0},2\eta(\widehat{x}))$ , $\frac{1}{2}d_{\mathbb{P}^k}(p,q)\leq||\sigma_{\widehat{x}}(p)-\sigma_{\widehat{x}}(q)||\leq\beta(\widehat{x})d_{\mathbb{P}^k}(p,q)$.
    \end{enumerate}
    \item There is a sequence of holomorphic maps $(f^{-n}_{\widehat{x}})_n$  such that $f^n\circ f^{-n}_{\widehat{x}}=\mathrm{Id}$ on\\ $B(x_{0},2\eta(\widehat{x}))$,  $f^{-n}_{\widehat{x}}(x_0) = x_{-n}$ and $\mathrm{Lip}(f^{-n}_{\widehat{x}})\leq \beta(\widehat{x})e^{-n(\lambda_k-\varepsilon)}$ for every $n\geq0$. 
    \item The following diagram commutes for any $n\geq N(\widehat{x})$:
    \begin{equation}\label{eq:diagramnormalforms}
        \xymatrix{
            B({x_{-n}},2\eta(\widehat{x}_{-n})) \ar[d]_{ \sigma_{\widehat{x}_{-n}} } & &  B({x_0},2\eta(\widehat{x})) \ar[ll]_{f^{-n}_{\widehat{x}}} \ar[d]^{ \sigma_{\widehat{x}}}  \\
            \mathbb{D}^k(\rho(\widehat{x}_{-n}))  & & \mathbb{D}^k(\rho({\widehat{x}}))  \ar[ll]^{R_{n,\widehat{x}}}
            }
    \end{equation}
    \item The map $R_{n,\widehat{x}}:\cmplex^k\to\cmplex^k$ is an invertible polynomial map and it has the form 
    \end{enumerate}
        \begin{equation}\label{eq:Rnxwithnoresonancesforlowerlyapunovexponents}
           \left(\alpha_{n,\widehat{x}}^{1}z_1 + P_{1,n}(z_2,\cdots,w_k),\cdots,\alpha_{n,\widehat{x}}^{r}z_r + P_{r,n}(w_{r+1},\cdots,w_{k}),\ \beta_{n,\widehat{x}}^{r+1}w_{r+1},\cdots,\beta_{n,\widehat{x}}^{k}w_{k}\right).
        \end{equation}
    \begin{enumerate}
        \item[5.] $e^{-n(\lambda_j+\varepsilon)}\leq |\alpha^j_{n,\widehat{x}}|\leq e^{-n(\lambda_j-\varepsilon)}$ and $e^{-n(\lambda_l+\varepsilon)}\leq |\beta^l_{n,\widehat{x}}| \leq e^{-n(\lambda_l-\varepsilon)}$.
        \item[6.] $P_{j,n}=\sum_{m=1}^{m_j}\gamma_{n,m,\widehat{x}}^j (z_{j+1},\cdots,w_k)^{C_{j,m}}$, where the length of $C_{j,m}\in\mathbb{N}^{k-j}$ is $\geq2$.
        \item[7.] Each $C_{j,m}$ is independent of $n$ and $|\gamma_{n,m,\widehat{x}}^j|\leq \Gamma(\widehat{x})e^{-n(\lambda_r-\varepsilon)}$.
        \item[8.] If $\lambda_{j}\geq\cdots\geq\lambda_k$ are not resonant, we can require that $\gamma_{n,m,\widehat{x}}^j=0$ for each $m$.
        \item[9.] {There exists $\displaystyle{C_0\!>\!0}$ such that for each $n\geq0$, there exists $0<r_n<\eta(\widehat{x})$ satisfying}
        \begin{equation}\label{eq:lemmedeproc}
            (f^{n})^*\omega_{\mathbb{P}^k}^{k-r}\geq C_0e^{2n(\lambda_{r+1}+\dots+\lambda_k-\varepsilon)}\omega_{\mathbb{P}^k}^{k-r}\ \mathrm{on}\ B(x_0,r_n).
        \end{equation}
    \end{enumerate}
In particular, if $\lambda_1\geq\cdots\geq\lambda_k$ are not resonant, then $R_{n,\widehat{x}}$ becomes a linear map of the form 
$$R_{n,\widehat{x}}=(\alpha_{n,\widehat{x}}^1z_1,\cdots,\alpha_{n,\widehat{x}}^r z_{r},\beta_{n,\widehat{x}}^{r+1}w_{r+1},\cdots,\beta_{n,\widehat{x}}^kw_k).$$
\end{thm}

We mention that the assumption $\lambda_r>\lambda_{r+1}$ is needed only in Item \textit{1.(a)} for the Oseledec bundle $V_s$, the other items hold without it. Similarly, the nonresonance of $\lambda_{r+1}\geq\cdots\geq\lambda_k$ is required only to get the linear form on $w_{r+1},\ldots,w_k$ in \eqref{eq:Rnxwithnoresonancesforlowerlyapunovexponents}, and the theorem has a similar version without this assumption. In particular, it allows a new proof of Theorem \ref{thm:ProofofDujardinTheorem} (originally due to Dujardin \cite[Thm.~3.6]{Duj12}) using the normal forms. The arguments are detailed in \cite{duptap23} for $\mathbb{P}^2$ and adapted here for $\mathbb{P}^k$.

\begin{thm}\label{thm:ProofofDujardinTheorem}
    Assume that there exists $r\in\{1,\cdots,k-1\}$ such that $\mu\ll T^r\wedge\omega_{\mathbb{P}^k}^{k-r}$. Then $\lambda_{r+1}=\cdots=\lambda_k=\frac{1}{2}\mathrm{Log}\ d$.
\end{thm}
\noindent\textbf{\underline{Proof:}} Let us denote $\Sigma_r:=\lambda_{r+1}+\dots+\lambda_k$, and let us denote $\sigma_{T^r}:=T^r\wedge \omega_{\mathbb{P}^k}^{k-r}$. Let $\psi$ be a positive measurable function such that $\mu=\psi\sigma_{T^r}$. For $\tau>0$, let $F_{\tau}$ be the set defined by
$$F_{\tau}:=\{\psi\geq\tau\}\cap\bigcap_{n\in\mathbb{N}}\mathrm{Leb}_{\sigma_{T^r}}(\psi\circ f^n),$$
where $\mathrm{Leb}_{\sigma_{T^r}}(\psi\circ f^n)$ is the set of Lebesgue points of the measure $(\psi\circ f^n)\sigma_{T^r}$. If $\tau>0$ is small enough, $\mu(F_{\tau})>0$ and there exists $\widehat{x}\in\mathcal{FN}\cap F_{\tau}$. Let $C_0>0$ and $r_n>0$ be given by Item \textit{9.} of Theorem \ref{thm:normalforms} (which also holds if $\lambda_r=\lambda_{r+1}$ or if $\lambda_{r+1},\ldots,\lambda_k$ are resonant, see \cite[Lemma 2.4]{duptap23}). Since $\mathcal{FN}$ avoids $\mathcal{C}=\bigcup_{n\in\mathbb{Z}}f^n(\mathrm{Crit}(f))$, for each $n$ there exists $\rho_n>0$ such that $f^n$ is injective on $B(x_0,\rho_n)$ and $f^n(B(x_0,\rho_n))\subset B(x_n,r_n)$.

Let us fix $n$ and let $B_{\rho}:=B(x_0,\rho)$ with $0<\rho\leq\rho_n$. Then $f^{n}$ is injective on $B_{\rho}$, thus we have 
$\mu(B_{\rho})=d^{-2n}\mu(f^{n}(B_{\rho}))$ i.e. $\int_{B_{\rho}}\psi\ \mathrm{d}\sigma_{T^r}=d^{-kn}\int_{f^{n}(B_{\rho})}\psi\ \mathrm{d}\sigma_{T^r}$. Using \eqref{eq:lemmedeproc} and ${f^{n}}^*T^r=d^{rn}T^r$ we deduce that 
$$\int_{B_{\rho}}\psi\ \mathrm{d}\sigma_{T^r}=d^{-kn}\int_{B_{\rho}}\psi\circ f^{n}\ \mathrm{d}(f^{n})^*\sigma_{T^r}\geq C_0e^{2n(\Sigma_r-\varepsilon)}d^{-(k-r)n}\int_{B_{\rho}}\psi\circ f^{n}\ \mathrm{d}\sigma_{T^r}.$$
By dividing by $\sigma_{T^r}(B_{\rho})$, we obtain
$$\frac{1}{\sigma_{T^r}(B_{\rho})}\int_{B_{\rho}}\psi\ \mathrm{d}\sigma_{T^r}\geq C_0e^{2n(\Sigma_r-\varepsilon)}d^{-(k-r)n}\frac{1}{\sigma_{T^r}(B_{\rho})}\int_{B_{\rho}}\psi\circ f^{n}\ \mathrm{d}\sigma_{T^r}.$$
Taking the limit as $\rho\to0^+$, we obtain by definition of $F_{\tau}$:
$$\psi(x_0)\geq C_0e^{2n(\Sigma_r-\varepsilon)}d^{-(k-r)n}\psi(x_{n})\geq C_0e^{2n(\Sigma_r-\varepsilon)}d^{-(k-r)n}\tau.$$
Since this is true for any $n$, we deduce that $\Sigma_r-\varepsilon\leq\frac{k-r}{2}\mathrm{Log}\ d$. Since this is true for all $\varepsilon\leq\varepsilon_0$, we get $\Sigma_r\leq\frac{k-r}{2}\mathrm{Log}\ d$, and so Briend-Duval inequality implies the result.\qed

\subsection{Unstable manifolds and Buzzi's partition}\label{sec:partitionbuzzi}

$\mathcal{FN}$ and the $\varepsilon$-tempered function $\eta(\widehat{x})$ can also be chosen to satisfy the following :

\begin{prop}
[Buzzi {\cite[Section 4]{buzzi2003}}, see also Dupont {\cite[Section 2.4]{dup06}}]
\label{prop:partitionPforinversebranchesI}
There exists a partition $\mathcal{P}=(\mathcal{P}_j)_{j\in\{1,\cdots,N\}}$ of a full $\mu-$measure set of $\mathbb{P}^k$, 
whose atoms are open sets of $\mathbb{P}^k$ and satisfy the two items below, where for each $x$ we have denoted $\mathcal{P}_{x}:=\mathcal{P}_j$ with $j$ such that $x\in\mathcal{P}_j$, and $\widehat{\mathcal{P}}_+:=\bigcap_{n\geq0}\widehat{f}^n(\pi_0^{-1}(\cup_j\mathcal{P}_j))$.
\begin{enumerate}
    \item The map $f$ is injective on each atom $\mathcal{P}_j$.
    \item For all $\widehat{x}\in\mathcal{FN}$, we have $\widehat{x}\in\widehat{\mathcal{P}}_{+}$ and $f^{-n}_{\widehat{x}}(B(x_0, 2\eta(\widehat{x}))) \subset \mathcal{P}_{x_{-n}}$ for every $n \geq 0$. 
\end{enumerate} 
\end{prop}

\begin{defn}\label{defn:invariantsetLambda_epsilon} \ 
\begin{enumerate}
    \item For $\widehat{x}\in\widehat{\mathcal{P}}_{+}$ we define the $\mathcal{P}-$address of $\widehat{x}$ by the sequence $(\mathcal{P}_{x_{-n}})_{n\geq0}$.
    \item Let $\widehat{\mathcal{P}}:=\bigcap_{n\in\mathbb{Z}}\widehat{f}^n(\pi_0^{-1}(\cup_j\mathcal{P}_j))\subset\widehat{\mathcal{P}}_{+}$ and let us define:
    $$\Lambda := \pi_0^{-1}(\mathrm{Supp}(\mu))\cap\pi_0^{-1}(A_{os}) \cap \mathcal{FN}\cap \widehat{\mathcal{P}}.$$ 
    $\Lambda$ is a totally invariant set of measure $\widehat{\mu}(\Lambda)=1$ on which $\vec v_s\circ\pi_0$, the $\mathcal{P}-$address and the inverse branches $f^{-n}_{\widehat{x}}$ are well defined.
\end{enumerate}
\end{defn}

\begin{defn}\label{defn:admissibleradius}
    Let $\widehat{x}\in\widehat{\mathcal{P}}_{+}\backslash\mathcal{C}$ be fixed, with $\mathcal{C}=\bigcup_{n \in \mathbb{Z}} f^n(\mathrm{Crit}(f))$. A pair $(r,\rho)$ such that $r\geq0,\ \rho>0$ is said to be admissible for $\widehat{x}$ if for every $n\geq0$, $f^{-n}_{\widehat{x}}$ is defined on $B(x_0,r)$ and satisfies $\mathrm{Lip}(f^{-n}_{\widehat{x}}|_{B(x_0,r)})\leq \rho e^{-n(\lambda_k-\varepsilon)}$ and $f^{-n}_{\widehat{x}}(B(x_0,r))\subset\mathcal{P}_{x_{-n}}$. A radius $r\geq0$ is said to be admissible for $\widehat{x}$ if there exists $\rho>0$ such that $(r,\rho)$ is admissible for $\widehat{x}$.
\end{defn}

Theorem \ref{thm:normalforms} and Proposition \ref{prop:partitionPforinversebranchesI} imply that for every $\widehat{x}\in\Lambda$ and for every $r\in[0,2\eta(\widehat{x})[$ and $\rho\geq\beta(\widehat{x})$, the pair $(r,\rho)$ is an admissible pair for $\widehat{x}$.

\begin{defn}\label{defn:unstablemanifolds}
    For $\widehat{x}\in\Lambda$ and for $r\geq0$ admissible for $\widehat{x}$, the unstable manifold associated to $(\widehat{x},r)$ is defined by:
    $$W^u(\widehat{x},r):=\left\{\widehat{z}\in \widehat{\mathbb{P}^k},\ \exists t_0\in B(x_0,r)\ :\ z_{-n}=f^{-n}_{\widehat{x}}(t_0),\ \forall n\in\mathbb{N}\right\}\subset\widehat{\mathcal{P}}_{+}.$$ 
    The $\mathcal{P}-$address is well defined and constant on $W^u(\widehat{x},r)$.
\end{defn}

\begin{lemme}\label{lemma:f-ncoincidesontheunstablemanifold} 
Let $\widehat{z}\in\Lambda$ be fixed and let $0\leq r<\eta(\widehat{z})$. 
\begin{enumerate}
    \item Let $\widehat{x}\in W^u(\widehat{z},r)\cap \Lambda$ and let $n\geq0$. Then there exists a holomorphic map $g_n:B(x_0,2\eta(\widehat{x}))\cup B(z_0,2\eta(\widehat{z}))\to\mathbb{P}^k$ such that $g_n=f^{-n}_{\widehat{x}}$ on $B(x_0,2\eta(\widehat{x}))$ and $g_n=f^{-n}_{\widehat{z}}$ on $B(z_0,2\eta(\widehat{z}))$.
    \item Let $\widehat{x}\in W^u(\widehat{z},r)$ be fixed. If  $\widehat{w}\in\widehat{\mathcal{P}}_{+}$ satisfies
    \begin{enumerate}
        \item $\widehat{x}$ and $\widehat{w}$ have the same $\mathcal{P}-$address,
        \item $w_0\in B(z_0,r)$,
    \end{enumerate}
    then $\widehat{w}\in W^u(\widehat{z},r)$. If moreover $w_0=z_0$ then $\widehat{w}=\widehat{z}$. 
\end{enumerate}
\end{lemme}
\noindent\textbf{\underline{Proof:}} 
Let us denote $V:=B(x_0,2\eta(\widehat{x}))\cup B(z_0,2\eta(\widehat{z}))$ and $U:=B(x_0,2\eta(\widehat{x}))\cap B(z_0,2\eta(\widehat{z}))$. Let $t_0$ be an element of $U$. By Proposition \ref{prop:partitionPforinversebranchesI} and Definition \ref{defn:invariantsetLambda_epsilon} we have for every $k \geq 0$ :
\begin{equation}\label{eq:P_x_-k=P_z_-kforallkgeq0}
    \mathcal{P}_{f^{-k}_{\widehat{x}}(t_0)}=\mathcal{P}_{x_{-k}}=\mathcal{P}_{z_{-k}}=\mathcal{P}_{f^{-k}_{\widehat{z}}(t_0)},
\end{equation}
where $\mathcal{P}_{x_{-k}}=\mathcal{P}_{z_{-k}}$ comes from $\widehat{x}\in W^u(\widehat{z},r)$.
If $f^{-j}_{\widehat{x}}(t_0)=f^{-j}_{\widehat{z}}(t_0)$ then $f(f^{-(j+1)}_{\widehat{x}}(t_0))$\\ $=f(f^{-(j+1)}_{\widehat{z}}(t_0))$. By using \eqref{eq:P_x_-k=P_z_-kforallkgeq0} and the fact that $f$ is injective on the atoms of $\mathcal{P}$ (see Proposition \ref{prop:partitionPforinversebranchesI}), we get $f^{-(j+1)}_{\widehat{x}}(t_0)=f^{-(j+1)}_{\widehat{z}}(t_0)$. An induction on $j\geq0$ thus implies $f^{-j}_{\widehat{x}}(t_0) = f^{-j}_{\widehat{z}}(t_0)$  for every $j\geq0$. Since this is true for any $t_0\in U$, $f^{-n}_{\widehat{x}}$ coincide with $f^{-n}_{\widehat{z}}$ on $U$ for every $n$. Therefore, we can define $g_n$ the map defined on $V$ by $f^{-n}_{\widehat{x}}$ on $B(x_0,2\eta(\widehat{x}))$ and by $f^{-n}_{\widehat{z}}$ on $B(z_0,2\eta(\widehat{z}))$. The conclusion of the first item then follows. The second item can be proved similarly.\qed

\section{Pesin boxes and rigidity along unstable manifolds}\label{sec:thePesinBoxA}

In this section we construct Pesin boxes $(P,r,\rho,\mathcal{T},c)$ (see Definition \ref{defn:pesinboxes}) for $\mu-$almost every $c\in\mathbb{P}^k$. In particular, we construct a special Pesin box $\mathcal{A}$ of $\widehat{\mu}-$positive measure on which we are able to control uniformly the dynamics. This set $\mathcal{A}$ plays a crucial role during the construction of the decreasing partition $\xi^u$ and the reconstruction of the topological entropy of $f$ in Section \ref{sec:partitionetaledrappier}. The set $\mathcal{A}$ will be also used to construct the family of conditional measures $q_{\widehat{x}}$ in Section \ref{definitionofthemeasuresq_xunderH_1}.

An other important part of this section is Theorem \ref{thm:localfoliations} which allows to patch the normal forms $dW_{\widehat{x}}=\left[dW_{\widehat{x}}^{r+1}\cdots dW_{\widehat{x}}^k\right]^{\mathrm{T}}$ along unstable manifolds inside Pesin boxes. As explained in the introduction, this property will be of crucial importance for proving that when $\lambda_{r+1} = \cdots = \lambda_k$, $q_{\widehat{x}} = q_{\widehat{y}}$ if $\widehat{x}$ and $\widehat{y}$ belong to the same atom of $\xi^u$. However, this result is also of independent interest, as it can be applied to the normal forms of any endomorphism $f$ whose Lyapunov exponents are not equal.

\subsection{Pesin boxes and Briend's theorem}\label{sec:defofPesinBoxes}

\begin{defn}\label{defn:pesinboxes}\ A Pesin box is a quintuplet $(P,r,\rho,\mathcal{T},c)$ where $c\in\mathbb{P}^k$ and
\begin{enumerate}
    \item $\mathcal{T}$ is a Borel subset of $\{\pi_0=c\}\cap\Lambda$.
    \item The pair $(r,\rho)$ is admissible in the sense of Definition \ref{defn:admissibleradius} for each $\widehat{z}\in\mathcal{T}$.
    \item $P=\bigsqcup_{\widehat{z}\in\mathcal{T}}W^u(\widehat{z},r)$. 
\end{enumerate}
\end{defn}
In the third item the unstable manifolds are pairwise disjoint, it is a consequence of the second item of Lemma \ref{lemma:f-ncoincidesontheunstablemanifold} and the fact that the unstable manifolds are centered at the same point $c$. Let $(P,r,\rho,\mathcal{T},c)$ be a Pesin box and let us define for any $\widehat{z}\in\mathcal{T}$:
$$\varphi_{\widehat{z}}  :  y_0\in B(c,r)\mapsto (f^{-n}_{\widehat{z}}(y_0))_{n\in\mathbb{Z}}\in W^u(\widehat{z},r) , $$ 
where $f^{n}_{\widehat{z}}=f^n$ if $n\geq0$. For any $\widehat{y}\in P$, let $\widetilde{\pi}(\widehat{y}):=\widehat{z}$ be the unique $\widehat{z}\in \mathcal{T}$ such that $\varphi_{\widehat{z}}(y_0)=\widehat{y}$. Then $P$ is homeomorphic to $B(c,r)\times \mathcal{T}$ via the following continuous bijective mappings:
$$
\Psi:\left\{
    \begin{array}{ll}
      B(c,r)\times\mathcal{T} & \longrightarrow P\\
      \textcolor{white}{0}    & \textcolor{white}{0}\\
      (y_0,\widehat{z})       & \longmapsto \varphi_{\widehat{z}}(y_0)
    \end{array}
\right.
\ \ , \ \ 
\Theta:\left\{
    \begin{array}{ll}
      P                       & \longrightarrow B(c,r)\times\mathcal{T}\\
      \textcolor{white}{0}    & \textcolor{white}{0}\\
      \widehat{y}             & \longmapsto (y_0,\widetilde{\pi}(\widehat{y}))
    \end{array}
\right.
$$
In particular, the Pesin box $P$ is a Borel subset of $\widehat{\mathbb{P}^k}$. The continuity of $\Psi$ and $\Theta$ can be checked by hand using the product topology on $\widehat{\mathbb{P}^k}$. This is implicit in {\cite{Bri01}}.

By Rokhlin's theorem {\cite[§3.1]{rokhlin1952fundamental}},   $\widehat{\mu}$ admits a decomposition on the fibers of $\pi_0$ :
\begin{equation*}
    \forall A\in\mathcal{B} (\widehat{\mathbb{P}^k})\ ,\ \widehat{\mu}(A) = \int_{\mathbb{P}^k}\mu_{\pi_0,c}(A\cap T_c)\ \mathrm{d}\mu(c),
\end{equation*}
where $\mu_{\pi_0,c}$ is the conditional measure of $\widehat{\mu}$ on the fiber $T_c:=\{\pi_0=c\}$. Actually, Briend proved the following result which specifies Rokhlin's formula on Pesin boxes:

\begin{thm}[Briend {\cite[Theorem 4.1]{Bri01}}]\label{thm:BriendforPesinBoxes}
For $\mu-$almost every $c\in\mathbb{P}^k$, there exist $r_c>0$ and $\rho_c>0$ such that $\widehat{\mu}$ is a product
$$\Theta_*\left(\widehat{\mu}|_P\right) = \mu|_{B(c,r)}\otimes\mu_{\pi_0,c}|_{\mathcal{T}}$$
on every Pesin box $(P,r,\rho,\mathcal{T},c)$ such that $r\leq r_c$, $\rho\geq\rho_c$ and $\displaystyle{\mathcal{T}\!\!\subset\!\!\{\widehat{x}\in\Lambda:\eta(\widehat{x})\!\geq\! r,\ \beta(\widehat{x})\!\leq\!\rho\}}$. 
\end{thm}

{\subsection{Construction of Pesin boxes}\label{sec:constructionofA}}

\begin{prop}\label{prop:constructionPesinBoxes}
    For $\widehat{\mu}-$almost every $\widehat{z}\in\Lambda$, there exists $(P,r,\rho,\mathcal{T},c)$ a Pesin box such that:
    \begin{enumerate}
        \item The point $c=z_0$ satisfies the content of Theorem \ref{thm:BriendforPesinBoxes}.
        \item $(r,\rho)$ is admissible in the sense of Definition \ref{defn:admissibleradius} for $\widehat{z}$, and we have $\widehat{z}\in\mathcal{T}$.
        \item $\mathcal{T}\subset\{\widehat{x}\in\Lambda:\eta(\widehat{x})\geq1/\rho,\ \beta(\widehat{x})\leq\rho\}$, $1/\rho\leq r_c$ and $\rho\geq\rho_c$, and $r=\frac{1}{\rho}$.
        \item For all $0< l\leq r$, if we denote $P_{c,l}:=\bigsqcup_{\widehat{x}\in\mathcal{T}}W^u(\widehat{x},l)\subset P$, then $(P_{c,l},l,\rho,\mathcal{T},c)$ is a Pesin box satisfying:
        \begin{equation}\label{eq:mu(P_c)=mu(boule)timesmu(K_C)>0}
            \widehat{\mu}(P_{c,l}) = \mu(B(c,l))\times \mu_{\pi_0,c}(\mathcal{T}) > 0.
        \end{equation}
    \end{enumerate}
    In particular, $W^u(\widehat{z},l)\subset P_{c,l}$ for each $0<l\leq r$. 
\end{prop}
\noindent\textbf{\underline{Proof:}} Let $Br$ be a set of full $\mu-$measure of points $c\in \mathrm{Supp}(\mu)$ that satisfy Theorem \ref{thm:BriendforPesinBoxes}. Let $\Lambda'\subset\Lambda\cap\pi_0^{-1}(Br)$ be a set of measure $1$ for $\widehat{\mu}$, such that the conditional measure $\mu_{\pi_0,z_0}$ exists and satisfies $\mu_{\pi_0,z_0}(\Lambda\cap T_{z_0})=1$ for every $\widehat{z}\in\Lambda'$. Let us fix now $\widehat{z}\in\Lambda'$ and denote $c:=z_0\in Br$. 

For $\rho$ such that $\eta(\widehat{z})\geq1/\rho$, $\beta(\widehat{z})\leq\rho$, $1/\rho\leq r_c$ and $\rho\geq \rho_c$ (recall that $r_c$ and $\rho_c$ are given by Theorem \ref{thm:BriendforPesinBoxes}), we define 
$$V:=\{\widehat{x}\in\Lambda:\eta(\widehat{x})\geq1/\rho,\ \beta(\widehat{x})\leq\rho\}\ni\widehat{z}.$$
By definition of $\Lambda'$, we have $\mu_{\pi_0,c}(V\cap T_c)>0$ if $\rho$ is large enough. Denote $r:=\frac{1}{\rho}$ and $\mathcal{T}:=V\cap T_c$. Observe that for every $l\in [0,r]$, the pair $(l,\rho)$ is admissible in the sense of Definition \ref{defn:admissibleradius} for each $\widehat{x}\in\mathcal{T}$. Thus $(P_{c,l},l,\rho,\mathcal{T},c)$ is a Pesin box with respect to Definition \ref{defn:pesinboxes}, where $P_{c,l}$ is defined by $P_{c,l}:=\bigsqcup_{\widehat{x}\in\mathcal{T}}W^u(\widehat{x},l)$. By Theorem \ref{thm:BriendforPesinBoxes} we have for each $l\in]0,r]$, $\widehat{\mu}(P_{c,l}) = \mu(B(c,l))\times \mu_{\pi_0,c}(V\cap T_c) > 0$. The conclusion then follows since $\widehat{\mu}(\Lambda')=1$.\qed\\

Let us now introduce the special Pesin box $\mathcal{A}$ that will be use to prove Theorem \ref{thm:TheoremB}. Let $(P,r,\rho,\mathcal{T},c)$ be a Pesin box, with $r=\frac{1}{\rho}$ and $\mathcal{T}\subset\{\widehat{z}\in\Lambda:\eta(\widehat{z})\geq1/\rho,\ \beta(\widehat{z})\leq\rho\}$, given by Proposition \ref{prop:constructionPesinBoxes}.

\begin{lemme}[{\cite[Lemma 2.2]{dup06}}]\label{lemma:slow} There exists $S\subset[1/4\rho, 1/2\rho]$ a Borel set of full $\mathrm{Leb}_{\reels}$-measure such that for every $s\in S$, there exists a measurable function $\gamma:\widehat{\mathbb{P}}^k\to]0,1]$ satisfying 
$$\forall\widehat{x}-\widehat{\mu}-a.e.,\ \forall n\in\mathbb{N},\ |d(x_{-n},c)-s|\geq \gamma(\widehat{x})e^{-n\varepsilon}.$$ 
Up to modify $\Lambda$, we can assume that this inequality holds for every $\widehat{x}\in\Lambda$.
\end{lemme}

Let us fix $\mathcal{R}\in S$ and define 
\begin{equation}\label{eq:realborelianA}
    \mathcal{A}:=P_{c,\mathcal{R}}=\bigsqcup_{\widehat{z}\in \mathcal{T}}W^u(\widehat{z},\mathcal{R})\subset \pi_0^{-1}(B(c,\mathcal{R})).
\end{equation}
$\widehat{\mu}(\mathcal{A})>0$ according to \eqref{eq:mu(P_c)=mu(boule)timesmu(K_C)>0}. This choice of $\mathcal{R}\in S$ given by Lemma \ref{lemma:slow} is used to prove Item \textit{3.} of Lemma \ref{lemma:6411}, and this item will be used to prove that the measure $(T^r\wedge dd^c|W_{\widehat{x}}^{r+1}|^2\wedge\cdots\wedge dd^c|W_{\widehat{x}}^{k}|^2)\circ \pi_0\lfloor_{\xi_{\widehat{x}}^u}$ does not vanish if $\widehat{x}\in\mathcal{A}\cap\Lambda$ (Proposition \ref{prop:62}).

\subsection{Extension of normal forms on Pesin boxes}

In this paragraph, we prove that the commutative diagram \eqref{eq:diagramnormalforms} can be extended to the ball $B(x_0, 2R)$ for each $\widehat{x} \in P \cap \Lambda$, provided that $(P, R, \rho, \mathcal{T}, c)$ is a Pesin box given by Proposition \ref{prop:constructionPesinBoxes} with $R \leq \frac{1}{2} \rho$. This applies in particular to the case $P = \mathcal{A}$ and $R = \mathcal{R}$. 

This extension property plays an important role in the article, as it allows us to define the measures $q_{\widehat{x}}$ in Section \ref{definitionofthemeasuresq_xunderH_1}, using the charts $(\sigma_{\widehat{x}})_{\widehat{x} \in \mathcal{A} \cap \Lambda}$ that will be extended to balls with uniform radius $\mathcal{R}$, see Proposition \ref{prop:extensionofnormalforms}. Proposition \ref{prop:extensionofnormalforms} is also used in the proof of Theorem \ref{thm:localfoliations}. Note that for this statement, we do not need any assumptions on the Lyapunov exponents, we just need the commutative diagram \eqref{eq:diagramnormalforms}.

\begin{prop}\label{prop:extensionofnormalforms} Let $\widehat{z}\in\Lambda$ and let $(P,R,\rho,\mathcal{T},c)$ be a Pesin box. Assume that $\widehat{z}\in\mathcal{T}$, $\frac{1}{\rho}\leq\eta(\widehat{z})$ and $R\in]0,\frac{1}{2\rho}]$. Then for each $\widehat{x}\in W^u(\widehat{z},R)\cap\Lambda$ the following points hold:
    \begin{enumerate}
        \item We have $B(x_0,2R)\subset B(z_0,2\eta(\widehat{z}))$.
        \item For each $n\geq0$, the map $f^{-n}_{\widehat{x}}:B(x_0,2\eta(\widehat{x}))\to\mathbb{P}^k$ extends holomorphically on $B(x_0,2\eta(\widehat{x}))\cup B(z_0,2\eta(\widehat{z}))$, the extension is unique thus it is also denoted $f^{-n}_{\widehat{x}}$.
        \item If $\eta(\widehat{x})<R$, the map $\sigma_{\widehat{x}}:B(x_0,2\eta(\widehat{x}))\to\mathbb{D}^k(\rho(\widehat{x}))$ extends holomorphically on 
        $B(x_0,2R)$, the extension is unique thus it is also denoted $\sigma_{\widehat{x}}$.
        \item There exists $\widetilde{\rho}(\widehat{x})>0$ such that $\sigma_{\widehat{x}}:B(x_0,2R)\to \mathbb{D}^k(\widetilde{\rho}(\widehat{x}))$
        and this map is injective, thus it defines holomorphic coordinates $\sigma_{\widehat{x}}=(Z_{\widehat{x}},W_{\widehat{x}})$ on $B(x_0,2R)$, where $Z_{\widehat{x}}=(Z^1_{\widehat{x}},\cdots,Z^r_{\widehat{x}})$ and $W_{\widehat{x}}=(W^{r+1}_{\widehat{x}},\cdots,W_{\widehat{x}}^k)$. 
        \item Up to increase $N(\widehat{x})$, for each $n\geq N(\widehat{x})$ we have:
        $$f^{-n}_{\widehat{x}}(B(x_0,2R))\subset B(x_{-n},\eta(\widehat{x}_{-n}))\ \mathrm{and}\ R_{n,\widehat{x}}\left(\mathbb{D}^k\left(\widetilde{\rho}(\widehat{x})\right)\right)\subset \mathbb{D}^k\left(\rho(\widehat{x}_{-n})\right),$$
        and the following diagram commutes:
        \begin{equation}\label{eq:extendeddiagramnormalforms}
        \xymatrix{
            B({x_{-n}},\eta(\widehat{x}_{-n})) \ar[d]_{ \sigma_{\widehat{x}_{-n}}} & &  B({x_0},2R) \ar[ll]_{f^{-n}_{\widehat{x}}} \ar[d]^{\sigma_{\widehat{x}}}  \\
            \mathbb{D}^k(\rho(\widehat{x}_{-n}))  & & \mathbb{D}^k(\widetilde{\rho}({\widehat{x}}))  \ar[ll]^{R_{n,\widehat{x}}}
            }
    \end{equation}
    \end{enumerate}
\end{prop}
\noindent\textbf{\underline{Proof:}}\\
\textit{1.} \& \textit{2.} Since $R\leq\frac{1}{2\rho}$ and $\frac{1}{\rho}\leq\eta(\widehat{z})$, $3R\leq\frac{3}{2\rho}\leq2\eta(\widehat{z})$, and since $x_0\in B(z_0,R)$, we deduce $B(x_0,2R)\subset B(z_0,2\eta(\widehat{z}))$. Moreover, according to Lemma \ref{lemma:f-ncoincidesontheunstablemanifold}, there exists a holomorphic mapping $g_n:B(x_0,2\eta(\widehat{x}))\cup B(z_0,2\eta(\widehat{z}))\to\mathbb{P}^k$ such that $g_n=f^{-n}_{\widehat{x}}$ on $B(x_0,2\eta(\widehat{x}))$ and $g_n=f^{-n}_{\widehat{z}}$ on $B(z_0,2\eta(\widehat{z}))$. In particular $g_n$ extends $f^{-n}_{\widehat{x}}$ on $B(x_0,2\eta(\widehat{x}))\cup B(z_0,2\eta(\widehat{z}))$.\\

\noindent\textit{3.} Let us denote $l_n:=\mathrm{Lip}(g_n|_{B(z_0,2\eta(\widehat{z}))})$, we have $l_n\leq \beta(\widehat{z})e^{-n(\lambda_k-\varepsilon)}$ by Theorem \ref{thm:normalforms}. For $n\geq0$, we have $B(z_{-n},2\eta(\widehat{z})l_n)\subset B(x_{-n},4\eta(\widehat{z})l_n)$, since $x_{-n}=g_n(x_0)$ and since $x_0\in B(z_0,2\eta(\widehat{z}))$. By Theorem \ref{thm:normalforms}, we also have $\eta(\widehat{x}_{-n})\geq e^{-n\varepsilon}\eta(\widehat{x})$ for each $n\geq 0$.
So we deduce that, if $N_1(\widehat{x})\geq N(\widehat{x})$ is large enough, for each $n\geq N_1(\widehat{x})$ we have
$g_{n}(B(z_0,2\eta(\widehat{z})))\subset B(x_{-n},\eta(\widehat{x}_{-n}))$. Recall that, according to Theorem \ref{thm:normalforms}, $R_{n,\widehat{x}}$ is invertible on $\cmplex^k$ with $R_{n,\widehat{x}}^{-1}$ being a polynomial map of $\cmplex^{k}$. Therefore, we can define for each $n\geq N_1(\widehat{x})$ a map  $\sigma_{n}:B(z_0,2\eta(\widehat{z}))\to\cmplex^k$ by:
$$\sigma_{n}(p):=R_{n,\widehat{x}}^{-1}\circ\sigma_{\widehat{x}_{-n}}\circ g_{n}(p),\ \forall p\in B(z_0,2\eta(\widehat{z})).$$
Let us fix $n_0\geq N_1(\widehat{x})$ an arbitrary integer. For each $p\in B(x_0,2\min\{\eta(\widehat{x}),R\})$ we have:
$$\sigma_{n_0}(p)=R_{n_0,\widehat{x}}^{-1}\circ R_{n_0,\widehat{x}}\circ\sigma_{\widehat{x}}(p)=\sigma_{\widehat{x}}(p),$$
by using that $g_n=f^{-n}_{\widehat{x}}$ on $B(x_0,2\eta(\widehat{x}))$ and that Diagram \eqref{eq:diagramnormalforms} is commutative. Then, if $\eta(\widehat{x})<R$, we observe that $\sigma_{n_0}$ extends $\sigma_{\widehat{x}}$ on $B(z_0,2\eta(\widehat{z}))$. It proves in particular that $\sigma_{\widehat{x}}$ extends (uniquely) as an injective holomorphic map on $B(x_0,2R)$.\\ 

\noindent\textit{4.} For each $n\geq N_1(\widehat{x})$, since $\sigma_{\widehat{x}_{-n}}\circ f^{-n}_{\widehat{x}}=R_{n,\widehat{x}}\circ\sigma_{\widehat{x}}$ on $B(x_0,2\eta(\widehat{x}))$ (by using again that Diagram \eqref{eq:diagramnormalforms} commutes), we also have
\begin{equation}\label{eq:relationcommutation}
    \sigma_{\widehat{x}_{-n}}\circ f^{-n}_{\widehat{x}}=R_{n,\widehat{x}}\circ\sigma_{n_0}\ \mathrm{on}\ B(x_0,2R)
\end{equation}
by analytic continuation if $\eta(\widehat{x})<R$. For $n\geq N_1(\widehat{x})$, in Item \textit{3.} we have established that $g_{n}(B(z_0,2\eta(\widehat{z})))\subset B(x_{-n},\eta(\widehat{x}_{-n}))$. Since $g_n$ extends $f^{-n}_{\widehat{x}}$ on $B(x_0,2R)$, and since $B(x_0,2R)\subset B(z_0,2\eta(\widehat{z}))$ by the first item, we have 
$$\sigma_{n_0}(B(x_0,2R))\subset R_{n_0,\widehat{x}}^{-1}\circ\sigma_{\widehat{x}_{-n_0}}(B(x_{-n_0},\eta(\widehat{x}_{-n_0})))\subset R_{n_0,\widehat{x}}^{-1}(\mathbb{D}^k(\rho(\widehat{x}_{-n_0}))),$$
since $\sigma_{\widehat{x}_{-n_0}}(B(x_{-n_0},\eta(\widehat{x}_{-n_0})))\subset\mathbb{D}^k(\rho(\widehat{x}_{-n_0}))$ by Theorem \ref{thm:normalforms}. Then taking $\widetilde{\rho}(\widehat{x})>0$ large enough we have $R_{n_0,\widehat{x}}^{-1}(\mathbb{D}^k(\rho(\widehat{x}_{-n_0})))\subset\mathbb{D}^k(\widetilde{\rho}(\widehat{x}))$, and thus we obtain that $\sigma_{n_0}(B(x_0,2R))\subset\mathbb{D}^k(\widetilde{\rho}(\widehat{x}))$.\\

\noindent\textit{5.} Recall that for $n\geq N_1(\widehat{x})$, we have proved that $g_{n}(B(z_0,2\eta(\widehat{z})))\subset B(x_{-n},\eta(\widehat{x}_{-n}))$ and that $g_n$ extends $f^{-n}_{\widehat{x}}$ on $B(x_0,2R)\subset B(z_0,2\eta(\widehat{z}))$, so we have $f_{\widehat{x}}^{-n}(B(x_0,2R))\subset B(x_{-n},\eta(\widehat{x}_{-n}))$. It remains to prove that \eqref{eq:extendeddiagramnormalforms} commutes for $n$ large enough. If $N_2(\widehat{x}) \geq N_1(\widehat{x})$ is large enough, one can infer that $R_{n,\widehat{x}}(\mathbb{D}^k(\widetilde{\rho}(\widehat{x}))) \subset \mathbb{D}^k(\rho(\widehat{x}_{-n}))$ for each $n \geq N_2(\widehat{x})$. Indeed, the coefficients of the polynomial maps defining $R_{n,\widehat{x}}$ are bounded, up to a constant, by $e^{-n(\lambda_k - \varepsilon)}$, according to Theorem \ref{thm:normalforms}. Thus, $R_{n,\widehat{x}}$ contracts $\mathbb{D}^k(\widetilde{\rho}(\widehat{x}))$ faster than $\rho(\widehat{x}_{-n})$ may decay to zero. The diagramme \eqref{eq:extendeddiagramnormalforms} then commutes for $n\geq N_2(\widehat{x})$, since we have \eqref{eq:relationcommutation}. 

{\subsection{Invariance of \texorpdfstring{$dW_{\widehat{x}}$}{TEXT} on unstable manifolds of {Pesin boxes}}\label{sec:InvarinceoddW_xonunstablemainifolds}}

We prove here Theorem \ref{thm:localfoliations}. Using Proposition \ref{prop:constructionPesinBoxes}, it implies in particular Theorem \ref{thm:Patching}. The assumptions of the theorem are satisfied for the Pesin box $(\mathcal{A},\mathcal{R},\rho,\mathcal{T},c)$.\\ 

\begin{thm}\label{thm:localfoliations}
    Assume that there exists $1\leq r\leq k-1$ such that {$\lambda_r>\lambda_{r+1}$}. Let $(P,R,\rho,\mathcal{T},c)$ be a Pesin box such that $\mathcal{T}\subset\{\widehat{z}\in\Lambda:\frac{1}{\rho}\leq\eta(\widehat{z})\}$ and $R\in]0,\frac{1}{2\rho}]$. Then for each $\widehat{z}\in\mathcal{T}$ and for each $\widehat{x},\widehat{y}\in W^u(\widehat{z},R)\cap\Lambda$ (e.g. $\widehat{y}:=\widehat{z}$) the following points hold:
    \begin{enumerate}
        \item The chart $\sigma_{\widehat{x}}$ defined on $B(x_0,2\eta(\widehat{x}))$ extends holomorphically into a chart $\sigma_{\widehat{x}}$ defined on $B(x_0,2R)$ defining holomorphic coordinates $\sigma_{\widehat{x}}=(Z_{\widehat{x}},W_{\widehat{x}})$ on $B(x_0,2R)$, where $Z_{\widehat{x}}=(Z^1_{\widehat{x}},\cdots,Z^r_{\widehat{x}})$ and $W_{\widehat{x}}=(W^{r+1}_{\widehat{x}},\cdots,W_{\widehat{x}}^k)$. The same property holds for $\widehat{y}$.
        \item $B(z_0,R)\subset\left(\displaystyle{\phantom{0}^{\phantom{0}^{\phantom{0}}}\!\!\!\!\!\!\!\!B}(x_0,2R)\cap B(y_0,2R)\right)$, thus $W_{\widehat{x}}$ and $W_{\widehat{y}}$ are defined on $B(z_0,R)$.
        \item There exists $C\in\mathrm{GL}_{k-r}({\mathcal{O}(B(z_0,R)}))$ (depending on $\widehat{x},\widehat{y}$) such that $dW_{\widehat{y}}=CdW_{\widehat{x}}$ on $B(z_0,R)$. Here $dW_{\widehat{x}}$ denotes the column vector $\left[dW_{\widehat{x}}^{r+1}\cdots dW_{\widehat{x}}^k\right]^{\mathrm{T}}$, and similarly for $dW_{\widehat{y}}$.
        \item If the Lyapunov exponents $\lambda_{r+1}\geq\cdots\geq\lambda_k$ are not resonant, then the lower triangular part and the diagonal of the matrix $C$ are constant functions. If moreover $\lambda_{r+1} = \cdots = \lambda_k$, then all the coefficients of $C$ are constant functions, hence $C \in \mathrm{GL}_{k-r}(\mathbb{C})$.
    \end{enumerate}
\end{thm}

The proof is divided in two parts. We give in a first part a proof of the theorem when the Lyapunov exponents are not resonant. We explain in a second part how one can adapt the proof when the Lyapunov exponents are resonant. We should also observe that the first item of the theorem is a direct consequence of Proposition \ref{prop:extensionofnormalforms}. 

{\subsubsection{Proof of Theorem {\ref{thm:localfoliations}} for non resonant Lyapunov exponents}}

The first item of Proposition \ref{prop:extensionofnormalforms} applied to $\widehat{x}$ and to $\widehat{y}$ implies:
\begin{equation}\label{eq:B_xetB_ysontdansB_z}
    B(x_0,2R)\subset B(z_0,2\eta(\widehat{z}))\ \mathrm{and}\ B(y_0,2R)\subset B(z_0,2\eta(\widehat{z})).
\end{equation}
Since $\widehat{x},\widehat{y}\in W^u(\widehat{z},R)$, we have also $B(z_0,R)\subset B(x_0,2R)$ and $B(z_0,R)\subset B(y_0,2R)$. Therefore we have:
\begin{equation}\label{eq:eta_xincludedinB(z_0,R)includedinBxinterBy}
     B(z_0, R) \subset \left(\displaystyle{\phantom{0}^{\phantom{0}^{\phantom{0}}}\!\!\!\!\!\!\!\!B}(x_0,2R)\cap B(y_0,2R)\right)\subset B(z_0,2\eta(\widehat{z})).
\end{equation}
According to Lemma \ref{lemma:f-ncoincidesontheunstablemanifold} (see also Proposition \ref{prop:extensionofnormalforms}), 
there exists a holomorphic map denoted $g_n:B(x_0,2\eta(\widehat{x}))\cup B(z_0,2\eta(\widehat{z}))\to\mathbb{P}^k$ such that $g_n=f^{-n}_{\widehat{x}}$ on $B(x_0,2\eta(\widehat{x}))$ and $g_n=f^{-n}_{\widehat{z}}$ on $B(z_0,2\eta(\widehat{z}))$, and similarly there exists a holomorphic map $h_n:B(y_0,2\eta(\widehat{y}))\cup B(z_0,2\eta(\widehat{z}))\to\mathbb{P}^k$ such that $h_n=f^{-n}_{\widehat{y}}$ on $B(y_0,2\eta(\widehat{x}))$ and $h_n=f^{-n}_{\widehat{z}}$ on $B(z_0,2\eta(\widehat{z}))$. In particular $g_n$ and $h_n$ coincide with $f^{-n}_{\widehat{z}}$ on $B(z_0,2\eta(\widehat{z}))$, and thus $f^{-n}_{\widehat{z}}$ extends $f^{-n}_{\widehat{x}}$ and $f^{-n}_{\widehat{y}}$ on $B(z_0,2\eta(\widehat{z}))$. So in what follows we can replace $f^{-n}_{\widehat{x}}$ and $f^{-n}_{\widehat{y}}$ by $f^{-n}_{\widehat{z}}$.

Let us denote $N_0:=\max\{N(\widehat{x}),N(\widehat{y})\}$. According to Proposition \ref{prop:extensionofnormalforms}, we have $f^{-n}_{\widehat{x}}(B(x_0,2R))\subset B(x_{-n},\eta(\widehat{x}_{-n}))$ (and similarly replacing $\widehat{x}$ by $\widehat{y}$) for each $n\geq N_0$. By using that $f^{-n}_{\widehat{z}}$ extends $f^{-n}_{\widehat{x}}$ and $f^{-n}_{\widehat{y}}$, and by using \eqref{eq:eta_xincludedinB(z_0,R)includedinBxinterBy} and the Diagram \eqref{eq:extendeddiagramnormalforms} we deduce that:
\begin{equation}\label{eq:f-n_z(eta_x)includedinBn}
    f^{-n}_{\widehat{z}}(B(z_0, R)) \subset B(x_{-n},\eta(\widehat{x}_{-n}))\cap B(y_{-n},\eta(\widehat{y}_{-n}))\ ,\ \forall n\geq N_0.
\end{equation}
For each $n\geq N_0$ we denote 
$$B_n:=B(x_{-n},2\eta(\widehat{x}_{-n}))\cap B(y_{-n},2\eta(\widehat{y}_{-n})).$$
Recall that $B(x_0,2R)\subset B(z_0,2\eta(\widehat{z}))$ by \eqref{eq:B_xetB_ysontdansB_z}, thus since $f^{-n}_{\widehat{x}}$ coincide with $f^{-n}_{\widehat{z}}$ on $B(z_0,2\eta(\widehat{z}))$, Proposition \ref{prop:extensionofnormalforms} ensures for any $n\geq N_0$:
$$\sigma_{\widehat{x}_{-n}}\circ f^{-n}_{\widehat{z}}=\sigma_{\widehat{x}_{-n}}\circ f^{-n}_{\widehat{x}}=R_{n,\widehat{x}}\circ\sigma_{\widehat{x}}\ \mathrm{on}\ B(x_0,2R).$$
Similarly we have $\sigma_{\widehat{y}_{-n}}\circ f^{-n}_{\widehat{z}}=R_{n,\widehat{y}}\circ\sigma_{\widehat{y}}$ on $B(y_0,2R)$ for each $n\geq N_0$. Using that the Lyapunov exponents are not resonant (see Theorem \ref{thm:normalforms}), an assumption which have not been used so far, we deduce from these equalities and from \eqref{eq:eta_xincludedinB(z_0,R)includedinBxinterBy} that the following equalities hold on $B(z_0,R)$ for any $n\geq N_0$, $j\in\{1,\cdots,r\}$ and $l\in\{r+1,\cdots,k\}$:
\begin{align}\label{eq:linearization_non_resonant_case}
    Z_{\widehat{x}_{-n}}^j\circ f^{-n}_{\widehat{z}} =\alpha_{n,\widehat{x}}^jZ_{\widehat{x}}^j\ \ & \mathrm{and}\ \ W_{\widehat{x}_{-n}}^l\circ f^{-n}_{\widehat{z}}=\beta_{n,\widehat{x}}^lW_{\widehat{x}}^l, \nonumber\\
    & \phantom{0} \nonumber\\
    W_{\widehat{y}_{-n}}^l\circ f^{-n}_{\widehat{z}} & =\beta_{n,\widehat{y}}^lW_{\widehat{y}}^l.
\end{align}
For each $n\geq N_0$, $j\in\{1,\cdots,r\}$ and $l,l'\in\{r+1,\cdots,k\}$ we define on $B_n$ the functions
\begin{align*}
    & C^{Z,lj}_n :=\frac{\partial W_{\widehat{y}_{-n}}^{l}}{\partial Z_{\widehat{x}_{-n}}^{j}} := \frac{\partial}{\partial z_j}\left[W^l\circ \sigma_{\widehat{y}_{-n}}\circ \sigma^{-1}_{\widehat{x}_{-n}}\right]\circ \sigma_{\widehat{x}_{-n}}\\
    & C^{W,ll'}_n :=\frac{\partial W_{\widehat{y}_{-n}}^l}{\partial W_{\widehat{x}_{-n}}^{l'}} = \frac{\partial}{\partial w_{l'}}\left[W^l\circ \sigma_{\widehat{y}_{-n}}\circ \sigma^{-1}_{\widehat{x}_{-n}}\right]\circ \sigma_{\widehat{x}_{-n}},
\end{align*}
where $W^l:\cmplex^k\to\cmplex$ is defined by $W^l(z_1,\cdots,z_r,w_{r+1},\cdots,w_k):=w_l$. Since we have the inclusion $f^{-n}_{\widehat{z}}(B(z_0,R))\subset B_n$ by \eqref{eq:f-n_z(eta_x)includedinBn},  we can also define the following map for any $l,l'\in\{r+1,\cdots,k\}$:
$$F_n^{ll'}:= C^{W,ll'}_{n}\circ f^{-n}_{\widehat{z}}:B(z_0,R)\to\cmplex.$$
We also denote $C_n^{ll'}:=C_n^{W,ll'}$ for simplicity.

\begin{lemme}\label{lemma:controlF_n}
    The following points hold for any $j\in\{1,\cdots,r\}$ and $l,l'\in\{r+1,\cdots,k\}$:
    \begin{enumerate}
        \item The functions $F_n^{ll'}|_{B(z_0,R)}$, $C^{Z,lj}_n|_{B_n}$, $C^{W,ll'}_n|_{B_n}$ are uniformly bounded by $2\beta(\widehat{y}_{-n})$ for $n\geq N_0$. 
        \item 
        There exists $N_1\geq N_0$ such that for all $n\geq N_1$:
        \begin{equation}\label{eq:controlfinsurF_n}
            \sup_{p,q\in B(z_0,R)}|F_n^{ll'}(p)-F_n^{ll'}(q)|\leq \frac{32R\beta(\widehat{x})\beta(\widehat{y})\beta(\widehat{z})}{\min\{\eta(\widehat{x}),\eta(\widehat{y})\}}e^{-n(\lambda_k-4\varepsilon)}.
        \end{equation}
    \end{enumerate}
\end{lemme}
\noindent\textbf{\underline{Proof:}} Denoting $H:=\sigma_{\widehat{y}_{-n}}\circ \sigma_{\widehat{x}_{-n}}^{-1}$ on $\sigma_{\widehat{x}_{-n}}(B_n),\ n\geq N_0$, we have $||dH||\leq \beta(\widehat{y}_{-n})\times 2$ (by Theorem \ref{sec:normalforms}). We can write $C_n^{ll'}$ as follows for $p\in B_n$:
$$C_n^{ll'}(p)=dW^l\circ d_{\sigma_{\widehat{x}_{-n}}(p)}H\cdot\vec{e}_{l'},$$
where $(\vec{e}_j)_{1 \leq j \leq k}$ denotes the canonical basis of $\mathbb{C}^k$. It implies $\sup_{B_n}|C_n^{ll'}|\leq 2\beta(\widehat{y}_{-n})$ for each $n\geq N_0$ (by the same arguments we prove $\sup_{B_n}|C^{Z,lj}_n|\leq 2\beta(\widehat{y}_{-n})$). It proves Item \textit{1}.

Let us prove now Item \textit{2}. By Theorem \ref{thm:normalforms} we have that $\mathrm{Lip}(f^{-n}_{\widehat{z}})\leq \beta(\widehat{z})e^{-n(\lambda_k-\varepsilon)}$ and that $\eta(\cdot)$ and $\beta(\cdot)$ are $\varepsilon-$tempered, thus since $\widehat{x},\widehat{y}\in W^u(\widehat{z},R)$, denoting $r_n:=\min\{\eta(\widehat{x}_{-n}),\eta(\widehat{y}_{-n})\}$, we deduce that there exists $N_1\geq N_0$ large enough such that:
$$d_{\mathbb{P}^k}(x_{-n},y_{-n})\leq \eta(\widehat{y}_{-n})\ \mathrm{and}\ f^{-n}_{\widehat{z}}(B(z_0,R))\subset B\left(x_{-n},\frac{r_n}{8\beta(\widehat{x}_{-n})}\right).$$ 
Let us fix $n\geq N_1$. From the Lipschitz estimate $||\sigma_{\widehat{x}_{-n}}(p)-\sigma_{\widehat{x}_{-n}}(q)||\geq \frac{1}{2} d_{\mathbb{P}^k}(p,q)$ for $p,q\in B(\widehat{x}_{-n},\eta(\widehat{x}_{-n}))$, see Theorem \ref{thm:normalforms}, we deduce that $\sigma_{\widehat{x}_{-n}}(B(x_{-n},r_n/2))\supset B_{\cmplex^k}(0,r_n/4)$ and that $\sigma_{\widehat{x}_{-n}}^{-1}$ is well defined on $B_{\cmplex^k}(0,r_n/4)$. We have $B(x_{-n},r_n/2)\subset B_n$ by using that $d_{\mathbb{P}^k}(x_{-n},y_{-n})\leq \eta(\widehat{y}_{-n})$, so we deduce that:
\begin{equation}\label{eq:distorsion}
    \sigma_{\widehat{x}_{-n}}^{-1}(B_{\cmplex^k}(0,r_n/4))\subset B(x_{-n},r_n/2)\subset B_n.
\end{equation}
By Cauchy inequality, we deduce from \eqref{eq:distorsion} and from the first item of the lemma that :
\begin{align*}
    \sup_{B_{\cmplex^k}(0,r_n/8)}||d(C_n^{ll'}\circ\sigma_{\widehat{x}_{-n}}^{-1})|| & \leq 8r_n^{-1}\sup_{B_n}|C_n^{ll'}|\leq 16r_n^{-1}\beta(\widehat{y}_{-n}).
\end{align*}
Moreover, since $f_{\widehat{z}}^{-n}(B(z_0,R))\subset B(x_{-n},r_n/8\beta(\widehat{x}_{-n}))$ and since $||d\sigma_{\widehat{x}_{-n}}||\leq\beta(\widehat{x}_{-n})$, we deduce that $\sigma_{\widehat{x}_{-n}}(f_{\widehat{z}}^{-n}(B(z_0,R)))\subset B_{\cmplex^k}(0,r_n/8)$. So, for each $p,q\in B(z_0,R)$, we have $|F_n^{ll'}(p)-F_n^{ll'}(q)| \leq \sup_{B_{\cmplex^k}(0,r_n/8)}||d(C_n^{ll'}\circ\sigma_{\widehat{x}_{-n}}^{-1})||\times ||\sigma_{\widehat{x}_{-n}}\circ f^{-n}_{\widehat{z}}(p)-\sigma_{\widehat{x}_{-n}}\circ f^{-n}_{\widehat{z}}(q)||$ and thus 
\begin{align*}
    |F_n^{ll'}(p)-F_n^{ll'}(q)| & \leq 16r_n^{-1}\beta(\widehat{y}_{-n})\beta(\widehat{x}_{-n})\times d_{\mathbb{P}^k}(f^{-n}_{\widehat{z}}(p),f^{-n}_{\widehat{z}}(q)).
\end{align*}
The result follows since $\mathrm{Lip}(f^{-n}_{\widehat{z}})\leq \beta(\widehat{z})e^{-n(\lambda_k-\varepsilon)}$ and since $\eta(\cdot)$ and $\beta(\cdot)$ are $\varepsilon-$tempered.\qed\\

Let us now express the column matrix of $1-$forms $dW_{\widehat{y}}=\left[dW^{r+1}_{\widehat{y}},\cdots,dW_{\widehat{y}}^k\right]^{\mathrm{T}}$ in the system of coordinates $(Z_{\widehat{x}},W_{\widehat{x}})$ on $B(z_0,R)$. To do so, we observe that for $n\geq N_0$ we have for each $l\in\{r+1,\cdots,k\}$:
$$dW_{\widehat{y}_{-n}}^l=\sum_{j=1}^{r}C_n^{Z,lj}dZ_{\widehat{x}_{-n}}^{j} + \sum_{l'=r+1}^{k}C^{W,ll'}_ndW_{\widehat{x}_{-n}}^{l'}\ \mathrm{on}\ B_n.$$
We can pull-back this equality by $f^{-n}_{\widehat{z}}:B(z_0,R)\to f^{-n}_{\widehat{z}}(B(z_0,R))$ since by \eqref{eq:f-n_z(eta_x)includedinBn} we have $f^{-n}_{\widehat{z}}(B(z_0,R))\subset B_n$. So, by using \eqref{eq:linearization_non_resonant_case} we have that for each $n\geq N_0$ and $l\in\{r+1,\cdots,k\}$:
\begin{equation}\label{eq:auniveaudes1formes}
    dW_{\widehat{y}}^l = \sum_{j=1}^{r}\frac{\alpha_{n,\widehat{x}}^j}{\beta_{n,\widehat{y}}^l}(C_n^{Z,lj}\circ f^{-n}_{\widehat{z}})dZ_{\widehat{x}}^j + \sum_{l'=r+1}^{k}\frac{\beta_{n,\widehat{x}}^{l'}}{\beta_{n,\widehat{y}}^{l}}(C_n^{W,ll'}\circ f^{-n}_{\widehat{z}})dW_{\widehat{x}}^{l'}\ \mathrm{on}\ B(z_0,R).
\end{equation}
By Theorem \ref{thm:normalforms} and by Lemma \ref{lemma:controlF_n} we observe that for each $n\geq N_0$, $j\in\{1,\cdots,r\}$ and $l\in\{r+1,\cdots,k\}$:
\begin{equation}\label{eq:controlsurC_n^Z}
    \left|\frac{\alpha_{n,\widehat{x}}^j}{\beta_{n,\widehat{y}}^l}\right|\sup_{B(z_0,R)}\left|C_{n}^{Z,lj}\circ f^{-n}_{\widehat{z}}\right|\leq 2\beta(\widehat{y}) e^{-n(\lambda_j-\lambda_l-3\varepsilon)}.
\end{equation}
Note that in \eqref{eq:auniveaudes1formes} we have expressed $dW_{\widehat{x}}^l$ in the basis $(dZ_{{\widehat{x}}}^1, \dots, dW_{\widehat{x}}^k)$, thus the coefficients $\frac{\alpha_{n,\widehat{x}}^j}{\beta_{n,\widehat{y}}^l}(C_n^{Z,lj}\circ f^{-n}_{\widehat{z}})$ and $\frac{\beta_{n,\widehat{x}}^{l'}}{\beta_{n,\widehat{y}}^{l}}(C_n^{W,ll'}\circ f^{-n}_{\widehat{z}})$ must be the coefficients in this basis. They are therefore elements of $\mathcal{O}(B(z_0, R))$ that are independent of $n$. Let us denote those associated with the $W-$coordinates as follows:
\begin{equation}\label{eq:b_x/b_yF_n(p)=dW_y+o(1)}
     C_{ll'}(p):=\frac{\beta_{n,\widehat{x}}^{l'}}{\beta_{n,\widehat{y}}^{l}}(C_n^{W,ll'}\circ f^{-n}_{\widehat{z}})=\frac{\beta_{n,\widehat{x}}^{l'}}{\beta_{n,\widehat{y}}^{l}}F_n^{ll'}(p),\ p\in B(z_0,R),\ \mathrm{does\ not\ depend\ on}\ n\geq N_0. 
\end{equation}

\begin{lemme}\label{lemma:laconvergenceuniformede_un|Fn|2}
    Assume that ${\lambda_{r}>\lambda_{r+1}}$ and that $\lambda_{r+1}\geq\cdots\geq\lambda_k$ are not resonant. Then for every $j\in\{1,\cdots,r\}$ and every $l,l'\in\{r+1,\cdots,k\}$ the following points hold:
    \begin{enumerate}
        \item The functions $\frac{\alpha_{n,\widehat{x}}^j}{\beta_{n,\widehat{y}}^l} (C_n^{Z, lj} \circ f^{-n}_{\widehat{z}})$, $n \geq N_0$, are identically null on $B(z_0, R)$.
        \item If $\lambda_{l'}\geq \lambda_l$, then $C_{ll'}\in\cmplex^*$ is constant on $B(z_0,R)$.
    \end{enumerate} 
    In particular, the coefficients $C_{ll'}, \ r+1 \leq l' \leq l \leq k$, are constant. Similarly, if $\lambda_{r+1} = \cdots = \lambda_k$, then all the coefficients $C_{ll'}, \ r+1 \leq l, l' \leq k$, are constant.
\end{lemme}
\noindent\textbf{\underline{Proof:}} The first item is a direct consequence of \eqref{eq:controlsurC_n^Z} since $\lambda_j>\lambda_l$ by assumption. Let us prove the second one. The fact that $\lambda_{l'}\geq \lambda_l$ implies $\left|\frac{\beta_{n,\widehat{x}}^{l'}}{\beta_{n,\widehat{y}}^l}\right|\leq e^{2n\varepsilon}$. So, using the second item of Lemma \ref{lemma:controlF_n} we deduce that, for any $n\geq N_1$ and for any $(p,q)\in B(z_0,R)^2$, we have
\begin{equation}\label{eq:controlintermediairedeF_n(p)-F_n(q)}
    \left|\frac{\beta_{n,\widehat{x}}^{l'}}{\beta_{n,\widehat{y}}^l}F_n^{ll'}(p)-\frac{\beta_{n,\widehat{x}}^{l'}}{\beta_{n,\widehat{y}}^l}F_n^{ll'}(q)\right|\leq \frac{32R\beta(\widehat{x})\beta(\widehat{y})\beta(\widehat{z})}{\min\{\eta(\widehat{x}),\eta(\widehat{y})\}}e^{-n(\lambda_k-6\varepsilon)}.
\end{equation}
For all $n\geq N_1$, since $C_{ll'}(z_0)=\frac{\beta_{n,\widehat{x}}^{l'}}{\beta_{n,\widehat{y}}^{l}}F_n^{ll'}(z_0)$ by \eqref{eq:b_x/b_yF_n(p)=dW_y+o(1)}, we deduce using \eqref{eq:controlintermediairedeF_n(p)-F_n(q)} that
$$\sup_{p\in B(z_0,R)}\left|\frac{\beta_{n,\widehat{x}}^{l'}}{\beta_{n,\widehat{y}}^l}F_n^{ll'}(p)-C_{ll'}(z_0)\right|\leq \frac{32R\beta(\widehat{x})\beta(\widehat{y})\beta(\widehat{z})}{\min\{\eta(\widehat{x}),\eta(\widehat{y})\}}e^{-n(\lambda_k-6\varepsilon)}.$$
It shows that $\left(\frac{\beta_{n,\widehat{x}}^{l'}}{\beta_{n,\widehat{y}}^{l}} F_n^{ll'}\right)_{n\geq N_0}$ converges uniformly on $B(z_0,R)$ to the constant $C_{ll'}=C_{ll'}(z_0)$.\qed\\

We are now able to finish the proof of Theorem \ref{thm:localfoliations}. Let us assume now that ${\lambda_r > \lambda_{r+1}}$ (and that $\lambda_{r+1}\geq\cdots\geq\lambda_k$ are not resonant). Then, going back to Equation \eqref{eq:auniveaudes1formes} and using the first item of Lemma \ref{lemma:laconvergenceuniformede_un|Fn|2}, and using the fact that $\frac{\beta_{n,\widehat{x}}^{l'}}{\beta_{n,\widehat{y}}^{l}} (C_n^{W, ll'} \circ f^{-n}_{\widehat{z}})(p) = C_{ll'}(p)$ is independent of $n$ by \eqref{eq:b_x/b_yF_n(p)=dW_y+o(1)}, we deduce that for all $p \in B(z_0, R)$ and for all $l \in \{r+1, \dots, k\}$ (where $n \geq N_0$):
\begin{align*}
    d_pW_{\widehat{y}}^l & =\sum_{j=1}^{r}\frac{\alpha_{n,\widehat{x}}^j}{\beta_{n,\widehat{y}}^l}(C_n^{Z,lj}\circ f^{-n}_{\widehat{z}})(p)d_pZ_{\widehat{x}}^j + \sum_{l'=r+1}^{k}\frac{\beta_{n,\widehat{x}}^{l'}}{\beta_{n,\widehat{y}}^{l}}(C_n^{W,ll'}\circ f^{-n}_{\widehat{z}})(p)d_pW_{\widehat{x}}^{l'}\\
    & =0+\sum_{l'=r+1}^{k}C_{ll'}(p)d_pW_{\widehat{x}}^{l'}.
\end{align*}
By defining the matrix $C := (C_{ll'})_{r+1 \leq l, l' \leq k}$, we observe that we have obtained
\begin{equation*}
        \begin{bmatrix}
        dW_{\widehat{y}}^{r+1}\\
        \vdots\\
        dW_{\widehat{y}}^{k} 
    \end{bmatrix}
    =C
    \begin{bmatrix}
        dW_{\widehat{x}}^{r+1}\\
        \vdots\\
        dW_{\widehat{x}}^{k}  
    \end{bmatrix}\ \mathrm{on}\ B(z_0,R).
    \end{equation*}
The matrix $C$ belongs to $\mathrm{GL}_{k-r}(\mathcal{O}(B(z_0, R)))$ since $(dZ_{\widehat{y}}^1, \dots, dW_{\widehat{y}}^k)$ and $(dZ_{\widehat{x}}^1, \dots, dW_{\widehat{x}}^k)$ are at each $p \in B(z_0, R)$ two bases of $T_p \mathbb{P}^k$. Lemma \ref{lemma:laconvergenceuniformede_un|Fn|2} also implies that the coefficients of the lower triangular part and the diagonal of this matrix are constant functions. This lemma also implies that, in the case $\lambda_{r+1} = \cdots = \lambda_k$, all the coefficients of the matrix $C$ are constant functions, so $C \in \mathrm{GL}_{k-r}(\mathbb{C})$. We finally obtain the third and fourth points of Theorem \ref{thm:localfoliations}.

\subsubsection{Proof of Theorem {\ref{thm:localfoliations}} for resonant Lyapunov exponents}\label{sec:remarkontheprooffortheresonantcase} 

To explain how the arguments adapt to the resonant case, we distinguish between two cases: the first case, where the exponents could be resonant, but the exponents $\lambda_{r+1} \geq \cdots \geq \lambda_k$ are not resonant, and the second case, where these $k-r$ smaller exponents may also be resonant.\\

\textbullet\textit{First case:} Assume that {$\lambda_r > \lambda_{r+1}$ and that $\lambda_{r+1}\geq\cdots\geq\lambda_k$ are not resonant}. When the Lyapunov exponents $\lambda_j$, $j\in\{1,\cdots,r\}$, satisfy some resonance relations of the form $\sum_{i=j+1}^kc_i\lambda_i=\lambda_j$, with $(c_i)_{i=j+1,\cdots,k}\in\mathbb{N}^{k-j}$ of length $\geq2$, the first $r$ components of the map $R_{n, \widehat{x}}$ may fail to be linear, see \eqref{eq:Rnxwithnoresonancesforlowerlyapunovexponents} for a precise expression of $R_{n, \widehat{x}}$ in this case.

In this situation, we then do not have exactly \eqref{eq:linearization_non_resonant_case} and thus we do not have exactly Equation \eqref{eq:auniveaudes1formes}, but these equations can be modified as follows. First in \eqref{eq:linearization_non_resonant_case} the only term which is affected is $Z_{\widehat{x}_{-n}}^j\circ f^{-n}_{\widehat{z}}$ which becomes on $B(z_0,R)$
\begin{align*}
    Z_{\widehat{x}_{-n}}^j\circ f^{-n}_{\widehat{z}} & = \alpha_{n,\widehat{x}}^{j}Z_{\widehat{x}}^{j}+\sum_{m=1}^{m_j}\gamma_{n,m,\widehat{x}}^{j}(Z_{\widehat{x}}^{j+1},\cdots,W_{\widehat{x}}^k)^{C_{j,m}}.
\end{align*}
Second we can then observe that \eqref{eq:auniveaudes1formes} becomes on $B(z_0,R)$
\begin{align}\label{eq:dw_y^ldanslecasresonant}
    dW_{\widehat{y}}^l & =\sum_{j=1}^{r}\frac{1}{\beta_{n,\widehat{y}}^l}\left(C_n^{Z,lj}\circ f^{-n}_{\widehat{z}}\right)\times\left(\alpha_{n,\widehat{x}}^{j}dZ_{\widehat{x}}^{j} + \sum_{m=1}^{m_j} \gamma_{n,m,\widehat{x}}^{j}d(Z_{\widehat{x}}^{j+1},\cdots,W_{\widehat{x}}^k)^{C_{j,m}}\right) \nonumber \\
    & + \sum_{l'=r+1}^{k}\frac{\beta_{n,\widehat{x}}^{l'}}{\beta_{n,\widehat{y}}^{l}}(C_n^{W,ll'}\circ f^{-n}_{\widehat{z}})dW_{\widehat{x}}^{l'}.
\end{align}
From this expression, using that $|C_n^{Z,lj}|_{B_n}|\leq 2\beta(\widehat{y})e^{n\varepsilon}$, and that $|\gamma_{n,m,\widehat{x}}^j|\leq \Gamma(\widehat{x})e^{-n(\lambda_r - \varepsilon)}$ (see Theorem \ref{thm:normalforms}), along with the fact that $\lambda_r > \lambda_{r+1}$, we conclude that the first sum on the right-hand side of this equality must vanish.\\

Concerning the second sum appearing in the second member of \eqref{eq:dw_y^ldanslecasresonant}, each coefficient in the sum must also be a function independent of $n$, which can be written as $C_{ll'} = \frac{\beta_{n,\widehat{x}}^{l'}}{\beta_{n,\widehat{y}}^{l}} (C_n^{W, ll'} \circ f^{-n}_{\widehat{z}}) = \frac{\beta_{n,\widehat{x}}^{l'}}{\beta_{n,\widehat{y}}^{l}} F_n^{ll'}$ as in \eqref{eq:b_x/b_yF_n(p)=dW_y+o(1)}. Thus, we can observe that we still have $dW_{\widehat{y}}^l = \sum_{l' = r+1}^k C_{ll'} dW_{\widehat{x}}^{l'}$, and by defining $C := (C_{ll'})_{r+1 \leq l, l' \leq k}$, we obtain $dW_{\widehat{y}} = C dW_{\widehat{x}}$. We can also observe that the arguments of Lemma \ref{lemma:laconvergenceuniformede_un|Fn|2} remain unchanged. Therefore, we can also conclude that the lower part and the diagonal of the matrix $C$ consist only of constant functions. Moreover, in the case where $\lambda_{r+1} = \cdots = \lambda_k$, all the coefficients of $C$ are constant functions. So, we can conclude that Theorem \ref{thm:localfoliations} also holds in this first case.\\

\textbullet\textit{Second case:} Assume now that $\lambda_r > \lambda_{r+1}$, but that $\lambda_{r+1} \geq \cdots \geq \lambda_k$ could be resonant. As stated below Theorem \ref{thm:normalforms}, we can adapt the statement of the theorem with $R_{n,\widehat{x}}$ having a similar form to the one given by \eqref{eq:Rnxwithnoresonancesforlowerlyapunovexponents}, i.e., we can write $R_{n,\widehat{x}}$ in the following form, where only the last component remains linear
\begin{equation}\label{eq:formepurementraisonnante}
    R_{n,\widehat{x}}=\left(\alpha_{n,\widehat{x}}^{1}z_1 + P_{1,n}(z_2,\cdots,w_k),\cdots,\beta_{n,\widehat{x}}^{k-1}w_{k-1} + P_{k-1,n}(w_k),\beta_{n,\widehat{x}}^{k}w_{k}\right).
\end{equation}
From this expression, we can observe that the last component is linear. So, one can prove similarly to the previous case that $dW^k_{\widehat{y}} = \sum_{l'=r+1}^k C_{kl'} dW_{\widehat{x}}^{l'}$, for some functions $C_{kl'} \in \mathcal{O}(B(z_0,R))$. Then, from the expression \eqref{eq:formepurementraisonnante}, one can prove with similar arguments as above that $dW_{\widehat{y}}^{k-1}$ is a combination of $dW_{\widehat{y}}^k$ and the forms $dW_{\widehat{x}}^l, \ r+1 \leq l \leq k$. Indeed, from \eqref{eq:formepurementraisonnante} applied to $\widehat{y}$, we have
$$W^{k-1}_{\widehat{y}_{-n}} \circ f^{-n}_{\widehat{z}} = \beta_{n,\widehat{y}}^{k-1} W_{\widehat{y}}^{k-1} + \gamma_{n,\widehat{y}}^{k-1} (W_{\widehat{y}}^k)^{c_{k-1}},$$
where $c_{k-1} \in \mathbb{N}_{\geq 2}$, and where $|\gamma_{n,\widehat{y}}^{k-1}| \leq \Gamma(\widehat{y}) e^{-n(\lambda_{k-1} - \varepsilon)}$. So, instead of \eqref{eq:auniveaudes1formes}, we will have
\begin{align*}
    dW_{\widehat{y}}^{k-1} + c_{k-1}\frac{\gamma^{k-1}_{n,\widehat{y}}}{\beta_{n,\widehat{y}}^{k-1} } (W_{\widehat{y}}^k)^{c_{k-1}-1} dW_{\widehat{y}}^k & = \frac{1}{\beta_{n,\widehat{y}}^{k-1}} (\text{combination of the}\ dZ_{\widehat{x}}^j,\ 1\leq j\leq r) \\
    & + \frac{1}{\beta_{n,\widehat{y}}^{k-1}} (\text{combination of the}\ dW_{\widehat{x}}^l,\ r+1\leq l\leq k).
\end{align*}
Similarly to what we observed in previous cases, the combination of the $dZ_{\widehat{x}}^j$ must vanish since it involves stronger contractions, leaving only the combination of the $dW_{\widehat{x}}^{l}$. Then, since $dW_{\widehat{y}}^k$ is itself a combination of the $dW_{\widehat{x}}^l, \ r+1 \leq l \leq k$, we deduce that $dW_{\widehat{y}}^{k-1} = \sum_{l'=r+1}^k C_{k-1,l'} dW_{\widehat{x}}^{l'}$, for some functions $C_{k-1,l'} \in \mathcal{O}(B(z_0, R))$.

By repeating these arguments inductively, going step by step up to $dW_{\widehat{y}}^{r+1}$, we can conclude that each form $dW_{\widehat{y}}^l, \ r+1 \leq l \leq k,$ is a combination of the forms $dW_{\widehat{x}}^l, \ r+1 \leq l \leq k$, i.e., there exists a matrix $C \in \mathrm{GL}_{k-r}(\mathcal{O}(B(z_0, R)))$ such that $dW_{\widehat{y}} = C dW_{\widehat{x}}$. Thus, we can conclude that Theorem \ref{thm:localfoliations} also holds in this second case. The only difference with the first case is that this time we cannot say whether certain coefficients of $C$ are constant.

\section{A decreasing partition \texorpdfstring{$\xi^u$}{TEXT} of the natural extension}\label{sec:partitionetaledrappier}

We construct the measurable partition $\xi^u$ that appears in Theorem \ref{thm:TheoremB}. We refer to {\cite{bedlyusmil1993, Cantat2001, dup06, rokhlin1952fundamental, Rokhlin1967LECTURESOT}} for accounts on the theory of measurable partitions and entropy. 
We mention that in \cite{bedlyusmil1993} (resp. in \cite{Cantat2001}) the theory of measurable partitions is used to prove the uniqueness of the measure of maximal entropy for polynomial automorphisms of $\cmplex^2$ (resp. for automorphisms of K3 surfaces). 

The arguments are borrowed from  {\cite{dup06}}. We use the partition $\mathcal{P}$ and the Pesin box $\mathcal{A}=\bigsqcup_{\widehat{z}\in\mathcal{T}}W^u(\widehat{z},\mathcal{R})$ respectively defined in \Dciteprop{\ref{prop:partitionPforinversebranchesI}} and in Sec. {\ref{sec:constructionofA}}, Equ. \eqref{eq:realborelianA}. We define the following measurable partitions of $\widehat{\mathbb{P}^k}$
$$ \xi:=\pi_0^{-1}(\mathcal{P})\bigvee\{\mathcal{A},\mathcal{A}^c\} \ \textrm { and } \ \xi^u:=\bigvee_{p\in\mathbb{N}}\widehat{f}^{p}(\xi).$$

\begin{lemme}\label{lemma:6411} The partition $\xi^u$  satisfies the following properties:
\begin{enumerate}
    \item {For all $\widehat{x}$ the elements of $\xi_{\widehat{x}}^u$ and $\widehat{x}$ have the same $\mathcal{P}-$address.}
    \item {$\forall\widehat{z}\in\mathcal{T},\forall\widehat{x}\in W^u(\widehat{z},\mathcal{R})$, $W^u(\widehat{z},\mathcal{R}) \supset \xi_{\widehat{x}}^u$.} 
    \item {There exists a measurable function $\eta'$ such that for every $\widehat{x}\in\Lambda$:} 
    {
    \begin{enumerate}
        \item $0<\eta'(\widehat{x})\leq\eta(\widehat{x})\leq 1$.
        \item $\forall n\geq0$, $f^{-n}_{\widehat{x}}(B(x_0,\eta'(\widehat{x})))\subset B(c,\mathcal{R})$ or $f^{-n}_{\widehat{x}}(B(x_0,\eta'(\widehat{x})))\subset B(c,\mathcal{R})^c$.
        \item $W^u(\widehat{x},\eta'(\widehat{x}))\subset\xi^u_{\widehat{x}}$.
    \end{enumerate}}
\end{enumerate}
\end{lemme}
\noindent\textbf{\underline{Proof :}} \
\begin{enumerate}
    \item[\textit{1.}] The $\mathcal{P}-$address is constant on the atoms of $\vee_p\widehat{f}^p(\pi_0^{-1}(\mathcal{P}))$, hence on those of $\xi^u$.
    \item[\textit{2.}] Let $\widehat{x}'\in\xi^u_{\widehat{x}}$, by the previous point $\widehat{x}'$ and $\widehat{x}$ have the same $\mathcal{P}-$address. One has $\widehat{x} \in \mathcal{A}$ by definition of $\mathcal{A}$, hence $\xi^u_{\widehat{x}}\subset \mathcal{A}$ because $\xi^u$ is thinner than $\{\mathcal{A},\mathcal{A}^c\}$. Thus $\widehat{x}'$ belongs to $\mathcal{A}$ and there is $\widehat{z}'\in\mathcal{T}$ such that $\widehat{x}'\in W^u(\widehat{z}',\mathcal{R})$. Since the $\mathcal{P}-$address is constant on unstable manifolds, we deduce that the $\mathcal{P}-$address of $\widehat{z}'$ is the one of $\widehat{x}'$, that is the one of $\widehat{x}$. Then, since $z_0'=c=z_0$, we can apply the point \textit{2.} of Lemma \ref{lemma:f-ncoincidesontheunstablemanifold} (applied with $\widehat{w}:=\widehat{z}'$) to conclude that $\widehat{z}'=\widehat{z}$ and thus $\widehat{x}'\in W^u(\widehat{z},\mathcal{R})$. 
    \item[\textit{3.}] Let $\eta'(\widehat{x}) := \min\left\{\mathcal{R},\eta(\widehat{x}),\frac{\gamma(\widehat{x})}{2\beta(\widehat{x})}\right\}$, 
    where $\gamma(\widehat{x})$ is the function given by Lemma \ref{lemma:slow} associated to the element $\mathcal{R}$ of $S$. This choice of $\mathcal{R}$ ensures that the diameter of $f^{-n}_{\widehat{x}}(B(x_0,\eta(\widehat{x})))$ decreases exponentially faster as $n\to+\infty$ than the distance between $x_{-n}$ and the boundary of the ball $B(c,\mathcal{R})$. It gives the point \textit{(b)}. We refer to {\cite[Lemma 4.2]{dup06}} for details and for a proof of \textit{(b)}\ $\Longrightarrow$\textit{(c)}.\qed\\
\end{enumerate}

The other important properties of the partition $\xi^u$ are gathered in the following theorem. Given $\zeta$ a measurable partition, we denote $\mathcal{M}(\zeta)$ the completion with respect to $\widehat{\mu}$ of the $\sigma-$algebra generated by the atoms of $\zeta$. Let $\mathcal{M}$ be the completion with respect to $\widehat{\mu}$ of the $\sigma-$algebra of Borel sets of $\widehat{\mathbb{P}^k}$. Denotes $\bigcup\xi^u$ the union of the atoms of $\xi^u$.

\newpage

\begin{thm}\label{thm:propertiesofeta}
The following properties hold for the partition $\xi^u$ :
\begin{enumerate}
    \item $\xi^u$ is decreasing : $\forall n\geq0$, $\forall\widehat{x}\in\bigcup\xi^u$, $({\widehat{f}}^{-n}\xi^u)_{\widehat{x}}\subset\xi^u_{\widehat{x}}$.
    \item $\pi_0$ is injective on $\xi^u_{\widehat{x}}$ for all $\widehat{x}\in\bigcup\xi^u$.
    \item For all $\widehat{x}\in\bigcup\xi^u$, $f^n$ is injective on $\pi_0\left(({\widehat{f}}^{-n}\xi^u)_{\widehat{x}}\right)$ for every $n\geq0$.
    \item $\forall n\geq0,\forall\widehat{x}\in\Lambda$, $\xi^u_{\widehat{x}}$ is a countable union of atoms of  ${\widehat{f}}^{-n}\xi^u$.
    \item $\bigvee_{n\geq 0}\mathcal{M}({\widehat{f}}^{-n}\xi^u)$ coincide with $\mathcal{M}$. 
   \item For every $n \geq 0$, $\int_{\widehat{\mathbb{P}^k}} -\LLog\ \mu_{\widehat{x}}({\widehat{f}}^{-n}\xi^u)_{\widehat{x}}\ \mathrm{d}\widehat{\mu}(x) =  \LLog\ d^{kn}$.     
\end{enumerate}
\end{thm}

The proof follows classical arguments, developed for instance in {\cite{dup06}}. Let us explain the last item. By definition $\int_{\widehat{\mathbb{P}^k}}-\LLog\ \mu_{\widehat{x}}(({\widehat{f}}^{-n}\xi^u)_{\widehat{x}})\ \mathrm{d}\widehat{\mu}(\widehat{x})$ is equal to the conditional entropy $H({\widehat{f}}^{-n}\xi^u|\xi^u)$. Moreover, following {\cite[Section 4.3]{dup06}} one can prove that  $H({\widehat{f}}^{-n}\xi^u|\xi^u)$ is equal to the relative entropy $h_{\widehat{\mu}}({\widehat{f}}^{-n},\xi_n)$, where $\xi_n:={\widehat{f}}^{-1}\xi\vee\cdots\vee {\widehat{f}}^{-n}\xi$. But  $\xi_n$ is a generator of  finite entropy for ${\widehat{f}}^{-n}$ {\cite[Proposition 4.1]{dup06}}, hence $h_{\widehat{\mu}}({\widehat{f}}^{-n},\xi_n)=h_{\widehat{\mu}}({\widehat{f}}^{-n})=n h_{\widehat{\mu}}({\widehat{f}})$ by Kolmogorov-Sinaï theorem {\cite[\textsection 9]{Rokhlin1967LECTURESOT}}. One finally obtains the last item {from $h_{\widehat{\mu}}({\widehat{f}}) = h_{\mu}(f)=\LLog\ d^k$.}

From Item \textit{3.} we deduce the following lemma used in the proof of Proposition \ref{lemma:66}.

\begin{lemme}\label{lemma:f-ncircfn=id} 
Let $\widehat{x}\in\mathcal{A}\cap\Lambda$ and $n\geq0$ be fixed. Let $\widehat{z}\in\mathcal{T}$ such that $\widehat{x}\in W^u(\widehat{z},\mathcal{R})$.
\begin{enumerate}
    \item $\pi_0(\xi_{\widehat{x}}^u)\subset B(z_0,2\eta(\widehat{z}))$ and $f^{-n}_{\widehat{z}}[\pi_0(\xi^u_{\widehat{x}})]= \pi_0[({\widehat{f}}^{-n}\xi^u)_{\widehat{x}_{-n}}]$.
    \item The map $f^{-n}_{\widehat{z}}$ extends the map $f^{-n}_{\widehat{x}}$ on $B(z_0,2\eta(\widehat{z}))\supset f^n(\pi_0[({\widehat{f}}^{-n}\xi^u)_{\widehat{x}_{-n}}])$.
    \item 
    Then $f^{-n}_{\widehat{x}}\circ f^n|_{\pi_0[({\widehat{f}}^{-n}\xi^u)_{\widehat{x}_{-n}}]} = \mathrm{Id}|_{\pi_0[({\widehat{f}}^{-n}\xi^u)_{\widehat{x}_{-n}}]}$.
\end{enumerate}
\end{lemme}
\noindent\textbf{\underline{Proof :}} According to Lemma \ref{lemma:6411} we have $\xi_{\widehat{x}}^u\subset W^u(\widehat{z},\mathcal{R})$, thus we deduce that the following equalities hold $f^{-n}_{\widehat{z}}[\pi_0(\xi_{\widehat{x}}^u)]=\{y_{-n},\widehat{y}\in\xi_{\widehat{x}}^u\}=\pi_0[\widehat{f}^{-n}(\xi_{\widehat{x}}^u)]=\pi_0[(\widehat{f}^{-n}\xi^u)_{\widehat{x}_{-n}}]$.

For the second item, recall that in Lemma \ref{lemma:f-ncoincidesontheunstablemanifold} (see also Proposition \ref{prop:extensionofnormalforms}) we have proved that $f^{-n}_{\widehat{z}}$ coincide with $f^{-n}_{\widehat{x}}$ on $B(x_0,2\eta(\widehat{x}))\cap B(z_0,2\eta(\widehat{z}))$. In particular, $f^{-n}_{\widehat{z}}$ extends holomorphically $f^{-n}_{\widehat{x}}$ on $B(z_0,2\eta(\widehat{z}))$. We should write $f^{-n}_{\widehat{x}}$ instead of $f_{\widehat{z}}^{-n}$ in what follows. So we have 
$$\pi_0(\xi_{\widehat{x}}^u)\subset \pi_0\left(W^u(\widehat{z},\mathcal{R})\right) =B(z_0,\mathcal{R})\subset B(z_0,2\eta(\widehat{z}))\subset\mathrm{Dom}(f^{-n}_{\widehat{x}}),$$
and thus $f^{n}(\pi_0[(f^{-n}\xi)_{\widehat{x}_{-n}}])=\pi_0[\xi_{\widehat{x}}^u]$ belongs to the domain of definition of the holomorphic map $f^{-n}_{\widehat{x}}$. 

Let us prove now the third item. At this stage the map $h:=f^{-n}_{\widehat{x}}\circ f^n|_{\pi_0[({\widehat{f}}^{-n}\xi^u)_{\widehat{x}_{-n}}]}$ is well defined. Let us verify that $h=\mathrm{Id}$. By definition one has :
\begin{equation}\label{eq:f^n(h(p))=f^n(p)}
    \forall p\in \pi_0[({\widehat{f}}^{-n}\xi^u)_{\widehat{x}_{-n}}],\ f^n(h(p)) = f^n(p).
\end{equation}
Since $f^{-n}_{\widehat{x}}$ and $f^{-n}_{\widehat{z}}$ coincide on $B(z_0,2\eta(\widehat{z}))\supset\pi_0(\xi^u_{\widehat{x}})$, we have by using the first item of this lemma :
$$h(\pi_0[({\widehat{f}}^{-n}\xi^u)_{\widehat{x}_{-n}}]) = f^{-n}_{\widehat{x}}(\pi_0[\xi^u_{\widehat{x}}]) = f^{-n}_{\widehat{z}}(\pi_0[\xi^u_{\widehat{x}}]) = \pi_0[({\widehat{f}}^{-n}\xi^u)_{\widehat{x}_{-n}}].$$
Then (\ref{eq:f^n(h(p))=f^n(p)}) and injectivity of $f^n$ on $\pi_0[({\widehat{f}}^{-n}\xi^u)_{\widehat{x}_{-n}}]$ (given by Theorem \ref{thm:propertiesofeta}) imply $h(p)=p$ for all $p\in\pi_0[({\widehat{f}}^{-n}\xi^u)_{\widehat{x}_{-n}}]$.\qed

\section{Building \texorpdfstring{$q_{\widehat{x}}$}{TEXT} using products of \texorpdfstring{$T^r$}{TEXT} with normal forms}\label{sec:buildingq_x}

We now assume that there exists $r\in\{1,\dots,k-1\}$ such that $\lambda_r>\lambda_{r+1}$, and that the Lyapunov exponents $(\lambda_{r+1},\dots,\lambda_k)$ are not resonant. {These two assumptions are important for constructing the measures $q_{\widehat{x}}$. To prove that these measures are the conditional measures on the atoms of $\xi^u$ of a probability $\widehat{q}$, we will need furthermore that the exponents $\lambda_{r+1} = \cdots = \lambda_k$ are equal.}

\subsection{Slices of \texorpdfstring{$T^r$}{TEXT} with respect to the submersions \texorpdfstring{$W_{\widehat{x}}^{r+1},\cdots,W^{k}_{\widehat{x}}$}{TEXT}}\label{definitionofmesq_x}

For $\widehat{x}\in\mathcal{A}\cap\Lambda$, Theorem \ref{thm:normalforms} and Proposition \ref{prop:extensionofnormalforms} provides holomorphic coordinates 
\begin{equation*}
\sigma_{\widehat{x}}:B(x_0,2\mathcal{R})\longrightarrow\cmplex^k, \ \ \sigma_{\widehat{x}}(x_0)=0 , \ \ \sigma_{\widehat{x}} = (Z_{\widehat{x}},W_{\widehat{x}})=(Z_{\widehat{x}}^1,\cdots,Z_{\widehat{x}}^{r},W_{\widehat{x}}^{r+1},\cdots,W_{\widehat{x}}^k). 
\end{equation*}

\begin{prop}\label{eq:pi(eta)includedinDom(f-n)} 
Let $\widehat{x}\in\mathcal{A}\cap\Lambda$. Let $\widehat{z}$ be the unique element of $\mathcal{T}$ such that $\widehat{x}\in W^u(\widehat{z},\mathcal{R})$. Recall that $\xi_{\widehat{x}}^u\subset W^u(\widehat{z},\mathcal{R})$, it yields by using Proposition \ref{prop:extensionofnormalforms}:
\begin{equation*}
    \pi_0\left(\xi^u_{\widehat{x}}\right)\subset B(z_0,\mathcal R) \subset B(x_0,2\mathcal{R})\subset\mathrm{Dom}\left(W_{\widehat{x}}\right).
\end{equation*}
\textcolor{black}{In particular, $\pi_0\left(\xi^u_{\widehat{x}}\right)\subset B(z_0,\mathcal R) \subset B(x_0,2\mathcal{R}) \cap B(y_0,2\mathcal{R})$ for every $\widehat y \in \xi^u_{\widehat{x}}\cap\Lambda$.}
\end{prop}

In consequence, the measures $T^r\wedge dd^c|W_{\widehat{x}}|^2$ (with ${\widehat{x}\in\mathcal{A}}\cap\Lambda$) are defined on balls $B(x_0,2\mathcal{R})$ containing $\pi_0(\xi_{\widehat{x}}^u)$, where we have denoted
$$dd^c|W_{\widehat{x}}|^2:=dd^c|W_{\widehat{x}}^{r+1}|^2\wedge\cdots\wedge dd^c|W_{\widehat{x}}^k|^{2}.$$
The following proposition ensures in particular that the measure $(T^r\wedge dd^c|W_{\widehat{x}}|^2)\circ \pi_0\lfloor_{\xi_{\widehat{x}}^u}$ does not vanish for each $\widehat{x}\in\mathcal{A}\cap\Lambda$. 

\begin{prop}\label{prop:62}
Let $\widehat{x}\in\Lambda$.
\begin{enumerate} 
\item If $\widehat{x}\in\mathcal{A}$, then $T^r\wedge dd^c|W_{\widehat{x}}|^2$ does not vanish on every ball $B(x_0,\eta)$ with ${0<\eta\leq \eta(\widehat{x})}$. 
\item For $n\geq0$ such that $\widehat{x}_{-n}\in\mathcal{A}$, we have $\left(T^r\wedge dd^c|W_{\widehat{x}_{-n}}|^2\right)\circ\pi_0\circ {\widehat{f}}^{-n}(\xi^u_{\widehat{x}})>0$.
\end{enumerate}
\end{prop}
\noindent\textbf{\underline{Proof :}} The first item is a consequence of Theorem \ref{thm:SupportofT^rwedgedd^c|w|^2}, since $T$ has locally Hölder continuous potentials, as stated in Proposition \ref{prop:HolderpotentialsfortheGreenCurrent}. So, let us prove the second item. We use the inclusion 
$$\xi^u_{\widehat{x}}\supset W^u(\widehat{x},\eta'(\widehat{x}))$$
given by the point \textit{3.(c)} of \Dcitelem{\ref{lemma:6411}}. We deduce by definition of $W^u(\widehat{x},\eta'(\widehat{x}))$ that :
$$\pi_0[{\widehat{f}}^{-n}(\xi^u_{\widehat{x}})] \supset \pi_0\left(\{\widehat{y}_{-n},\ \widehat{y}\in W^u(\widehat{x},\eta'(\widehat{x}))\}\right) = f^{-n}_{\widehat{x}}(B(x_0,\eta'(\widehat{x}))) \supset B(x_{-n},\eta''),$$
where $0<\eta''\leq\eta(\widehat{x}_{-n})$ is small enough to get the last inclusion. 
Since $\pi_0(\widehat{f}^{-n}(\xi_{\widehat{x}}^u))\subset\pi_0(\xi_{\widehat{x}_{-n}}^u)\subset\mathrm{Dom}(W_{\widehat{x}_{-n}})$
by Proposition \ref{eq:pi(eta)includedinDom(f-n)}, the quantity $\left(T^r\wedge dd^c|W_{\widehat{x}_{-n}}|^2\right)[\pi_0({\widehat{f}}^{-n}(\xi^u_{\widehat{x}}))]$ is well defined and bounded from below by $\left(T^r\wedge dd^c|W_{\widehat{x}_{-n}}|^2\right)(B(x_{-n},\eta''))$. The first item applied to $\widehat{x}_{-n}$ and $\eta''$ then allows to conclude.\qed

{\subsection{Recurrent properties}\label{sec:recproperties}}

Recall that $V_s\subset T\mathbb{P}^k$ is the bundle given by Oseledec Theorem \ref{thm:oseledec}.\\

\begin{defn}\label{eq:M(x)andDelta(x,p)} For each $\widehat{x}\in\Lambda$ we define:
    \begin{enumerate}
        \item $M(\widehat{x}):=|\mathrm{det}_{\cmplex}(d_0\sigma_{\widehat{x}}^{-1}:\{0\}^{r}\times\cmplex^{k-r}\to V_s(x_0))|$.
        \item $\Delta(x_0,p):=|\mathrm{det}_{\cmplex}(d_{x_0}f^p:V_s(x_0)\to V_s(f^px_0))|$ for all $p\in\mathbb{N}$. 
        \item $\beta_{n,\widehat{x}}(r):=\prod_{l=r+1}^{k}\beta_{n,\widehat{x}}^l$, where the $\beta_{n,\widehat{x}}^l$ are given by Theorem \ref{thm:normalforms}.
    \end{enumerate}
\end{defn}

The following proposition ensures that the measures $q_{\widehat{x}}$ will be well defined, see Remark \ref{rmq:remarkonthedefinitionofthemeasuresq_x}. Up to a modification of $\Lambda$ (defined in Definition \ref{defn:invariantsetLambda_epsilon}) using Poincaré recurrence theorem, the set of integers  $\mathrm{Rec}_{\widehat{x}}(\mathcal{A}):=\{ p \geq 0, \ \widehat{x}_{-p}\in\mathcal{A}\}$ is infinite for every $\widehat{x}\in\Lambda$.

\begin{prop}\label{lemma:66} Let $\widehat{x}\in \Lambda$, $(p,q)\in\mathrm{Rec}_{\widehat{x}}(\mathcal{A})^2$ and $A\subset\xi^u_{\widehat{x}}$ be a Borel set. According to Proposition \ref{eq:pi(eta)includedinDom(f-n)} we have $\pi_0(\widehat{f}^{-p}(A))\subset \mathrm{Dom}(W_{\widehat{x}_{-p}})$, so let us define:
$$E^p:=d^{rp}\Delta(x_{-p},p)^2M(\widehat{x}_{-p})^2(T^r\wedge dd^c|W_{\widehat{x}_{-p}}|^2)\circ\pi_0\circ {\widehat{f}}^{-p}\left(A\right),$$ and let $E^q$ be defined similarly replacing $p$ by $q$. Then $E^p=E^q$.
\end{prop}

Before going into the proof of Proposition \ref{lemma:66}, we need two lemmas.

\begin{lemme}\label{lemma:63} For each $\textcolor{black}{\widehat{x}\in\mathcal{A}\cap}\Lambda$ and for each $n\geq N(\widehat{x})$ we have:
$$(f^{-n}_{\widehat{x}})^*dd^c|W_{\widehat{x}_{-n}}|^2=|\beta_{n,\widehat{x}}(r)|^2dd^c|W_{\widehat{x}}|^2\ \mathrm{on}\ B(x_0,2\mathcal{R}).$$
\end{lemme}
\noindent\textbf{\underline{Proof :}} The commutative diagram \eqref{eq:extendeddiagramnormalforms} of Proposition \ref{prop:extensionofnormalforms} provides for $n\geq N(\widehat{x})$ and for $l\in\{r+1,\cdots,k\}$:
\begin{align*}
    W_{\widehat{x}_{-n}}^l\circ f^{-n}_{\widehat{x}} & = W^l\circ\left(\sigma_{\widehat{x}_{-n}}\circ f_{\widehat{x}}^{-n}\right) = \left(W^l\circ R_{n,\widehat{x}}\right)\circ\sigma_{\widehat{x}} = \beta_{n,\widehat{x}}^l\times W_{\widehat{x}}^l\ \mathrm{on}\ B(x_0,2\mathcal{R}),
\end{align*}
where $W^l : \mathbb{C}^k \to \mathbb{C}$ denotes the projection onto the $l$-th coordinate. The conclusion follows by taking wedge products.\qed

\begin{lemme}\label{lemma:cocyclerelationDelta} For every $\widehat{x}\in\Lambda$ we have:
\begin{enumerate}
    \item For all $n,k \geq 0$, $\Delta(x,n)\Delta(x_{n},k)=\Delta(x,n+k)$.
    \item For all $n\geq N(\widehat{x})$, $\Delta(x_{-n},n)|\beta_{n,\widehat{x}}(r)| = \frac{M(\widehat{x})}{M(\widehat{x}_{-n})}$.
\end{enumerate}
\end{lemme}
\noindent\textbf{\underline{Proof :}}
The first item comes from the chain rule formula. For the second item, we proceed as follows. First, since $R_{n,\widehat{x}}$ has a triangular form \eqref{eq:Rnxwithnoresonancesforlowerlyapunovexponents}, we have for $(z,w)\in\cmplex^k$ that $d_0 R_{n,\widehat{x}}(z,w)=(\alpha^1_{n,\widehat{x}}z_1, \dots, \beta^k_{n,\widehat{x}}w_k)$. So, denoting $L_r := \{0\}^r \times \mathbb{C}^{k-r}$, we have
$$|\beta_{n,\widehat{x}}(r)|=\left|\prod_{l=r+1}^{k}\beta_{n,\widehat{x}}^l\right|=|\mathrm{det_{\cmplex}}(d_0R_{n,\widehat{x}}:L_r\to L_r)|.$$
Second, the commutative diagram \eqref{eq:diagramnormalforms} of Theorem \ref{thm:normalforms} gives for $n\geq N(\widehat{x})$ :
$$|\mathrm{det_{\cmplex}}(d_0R_{n,\widehat{x}}:L_r\to L_r)|= \left|\mathrm{det}_{\cmplex}(d_{x_{-n}}\sigma_{\widehat{x}_{-n}}\circ d_{x_{0}}f^{-n}_{\widehat{x}}\circ d_{0}\sigma_{\widehat{x}}^{-1}:L_r\to L_r)\right|.$$
One concludes by Theorems \ref{thm:oseledec} and \ref{thm:normalforms} which yield that $|\beta_{n,\widehat{x}}(r)|$ is equal to the product 
{\footnotesize
\begin{align*}
    |\mathrm{det}_{\cmplex}(d_{x_{-n}}\sigma_{\widehat{x}_{-n}}:V_s(x_{-n})\to L_r)|
    \times |\mathrm{det}_{\cmplex}(d_{x_0}f^{-n}_{\widehat{x}}:V_s(x_{0})\to V_s(x_{-n}))|
    \times |\mathrm{det}_{\cmplex}(d_{0}\sigma_{\widehat{x}}^{-1}:L_r\to V_s(x_0))|,
\end{align*}
}%
and thus $|\beta_{n,\widehat{x}}(r)|=M(\widehat{x}_{-n})^{-1}\Delta(x_{-n},n)^{-1}M(\widehat{x})$.\qed\\
      
\noindent\textbf{\underline{Proof of Proposition \ref{lemma:66}:}} For proving $E^p=E^q$ for any $p,q\in\mathrm{Rec}_{\widehat{x}}(\mathcal{A})$, it suffices to show the following property $(P)$:
$$\left[\displaystyle{\phantom{0^{0^{0}}}\!\!\!\!\!\!\!\forall} p,q\in\mathrm{Rec}_{\widehat{x}}(\mathcal{A}),\ p-q\geq N(\widehat{x}_{-q}) \Longrightarrow E^p=E^q\right] \ :\ (P)$$
Indeed, since $\mathrm{Rec}_{\widehat{x}}(\mathcal{A})$ is infinite, for any $p,q\in\mathrm{Rec}_{\widehat{x}}(\mathcal{A})$ we can find $s\in\mathrm{Rec}_{\widehat{x}}(\mathcal{A})$ such that $s-p\geq N(\widehat{x}_{-p})$ and $s-q\geq N(\widehat{x}_{-q})$. Thus if $(P)$ is true, $E^s=E^p$ and $E^s=E^q$, which implies $E^p=E^q$. 
Let us prove now the property $(P)$.
Let us fix $p,q\in\mathrm{Rec}_{\widehat{x}}(\mathcal{A})$ such that $n:=p-q\geq N(\widehat{x}_{-q})$. Let us denote $F^p := (T^r \wedge dd^c |W_{\widehat{x}_{-p}}|^2) \circ \pi_0 \circ {\widehat{f}}^{-p}(A)$, and similarly, let us denote $F^q$ by replacing $p$ with $q$. Proving $E^p=E^q$ amounts to show the equality :
\begin{equation}\label{eq:d^pM^2Delta^2E_1=d^qM^2Delta^2E_2}
    d^{rp}\Delta(x_{-p},p)^2M(\widehat{x}_{-p})^2F^p=d^{rq}\Delta(x_{-q},q)^2M(\widehat{x}_{-q})^2F^q.
\end{equation}
Since $\widehat{x}_{-q}\in\mathcal{A}$ and ${\widehat{f}}^{-p}(\xi^u_{\widehat{x}})\subset \left({\widehat{f}}^{-n}\xi^u\right)_{\widehat{x}_{-q-n}}$,
one gets by \Dcitelem{\ref{lemma:f-ncircfn=id}} (applied to $\widehat{x}_{-q}$ with $n=p-q$) :
\begin{equation}\label{eq:f^-(p-q)circf^(p-q)=Id}
    f^{-(p-q)}_{\widehat{x}_{-q}}\circ f^{p-q}|_{\pi_0({\widehat{f}}^{-p}(\xi^u_{\widehat{x}}))} = \mathrm{Id}|_{\pi_0({\widehat{f}}^{-p}(\xi^u_{\widehat{x}}))}.
\end{equation}
Let us denote $U:=\pi_0({\widehat{f}}^{-p}(\xi^u_{\widehat{x}}))$ and $f^{p-q}_U:=f^{p-q}|_U$. By \eqref{eq:f^-(p-q)circf^(p-q)=Id} and by $(f^{p-q}_U)^*T^r = d^{r(p-q)} T^r$, we obtain 
\begin{align*}
                                   F^p & = \left[T^r\wedge\ (f^{-(p-q)}_{\widehat{x}_{-q}}\circ f^{p-q}_U)^* dd^c|W_{\widehat{x}_{-p}}|^2\right]\circ\pi_0\circ {\widehat{f}}^{-p}(A)  \\
                                   & = d^{r(q-p)}\left[(f^{p-q}_U)^*T^r\wedge\ (f^{p-q}_U)^* dd^c\left|W_{\widehat{x}_{-p}}\circ f^{-(p-q)}_{\widehat{x}_{-q}}\right|^2\right]\circ\pi_0\circ {\widehat{f}}^{-p}(A) .
\end{align*}
Then by injectivity of $f^{p-q}_U$ we have :
$$  F^p = d^{r(q-p)}\left[T^r\wedge\ dd^c\left|W_{\widehat{x}_{-p}}\circ f^{-(p-q)}_{\widehat{x}_{-q}}\right|^2\right]\circ(f^{p-q}_U)\circ\pi_0\circ {\widehat{f}}^{-p}(A). $$
Using the relation $f\circ\pi_0=\pi_0\circ {\widehat{f}}$, we have $(f_U^{p-q})\circ\pi_0\circ {\widehat{f}}^{-p}(A)=\pi_0\circ {\widehat{f}}^{-q}(A)$. Moreover, by Lemma \ref{lemma:63} with $\widehat{x}$ replaced by $\widehat{x}_{-q}\in\mathcal{A}$ and $n=p-q\geq N(\widehat{x}_{-q})$, we get : 
$$dd^c|W_{\widehat{x}_{-p}}\circ f^{-(p-q)}_{\widehat{x}_{-q}}|^2= |\beta_{p-q,\widehat{x}_{-q}}(r)|^2\times dd^c|W_{\widehat{x}_{-q}}|^2\ \mathrm{on}\ B(x_{-q},2\mathcal{R}).$$
\textcolor{black}{Using again Proposition \ref{eq:pi(eta)includedinDom(f-n)} we have $\pi_0(\widehat{f}^{-q}(A))\subset\pi_0(\xi_{\widehat{x}_{-q}}^u)\subset B(x_{-q},2\mathcal{R})$}, and thus we have $dd^c|W_{\widehat{x}_{-p}}\circ f^{-(p-q)}_{\widehat{x}_{-q}}|^2=|\beta_{p-q,\widehat{x}_{-q}}(r)|^2\times dd^c|W_{\widehat{x}_{-q}}|^2$ on $\pi_0(\widehat{f}^{-q}(A))$. We then deduce: 
$$F^p=d^{r(q-p)}|\beta_{p-q,\widehat{x}_{-q}}(r)|^2\times F^q.$$
We replace now $|\beta_{p-q,\widehat{x}_{-q}}(r)|^2$ by applying the second item of Lemma \ref{lemma:cocyclerelationDelta} to $\widehat{x}_{-q}$ with $n=p-q\geq N(\widehat{x}_{-q})$. We get 
$$d^{rp}M(\widehat{x}_{-p})^2F^p=d^{rq}M(\widehat{x}_{-q})^2F^q\times \frac{1}{\Delta(x_{-p},p-q)^2}.$$
To conclude, we use $\Delta(x_{-p},p-q)\Delta(x_{-q},q)=\Delta(x_{-p},p)$ which comes from Lemma \ref{lemma:cocyclerelationDelta} applied to $\widehat {x}_{-p}$, $n=p-q$ and $k=q$.\qed 

\subsection{The measures \texorpdfstring{$q_{\widehat{x}}$}{TEXT} and \texorpdfstring{$\widehat{q}$}{TEXT}}\label{definitionofthemeasuresq_xunderH_1}

Recall that we have defined $dd^c|W_{\widehat{x}}|^2:=dd^c|W_{\widehat{x}}^{r+1}|^2\wedge\cdots\wedge dd^c|W_{\widehat{x}}^k|^2$.

\begin{defn}\label{defn:64} For every $\widehat{x}\in\Lambda$ we define 
\begin{equation*}
        \widetilde{q}_{\widehat{x}} := d^{rp}\Delta(x_{-p},p)^2M(\widehat{x}_{-p})^2\times \left(T^r\wedge dd^c|W_{\widehat{x}_{-p}}|^2\right)\circ\pi_0\circ {\widehat{f}}^{-p}\ \mathrm{on}\ \xi^u_{\widehat{x}},
\end{equation*}
where $p\geq0$ is any integer such that $\widehat{x}_{-p}$ belongs to $\mathcal{A}$, and where $M(\cdot)$ and $\Delta(\cdot,\cdot)$ are defined in Definition \ref{eq:M(x)andDelta(x,p)}. The probability measure $q_{\widehat{x}}$ is then defined by (where $\mathcal{B}$ is the $\sigma$-algebra of Borel sets)
$$\forall A\in\mathcal{B},\ q_{\widehat{x}}(A):=\frac{\widetilde{q}_{\widehat{x}}(A\cap \xi^u_{\widehat{x}})}{L(\widehat{x})},$$
the normalization $L(\widehat{x}) := \widetilde{q}_{\widehat{x}}(\xi^u_{\widehat{x}})$ being $>0$  by \Dciteprop{\ref{prop:62}}.
\end{defn}
\begin{rmq}\label{rmq:remarkonthedefinitionofthemeasuresq_x} \ 
\begin{enumerate} 
    \item In the definition of $\widetilde{q}_{\widehat{x}}$ we need to introduce $p\in\mathrm{Rec}_{\widehat{x}}(\mathcal{A})$ to guarantee that $\widetilde{q}_{\widehat{x}}$ is well defined. Indeed, according to Proposition \ref{eq:pi(eta)includedinDom(f-n)}, $T^r\wedge dd^c|W_{\widehat{x}_{-p}}|^2$ is well defined on $B(x_{-p},2\mathcal{R})$ and this ball contains $\pi_0(\xi_{\widehat{x}_{-p}}^u)$ which itself contains $\pi_0(\widehat{f}^{-p}(\xi_{\widehat{x}}^u))$ by the decreasing property of $\xi^u$. Moreover, Proposition \ref{lemma:66} ensures that the definition of $\widetilde{q}_{\widehat{x}}$ provided here does not depend on the choice of $p\in\mathrm{Rec}_{\widehat{x}}(\mathcal{A})$. 
    \item By \Dcitethm{\ref{thm:propertiesofeta}}, $\pi_0$ is injective on the atoms of ${\widehat{f}}^{-p}\xi^u$, thus $\widetilde{q}_{\widehat{x}}$ is a measure. 
    \item The definition of $q_{\widehat{x}}$ given here only requires that $\lambda_r > \lambda_{r+1}$ and that the exponents $\lambda_{r+1} \geq \cdots \geq \lambda_k$ are not resonant. However, to justify that the $q_{\widehat{x}}$ are conditional measures of a probability $\widehat{q}$, our arguments require that $\lambda_{r+1} = \cdots = \lambda_k$, cf. Proposition \ref{prop:q_x=q_ysieta_x=eta_y} just below.
    \item If $\widehat{x}\in \mathcal{A}$, then  $\widetilde{q}_{\widehat{x}} = M(\widehat{x})^2 
    \left( T^r\wedge dd^c |W_{\widehat{x}}|^2 \right) \circ \pi_0$ on $\xi^u_{\widehat{x}}$.
    \item By denoting $\mathcal{X}:=\bigcup_{\widehat{x}\in\mathcal{FN}}\{\widehat{x}\}\times B(x_0,2\eta(\widehat{x}))\subset\widehat{\mathbb{P}^k}\times\mathbb{P}^k$, and by looking at the proof of Theorem \ref{thm:normalforms} in \cite{bdm07}, one observes that the map $\mathcal{X}\to\cmplex^k$ defined by $(\widehat{x},z)\mapsto\sigma_{\widehat{x}}(z)$ is measurable. Moreover, the $dd^c$ does not affect this measurability, and $(\widehat{x},z)\in\mathcal{X}\mapsto dd^c_z|W_{\widehat{x}}|^2$ is also measurable from $\mathcal{X}$ to $\bigwedge^{k-r,k-r}T^*\mathbb{P}^k$.
    \item In particular, one can check that $\widehat{x}\to L(\widehat{x})$ is measurable, and that for a fixed Borel set $A$, the map $\widehat{x}\mapsto q_{\widehat{x}}(A\cap\xi^u_{\widehat{x}})$ is also measurable.
\end{enumerate} 
\end{rmq}

\newpage

\begin{prop}\label{prop:q_x=q_ysieta_x=eta_y}
    Assume that $\lambda_r>\lambda_{r+1}=\cdots=\lambda_k$. 
    \begin{enumerate}
        \item Let $\widehat{x}\in\mathcal{A}\cap\Lambda$. Let $\widehat z$ be the unique element of $\mathcal T$ such that $\widehat{x}\in W^u(\widehat{z},\mathcal{R})$. Let $\widehat{y}\in\Lambda$ and assume $\widehat{y}\in\xi_{\widehat{x}}^u$. Then the following points hold:
        \begin{enumerate}
            \item There exists $C\in\mathrm{GL}_{k-r}(\cmplex)$ such that $dW_{\widehat{y}}=CdW_{\widehat{x}}\ \mathrm{on}\ B(z_0,\mathcal{R})$.
            \item $q_{\widehat{x}}=q_{\widehat{y}}$ on $\xi^u_{\widehat{x}}$.
        \end{enumerate}
        \item More generally, for all $\widehat{x}\in\Lambda$, we have
            \begin{equation}\label{eq:q_x=q_yoneta_x2}
                \textcolor{black}{\forall\widehat{y}\in\Lambda,\ \widehat{y}\in\xi_{\widehat{x}}^u\Longrightarrow q_{\widehat{x}}=q_{\widehat{y}}.}
            \end{equation}
    \end{enumerate}
    By Remark \ref{rmq:remarkonthedefinitionofthemeasuresq_x} the function $\widehat{x}\mapsto q_{\widehat{x}}(A\cap\xi^u_{\widehat{x}})$ is measurable for all Borel set $A$, thus the following probability measure is well defined
    $$\forall A\in\mathcal{B},\ \widehat{q}(A) := \int_{\widehat{\mathbb{P}^k}}q_{\widehat{x}}(A\cap \xi^u_{\widehat{x}})\ \mathrm{d}\widehat{\mu}(\widehat{x}).$$
    According to \eqref{eq:q_x=q_yoneta_x2}, the measures $q_{\widehat{x}}$ are the conditional measures of $\widehat{q}$ on the atoms $\xi_{\widehat{x}}^u$.
\end{prop}
\noindent\textbf{\underline{Proof:}}\ \\
\noindent\textit{1.} Let $\widehat{x}\in\mathcal{A}\cap\Lambda$, and let $\widehat{z}$ be the unique element of $\mathcal{T}$ such that $\widehat{x}\in W^u(\widehat{x},\mathcal{R})$. Let $\widehat{y}\in\Lambda$ such that $\widehat{y}\in\xi_{\widehat{x}}^u$, we have to show that $q_{\widehat{x}}=q_{\widehat{y}}$. Lemma \ref{lemma:6411} yields $\xi_{\widehat{x}}^u\subset W^{u}(\widehat{z},\mathcal{R})$ which implies $\widehat{y}\in W^{u}(\widehat{z},\mathcal{R})\subset\mathcal{A}$, thus by the fourth point of Remark \ref{rmq:remarkonthedefinitionofthemeasuresq_x} we have:
$$q_{\widehat{y}}=\frac{M(\widehat{y})^2}{L(\widehat{y})}(T^r\wedge dd^c|W_{\widehat{y}}|^2)\circ\pi_0\lfloor_{\xi^u_{\widehat{y}}}\ \mathrm{and}\ q_{\widehat{x}}=\frac{M(\widehat{x})^2}{L(\widehat{x})}(T^r\wedge dd^c|W_{\widehat{x}}|^2)\circ\pi_0\lfloor_{\xi^u_{\widehat{x}}}.$$
Since $(\mathcal{A},\mathcal{R},\rho,\mathcal{T},c)$ is a Pesin box which satisfies $\mathcal{T}\subset\left\{\widehat{t}\in\Lambda:\frac{1}{\rho}\leq\eta\left(\hat{t}\right)\right\}$ and $0<\mathcal{R}\leq\frac{1}{2\rho}$ by construction in Section \ref{sec:constructionofA}, and since {$\lambda_r>\lambda_{r+1}$}, we can use Theorem \ref{thm:localfoliations} which gives the existence of $C\in\mathrm{GL}_{k-r}({\mathcal{O}(B(z_0,\mathcal{R}))})$ such that $dW_{\widehat{y}}=CdW_{\widehat{x}}$ on $B(z_0,\mathcal{R})$. Then, an explicit computation gives
$$dd^c|W_{\widehat{y}}|^2=\left|\mathrm{det}(C)\right|^2dd^c|W_{\widehat{x}}|^2\ \mathrm{on}\ B(z_0,\mathcal{R}).$$
Thus, by using that $\xi^u_{\widehat{x}}=\xi^u_{\widehat{y}}$ and by using that $B(z_0,\mathcal{R})\supset\pi_0(\xi_{\widehat{x}}^u)$ by Proposition \ref{eq:pi(eta)includedinDom(f-n)}, we get ({where $\mathrm{det}(C)\in\mathcal{O}^*(B(z_0,\mathcal{R}))$})
$$q_{\widehat{y}}=\frac{M(\widehat{y})^2L(\widehat{x})}{M(\widehat{x})^2L(\widehat{y})}\left|\mathrm{det}(C\circ\pi_0)\right|^2\times q_{\widehat{x}}\ \mathrm{on}\ \xi_{\widehat{x}}^u.$$
At this stage, we cannot conclude that $q_{\widehat{y}} = q_{\widehat{x}}$ because $|\mathrm{det}(C \circ \pi_0)|^2$ is not necessarily constant on $\xi^u_{\widehat{x}}$. This is where the assumption $\lambda_{r+1} = \cdots = \lambda_k$ comes into play: according to Theorem \ref{thm:localfoliations}, $\lambda_{r+1} = \cdots = \lambda_k$ implies that $C \in \mathrm{GL}_{k-r}(\mathbb{C})$ is constant, and thus $|\mathrm{det}(C)|^2 \in \mathbb{C}^*$ is also constant. We then obtain $q_{\widehat{y}} = \frac{M(\widehat{y})^2 L(\widehat{x})}{M(\widehat{x})^2 L(\widehat{y})} \left| \mathrm{det}(C) \right|^2 \times q_{\widehat{x}}$ on $\xi^u_{\widehat{x}}$, with $\frac{M(\widehat{y})^2 L(\widehat{x})}{M(\widehat{x})^2 L(\widehat{y})} \left| \mathrm{det}(C) \right|^2 \in \mathbb{C}^*$ being constant. Since $q_{\widehat{x}}$ and $q_{\widehat{y}}$ are probability measures on $\xi_{\widehat{x}}^u$, we finally deduce that this constant is equal to $1$ and that $q_{\widehat{y}} = q_{\widehat{x}}$.\\

\noindent\textit{2.} We have to prove that $\mathrm{\eqref{eq:q_x=q_yoneta_x2}}$ holds for each $\widehat{x}\in\Lambda$. By the preceding item of the proposition, we have proved that \eqref{eq:q_x=q_yoneta_x2} is true for all $\widehat{x}\in\mathcal{A}\cap\Lambda$, so let us fix now $\widehat{x}\in \Lambda$. We recall that $\mathrm{Rec}_{\widehat{x}}(\mathcal{A})=\{n\geq0:\widehat{x}_{-n}\in\mathcal{A}\}$ is not empty, so let us fix $n\in\mathrm{Rec}_{\widehat{x}}(\mathcal{A})$. Let us fix $\widehat{y}\in\Lambda$ and let us assume that $\widehat{y}\in \xi_{\widehat{x}}^u$.
By the decreasing property of $\xi^u$, $\widehat{y}_{-n}\in\xi_{\widehat{x}_{-n}}^u\cap\Lambda$, so by the preceding item of the proposition we have $q_{\widehat{y}_{-n}}=q_{\widehat{x}_{-n}}$, which can be written as
\begin{equation}\label{eq:qql}
\frac{1}{L(\widehat{y}_{-n})} \widetilde{q}_{\widehat{y}_{-n}}=\frac{1}{L(\widehat{x}_{-n})} \widetilde{q}_{\widehat{x}_{-n}} .
\end{equation}
Recall that $\xi_{\widehat{x}_{-n}}^u\subset\mathcal{A}$ (by the second item of Lemma \ref{lemma:6411} since $\widehat{x}_{-n}\in\mathcal{A}$), thus $\widehat{y}_{-n}\in\mathcal{A}$. Now using the definition of  $q_{\widehat{y}}$ we have :
$$L(\widehat{y})q_{\widehat{y}} = d^{rn}\Delta(y_{-n},n)^2 \times  M(\widehat{y}_{-n})^2\left(T^r\wedge dd^c|W_{\widehat{y}_{-n}}|^2\right)\circ\pi_0\circ {\widehat{f}}^{-n}\lfloor_{\xi^u_{\widehat{y}}} .$$
Since ${\widehat{f}}^{-n}(\xi^u_{\widehat{y}}) \subset \xi^u_{\widehat{y}_{-n}}$, one can replace $\pi_0\circ {\widehat{f}}^{-n}\lfloor_{\xi^u_{\widehat{y}}}$ by $\pi_0 \lfloor_{\xi^u_{\widehat{y}_{-n}}} \circ {\widehat{f}}^{-n}\lfloor_{\xi^u_{\widehat{y}}}$. Moreover, since $\widehat{y}_{-n}\in\mathcal{A}$ we get by the fourth item of Remark \ref{rmq:remarkonthedefinitionofthemeasuresq_x} that $\widetilde{q}_{\widehat{y}_{-n}} =  M(\widehat{y}_{-n})^2\left(T^r\wedge dd^c|W_{\widehat{y}_{-n}}|^2\right)\circ \pi_0 \lfloor_{\xi^u_{\widehat{y}_{-n}}}.$
We finally deduce
$$L(\widehat{y})q_{\widehat{y}} = d^{rn}\Delta(y_{-n},n)^2\times \widetilde{q}_{\widehat{y}_{-n}}\circ {\widehat{f}}^{-n}\lfloor_{\xi^u_{\widehat{y}}},$$
and  similarly $L(\widehat{x})q_{\widehat{x}} = d^{rn}\Delta(x_{-n},n)^2\times \widetilde{q}_{\widehat{x}_{-n}}\circ {\widehat{f}}^{-n}\lfloor_{\xi^u_{\widehat{x}}}$. \Dciteequa{\ref{eq:qql}} yields using $\xi_{\widehat{x}}^u=\xi_{\widehat{y}}^u$:
$$q_{\widehat{y}} = \frac{\Delta(y_{-n},n)^2L(\widehat{y}_{-n})L(\widehat{x})}{\Delta(x_{-n},n)^2L(\widehat{x}_{-n})L(\widehat{y})}\times q_{\widehat{x}}.$$
At last, the quotient is equal to $1$ and $q_{\widehat{y}} = q_{\widehat{x}}$ since these are probability measures $\mathrm{on}\ \xi^u_{\widehat{x}}$.\qed

{\section{Proofs of Theorem {\ref{thm:TheoremB}} and Theorem {\ref{thm:TheoremA}}}\label{sec:constructionofq_x}}

\subsection{Proof of Theorem {\ref{thm:TheoremB}} Items \textit{1.} to \textit{5.}}

We assume $\lambda_r>\lambda_{r+1}$ and that $\lambda_{r+1}\geq\cdots\geq\lambda_k$ are not resonant. The formula \textit{5.} of Theorem \ref{thm:TheoremB} is contained in Theorem \ref{thm:propertiesofeta}. Let us now prove that formula \textit{4.} of Theorem \ref{thm:TheoremB} holds for every $n\geq0$:
\begin{equation}\label{eq:laformuledutheoremB}
    \LLog\ d^{rn} + 2n(\lambda_{r+1}+\cdots+\lambda_k) = \int_{\widehat{\mathbb{P}^k}}-\LLog\ q_{\widehat{x}}\left({\widehat{f}}^{-n}\xi^u\right)_{\widehat{x}}\ \mathrm{d}\widehat{\mu}(\widehat{x}).
\end{equation}

\begin{prop}\label{prop:z20}
For all $\widehat{x}\in\Lambda$ and for every $n \geq 0$:
\begin{equation}\label{eq:L(x)=q_x(g-neta)=L(x_n)/dnDelta(x,n)2}
    L(\widehat{x})q_{\widehat{x}}\left(\left({\widehat{f}}^{-n}\xi^u\right)_{\widehat{x}}\right) = \frac{1}{d^{rn}}\frac{L(\widehat{x}_n)}{\Delta(x_0,n)^2}.
\end{equation}
\end{prop}
\noindent\textbf{\underline{Proof :}} Let $\widehat{x}\in\Lambda$ and $p\in\mathbb{N}$ such that $\widehat{x}_{-p}\in\mathcal{A}$. Let us fix $n\in\mathbb{N}$ and let $q:=n+p$ so that $\widehat{x}_{n-q}=\widehat{x}_{-p}\in\mathcal{A}$. By definition of $q_{\widehat{x}}$ and using $({\widehat{f}}^{-n}\xi^u)_{\widehat{x}}\subset\xi^u_{\widehat{x}}$ we have:
$$L(\widehat{x})q_{\widehat{x}}\left(\left({\widehat{f}}^{-n}\xi^u\right)_{\widehat{x}}\right) = d^{rp}\Delta(x_{-p},p)^2M(\widehat{x}_{-p})^2\times \left[T^r\wedge dd^c|W_{\widehat{x}_{-p}}|^2\right]\circ\pi_0\circ {\widehat{f}}^{-p}\left(\left({\widehat{f}}^{-n}\xi^u\right)_{\widehat{x}}\right).$$
Using $p=q-n$ and ${\widehat{f}}^{-p}(({\widehat{f}}^{-n}\xi^u)_{\widehat{x}})= {\widehat{f}}^{-q}(\xi^u_{\widehat{x}_{n}})$ we have :
$$L(\widehat{x})q_{\widehat{x}}\left(\left({\widehat{f}}^{-n}\xi^u\right)_{\widehat{x}}\right) = \frac{1}{d^{rn}} \Delta(x_{n-q},p)^2 d^{rq} M(\widehat{x}_{n-q})^2\left[T^r\wedge dd^c\left|W_{\widehat{x}_{n-q}}\right|^2\right]\circ\pi_0\circ {\widehat{f}}^{-q}\left(\xi^u_{\widehat{x}_n}\right).$$
By definition of $L(\widehat{x}_n)$ (use Definition \ref{defn:64} replacing in it $\widehat{x}$ by $\widehat{x}_n$ and $p$ by $q$) we have 
$$L(\widehat{x})q_{\widehat{x}}\left(\left({\widehat{f}}^{-n}\xi^u\right)_{\widehat{x}}\right) = \frac{\Delta(x_{n-q},p)^2}{d^{rn}\Delta(x_{n-q},q)^2}\times L(\widehat{x}_n).$$
The first item of \Dcitelem{\ref{lemma:cocyclerelationDelta}} applied to $\widehat{x}_{n-q}$, $p$ and $q-p$ completes the proof.\qed \\

By taking$\LLog$ in Equation (\ref{eq:L(x)=q_x(g-neta)=L(x_n)/dnDelta(x,n)2}), we obtain :
\begin{equation}\label{eq:logdn+2logDelta=-logq_x+logl(x_n)/L(x)}
    \LLog\ d^{rn} + 2\LLog\ \Delta(x_0,n) = -\LLog\ q_{\widehat{x}}\left({\widehat{f}}^{-n}\xi^u\right)_{\widehat{x}} + \LLog\frac{L(\widehat{x}_n)}{L(\widehat{x})}.
\end{equation}
The third item of \Dcitethm{\ref{thm:oseledec}} asserts that $\int_{\widehat{\mathbb{P}^k}} \LLog\ \Delta(x_0,n)\ \mathrm{d}\widehat{\mu}(\widehat{x}) = n(\lambda_{r+1}+\cdots+\lambda_k)$. To finish the proof of Formula (\ref{eq:laformuledutheoremB}), it remains to show that $h_n:=\LLog\left(\frac{L\circ {\widehat{f}}^n}{L}\right) \in L^1(\widehat{\mu})$ and satisfies $\int h_n\ \mathrm{d}\widehat{\mu}=0$. We use the following classical lemma, the original statement  is stated with$\LLog^{-}$ instead of$\LLog^{+}$, but the proof also works with$\LLog^{+}$.

\begin{lemme}[Ledrappier-Strelcyn {\cite[Proposition 2.2]{ledstr82}}]\label{lemma:Ledrappier-Strelcynlemma} Let $n\geq0$ and let $\varphi$ be a positive finite measurable function on $\widehat{\mathbb{P}^k}$. If $\LLog^{+}\left(\frac{\varphi\circ {\widehat{f}}^n}{\varphi}\right)\in L^1(\widehat{\mu})$ then : 
$$\LLog\left(\frac{\varphi\circ {\widehat{f}}^n}{\varphi}\right)\in L^1(\widehat{\mu})\ \mathrm{and}\ \int\LLog\left(\frac{\varphi\circ {\widehat{f}}^n}{\varphi}\right)\ \mathrm{d}\widehat{\mu} = 0.$$
\end{lemme}
The assumptions of this lemma are satisfied. Indeed, Equation (\ref{eq:logdn+2logDelta=-logq_x+logl(x_n)/L(x)}) implies for all $\widehat{x}$ in $\{h_n\geq 0\}\cap\Lambda$:
$$0\leq h_n(\widehat{x}) = \LLog(d^{rn}\Delta(x_0,n)^2) + \LLog\ q_{\widehat{x}}\left(\left({\widehat{f}}^{-n}\xi^u\right)_{\widehat{x}}\right) \leq \LLog(d^{rn}\Delta(x_0,n)^2).$$
Moreover, $\widehat{x}\mapsto\LLog\ \Delta(x_0,n)$ belongs to $L^1(\widehat{\mu})$ by \Dcitethm{\ref{thm:oseledec}}, and the function $\widehat{x}\mapsto L(\widehat{x})$ is measurable by Remark \ref{rmq:remarkonthedefinitionofthemeasuresq_x}. So, we can apply \Dcitelem{\ref{lemma:Ledrappier-Strelcynlemma}} to get $h_n\in L^1(\widehat{\mu})$ and $\int h_n\ \mathrm{d}\widehat{\mu}=0$ as desired, and therefore we obtain \eqref{eq:laformuledutheoremB} from \eqref{eq:logdn+2logDelta=-logq_x+logl(x_n)/L(x)} by integrating.\\

To complete the proof of the first five items of Theorem \ref{thm:TheoremB}, we now prove that $(\pi_0)_* q_{\widehat{x}} \ll (T^r \wedge \omega_{\mathbb{P}^k}^{k-r}) =: \sigma_{T^r}$ for every $\widehat{x} \in \Lambda$, and that $(\pi_0)_*\widehat{q}\ll T^r\wedge\omega_{\mathbb{P}^k}^{k-r}$. 

\begin{lemme}\label{lemma:pi_0_*q_x<<Twedgeomega} The following points hold:
\begin{enumerate}
    \item Let $N\geq 0$ and let $\widehat{x}\in {\widehat{f}}^N(\mathcal{A})\cap\Lambda$. Then $(\pi_0)_*q_{\widehat{x}}\ll \sigma_{T^r}$.
    \item For $\widehat{\mu}$-almost every $\widehat{x}$, $(\pi_0)_*q_{\widehat{x}}\ll \sigma_{T^r}$.
    \item We have also $(\pi_0)_*\widehat{q}\ll \sigma_{T^r}$.
\end{enumerate}
\end{lemme}
\noindent{\textbf{\underline{Proof:}}} Recall that we have assumed $\#\mathrm{Rec}_{\widehat{x}}(\mathcal{A})=+\infty$ for each $\widehat{x}\in\Lambda$, using Poincaré recurrence theorem (see Section \ref{sec:recproperties}). In particular, we have $\Lambda\subset\bigcup_{N\geq0}\widehat{f}^N(\mathcal{A})$, therefore, the first item of the lemma implies the second since $\widehat{\mu}\left(\Lambda\right)=1$. We thus prove only the first and third items.
To do so, let us remark that for every $\widehat{x}\in\Lambda$, the measure $T^r\wedge dd^c|W_{\widehat{x}}|^2$ is absolutely continuous with respect to the trace measure $\sigma_{T^r}=T^r\wedge\omega_{\mathbb{P}^k}^{k-r}$ on the domain of definition of $W_{\widehat{x}}$. So, for $\widehat{x}\in\mathcal{A}\cap\Lambda$, we can define the Radon-Nikodym derivative ($W_{\widehat{x}}$ is defined on $B(x_0,2\mathcal{R})$ by Proposition \ref{eq:pi(eta)includedinDom(f-n)}):
\begin{equation}\label{defn:RadonNikodymderivativesD_x}
    D_{\widehat{x}}:=\frac{d(T^r\wedge dd^c|W_{\widehat{x}}|^2)}{d\sigma_{T^r}}\ \mathrm{on}\ B(x_0,2\textcolor{black}{\mathcal{R}}).
\end{equation}

\noindent\textit{1.} Let us fix a Borel set $B$ of $\mathbb{P}^k$, and let $C(N,\widehat{x}):= d^{rN}\Delta(x_{-N},N)^2M(\widehat{x}_{-N})^2 / L(\widehat{x})$ be a constant so that
$$q_{\widehat{x}}(\pi_0^{-1}(B))  = C(N,\widehat{x})\times(T^r\wedge dd^c|W_{\widehat{x}_{-N}}|^2)\circ\pi_0 \circ {\widehat{f}}^{-N}[\xi^u_{\widehat{x}}\cap \pi_0^{-1}(B)].$$
Now observe that ${\widehat{f}}^{-N}[\xi^u_{\widehat{x}}\cap \pi_0^{-1}B]  = ({\widehat{f}}^{-N}\xi^u)_{\widehat{x}_{-N}} \cap ({\widehat{f}}^{-N}\pi_0^{-1}B) =  ({\widehat{f}}^{-N}\xi^u)_{\widehat{x}_{-N}} \cap \pi_0^{-1}(f^{-N}B)$. Recall that $f^N$ is injective on $\pi_0\left[({\widehat{f}}^{-N}\xi^u)_{\widehat{x}_{-N}}\right]$ by \Dcitethm{\ref{thm:propertiesofeta}}. Let $g_N$ denote the inverse of the restriction of $f^N$ to $\pi_0\left[({\widehat{f}}^{-N}\xi^u)_{\widehat{x}_{-N}}\right]$, and let $B':=B\cap \pi_0(\xi^u_{\widehat{x}})$. We obtain
$$q_{\widehat{x}}(\pi_0^{-1}(B))/C(N,\widehat{x})\leq(T^r\wedge dd^c|W_{\widehat{x}_{-N}}|^2)\left(g_N (B')\right) =:  l_{N,\widehat{x}}(B').$$
Using Definition \eqref{defn:RadonNikodymderivativesD_x}, we have
$$l_{N,\widehat{x}}(B')  = \int_{g_N(B')} D_{\widehat{x}_{-N}}\ \mathrm{d}\sigma_{T^r}.$$
Now it is sufficient to show that $\sigma_{T^r}(B)=0\Longrightarrow l_{N,\widehat{x}}(B')=0$. Let $E_N$ be the continuous function on $\mathbb{P}^k$ defined on $\mathbb{P}^k\backslash\mathrm{Crit}\left(f^N\right)$ by $E_N(p):=||(d_pf^N)^{-1}||^{-2}$ and extended by $0$ on $\mathrm{Crit}(f^N)$.
It satisfies $E_N \times\omega_{\mathbb{P}^k} \leq \left(f^N\right)^*\omega_{\mathbb{P}^k}$ and therefore 
$$\sigma_{T^r} \leq E_N^{-(k-r)} \times T^r\wedge (f^N)^*\omega_{\mathbb{P}^k}^{k-r} = E_N^{-(k-r)} \times d^{-rN} \left(f^N\right)^* (T^r\wedge \omega_{\mathbb{P}^k}^{k-r}) \textrm{ on } \mathbb{P}^k\backslash\mathrm{Crit}(f^{{N}}) ,$$ 
where we used $(f^N)^*T^r=d^{rN} T$ for the equality. Recall that $T$ has continuous local potentials, therefore by the Chern-Levine-Nirenberg inequality, the measures $\sigma_{T^r}$ and $(f^N)^* (T^r \wedge \omega_{\mathbb{P}^k}^{k-r})$ do not charge analytic subsets of $\mathbb{P}^k$ of codimension $\geq 1$. In particular, these measures do not charge $\mathrm{Crit}(f^N)$, so the preceding estimate actually holds on $\mathbb{P}^k$. Hence
$$l_{N,\widehat{x}}(B')\leq d^{-rN}\int_{g_N(B')} E_N^{-(k-r)} D_{\widehat{x}_{-N}}\ \mathrm{d}\left(f^N\right)^*\left[T^r \wedge \omega_{\mathbb{P}^k}^{k-r}\right].$$
Since $f^N:g_N(B')\longrightarrow B'$ is invertible with inverse map $g_N|_{B'}$, we get 
$$l_{N,\widehat{x}}(B')\leq d^{-rN}\int_{B'} (E_N^{-(k-r)}D_{\widehat{x}_{-N}}) \circ g_N \ \mathrm{d}\left[T^r \wedge \omega_{\mathbb{P}^k}^{k-r}\right].$$
That implies $\sigma_{T^r}(B)=0 \Longrightarrow l_{N,\widehat{x}}(B')=0$, as desired.\\

\noindent\textit{3.} By Proposition \ref{prop:q_x=q_ysieta_x=eta_y}, $(\pi_0)_*\widehat{q}(B)=\int q_{\widehat{x}}(\pi_0^{-1}(B)\cap\xi^u_{\widehat{x}})\ \mathrm{d}\widehat{\mu}(\widehat{x})=\int({\pi_0}_*q_{\widehat{x}})(B)\ \mathrm{d}\widehat{\mu}(\widehat{x})$ for any Borel set $B$, where the second equality comes from the fact that $q_{\widehat{x}}(\cdot)=q_{\widehat{x}}(\cdot\cap\xi^u_{\widehat{x}})$. The conclusion then follows since the measures $(\pi_0)_*q_{\widehat{x}}$ are absolutely continuous with respect to $\sigma_{T^r}$ by the previous item.\qed

\subsection{Proof of Theorem {\ref{thm:TheoremB}} Item \textit{6.} and proof of Theorem {\ref{thm:TheoremA}}}\label{sec:proofmainresult}

The third item of the following theorem implies in particular Theorem \ref{thm:TheoremA}. The probability measure $\widehat{q}$  is defined in \Dciteprop{\ref{prop:q_x=q_ysieta_x=eta_y}}. 

\begin{thm}\label{thm:q_x((g-neta)_x)=mu_x((g-neta)_x)} Assume $\lambda_r>\lambda_{r+1}=\cdots=\lambda_k=\lambda_{\mathrm{min}}=\frac{1}{2}\LLog\ d$. 
\begin{enumerate}
\item For $\widehat{\mu}-$almost every $\widehat{x}$ and for every $n\geq0$, $q_{\widehat{x}}\left({\widehat{f}}^{-n}\xi^u\right)_{\widehat{x}} = \mu_{\widehat{x}}\left({\widehat{f}}^{-n}\xi^u\right)_{\widehat{x}}$.
\item For $\widehat{\mu}-$almost every $\widehat{x}$, $q_{\widehat{x}}=\mu_{\widehat{x}}$ and so $\widehat{q} = \widehat{\mu}$. 
\item We have $\mu\ll T^r\wedge\omega_{\mathbb{P}^k}^{k-r}$.
\end{enumerate}
\end{thm}
\noindent\textbf{\underline{Proof:}} The second item is a classical consequence of the first one, using that the $\sigma-$algebra $\bigvee_{n\geq 0}\mathcal{M}({\widehat{f}}^{-n}\xi^u)$ coincide with the $\sigma-$algebra $\mathcal{M}$ (Item \textit{5}. of Theorem \ref{thm:propertiesofeta}), see for instance {\cite[Lemme 5.6]{dup06}}. The third item is deduced from the second, since by Lemma \ref{lemma:pi_0_*q_x<<Twedgeomega} we have $(\pi_0)_* \widehat{q} \ll T^r \wedge \omega_{\mathbb{P}^k}^{k-r}$ (and since $(\pi_0)_*\widehat{\mu}=\mu$). So, let us prove the first item. By Theorem \ref{thm:propertiesofeta} we know that for all $n\geq0$:
$$\int_{\widehat{\mathbb{P}^k}}-\LLog\ \mu_{\widehat{x}}\left({\widehat{f}}^{-n}\xi^u\right)_{\widehat{x}}\ \mathrm{d}\widehat{\mu}(\widehat{x}) = \LLog\ d^{rn} + 2n(k-r)\frac{1}{2}\mathrm{Log}\ d.$$
But according to \eqref{eq:laformuledutheoremB} we have also for all $n\geq0$:
$$\int_{\widehat{\mathbb{P}^k}}-\LLog\ q_{\widehat{x}}\left({\widehat{f}}^{-n}\xi^u\right)_{\widehat{x}}\ \mathrm{d}\widehat{\mu}(\widehat{x}) = \LLog\ d^{rn} + 2n(k-r)\lambda_{\mathrm{min}}.$$
Since $\lambda_{\mathrm{min}}=\frac{1}{2}\mathrm{Log}\ d$ we get $2n(k-r)\frac{1}{2}\mathrm{Log}\ d=2n(k-r)\lambda_{\mathrm{min}}$ and thus:
$$\int_{\widehat{\mathbb{P}^k}}-\LLog\ q_{\widehat{x}}\left({\widehat{f}}^{-n}\xi^u\right)_{\widehat{x}}\ \mathrm{d}\widehat{\mu}(\widehat{x}) = \int_{\widehat{\mathbb{P}^k}}-\LLog\ \mu_{\widehat{x}}\left({\widehat{f}}^{-n}\xi^u\right)_{\widehat{x}}\ \mathrm{d}\widehat{\mu}(\widehat{x}).$$
By Jensen inequality we deduce :
\begin{equation} \label{eq:jensen}
    0=\int_{\widehat{\mathbb{P}^k}}\LLog\ \frac{q_{\widehat{x}}}{\mu_{\widehat{x}}}\left({\widehat{f}}^{-n}\xi^u\right)_{\widehat{x}}\ \mathrm{d}\widehat{\mu}(\widehat{x})\leq \LLog\ \int_{\widehat{\mathbb{P}^k}} \frac{q_{\widehat{x}}}{\mu_{\widehat{x}}}\left({\widehat{f}}^{-n}\xi^u\right)_{\widehat{x}}\ \mathrm{d}\widehat{\mu}(\widehat{x}).
\end{equation}
By definition of the conditional measures of $\widehat{\mu}$ with respect to $\xi^u$, we get 
$$\int_{\widehat{\mathbb{P}^k}}\frac{q_{\widehat{x}}}{\mu_{\widehat{x}}}\left(\left({\widehat{f}}^{-n}\xi^u\right)_{\widehat{x}}\right)\ \mathrm{d}\widehat{\mu}(\widehat{x})  = \int_{\widehat{\mathbb{P}^k}}\int_{\xi^u_{\widehat{x}}}\frac{q_{\widehat{y}}}{\mu_{\widehat{y}}}\left(\left({\widehat{f}}^{-n}\xi^u\right)_{\widehat{y}}\right)\ \mathrm{d}\mu_{\widehat{x}}(\widehat{y})\ \mathrm{d}\widehat{\mu}(\widehat{x}).$$
Recall now that, since $\lambda_{r+1} = \cdots = \lambda_k$, the measures $q_{\widehat{x}}$ are the conditional measures of $\widehat{q}$, and thus $q_{\widehat{x}}$ depends only on $\xi^u_{\widehat{x}}$ and not on $\widehat{x}$. More precisely, according to Proposition \ref{prop:q_x=q_ysieta_x=eta_y}, for all $\widehat{x} \in \Lambda$ and for all $\widehat{y} \in \Lambda \cap \xi_{\widehat{x}}^u$, $q_{\widehat{x}} = q_{\widehat{y}}$. Since $\widehat{\mu}(\Lambda) = 1$, we must have $\mu_{\widehat{x}}(\Lambda) = 1$ for $\widehat{\mu}$-almost every $\widehat{x}$. Thus, we deduce that for $\widehat{\mu}$-almost every $\widehat{x}$:
$$\int_{\xi^u_{\widehat{x}}}\frac{q_{\widehat{y}}}{\mu_{\widehat{y}}}\left(\left({\widehat{f}}^{-n}\xi^u\right)_{\widehat{y}}\right)\ \mathrm{d}\mu_{\widehat{x}}(\widehat{y}) = \int_{\xi^u_{\widehat{x}}}\frac{q_{\widehat{x}}}{\mu_{\widehat{x}}}\left(\left({\widehat{f}}^{-n}\xi^u\right)_{\widehat{y}}\right)\ \mathrm{d}\mu_{\widehat{x}}(\widehat{y}). $$
Now by \Dcitethm{\ref{thm:propertiesofeta}}, $\xi^u_{\widehat{x}}$ is a countable union $\bigsqcup_{j\in\mathbb{N}}A_{j,n}^{\widehat{x}}$ of atoms of  ${\widehat{f}}^{-n}\xi^u$. We infer 
\begin{align*}
    \int_{\xi^u_{\widehat{x}}}\frac{q_{\widehat{x}}}{\mu_{\widehat{x}}}\left(\left({\widehat{f}}^{-n}\xi^u\right)_{\widehat{y}}\right)\ \mathrm{d}\mu_{\widehat{x}}(\widehat{y}) 
                                             & = \sum_{j=0}^{+\infty}\int_{A_{j,n}^{\widehat{x}}}\frac{q_{\widehat{x}}}{\mu_{\widehat{x}}}\left(A_{j,n}^{\widehat{x}}\right)\ \mathrm{d}\mu_{\widehat{x}}(\widehat{y}) = \sum_{j=0}^{+\infty}\frac{q_{\widehat{x}}}{\mu_{\widehat{x}}}\left(A_{j,n}^{\widehat{x}}\right)\times\mu_{\widehat{x}}\left(A^{\widehat{x}}_{j,n}\right)
\end{align*}
which is equal to $q_{\widehat{x}}\left(\bigsqcup_{j}A^{\widehat{x}}_{j,n}\right)=1$. Finally, 
$\int_{\widehat{\mathbb{P}^k}}\frac{q_{\widehat{x}}}{\mu_{\widehat{x}}} \left({\widehat{f}}^{-n}\xi^u\right)_{\widehat{x}}\ \mathrm{d}\widehat{\mu}(\widehat{x}) = 1$, which implies the equality in \Dciteequa{\ref{eq:jensen}}. The strict concavity of$\LLog$ completes the proof.\qed

\bibliographystyle{abbrv}
\begin{otherlanguage}{english}
    \bibliography{biblio}
\end{otherlanguage}

$ $ \\
\noindent {\footnotesize V. Tapiero}\\
{\footnotesize Universit\'e de Rennes}\\
{\footnotesize CNRS, IRMAR - UMR 6625}\\
{\footnotesize F-35000 Rennes, France}\\
{\footnotesize virgile.tapiero@univ-rennes1.fr}\\

\end{document}